\documentstyle[amstex]{article}
\numberwithin{equation}{section}
\newtheorem{thm}{Theorem}[section]
\newtheorem{lem}{Lemma}[section]
\newtheorem{defn}{Definition}[section]
\newtheorem{prop}{Proposition}[section]
\newtheorem{sublem}{Sublemma}[section]
\newtheorem{cor}{Corollary}[section]
\begin{document}
\Large

\title{ \Large \bf Infinite sequences of almost Kaehler manifolds with \\ 
high symmetry, their perturbations  
and \\ 
pseudo holomorphic curves }  
\author{ \Large \bf Tsuyoshi Kato}

\date{}
\maketitle
\begin{abstract}
We study  analysis over infinite dimensional manifolds
consisted by sequences of almost Kaehler manifolds.
We develop moduli theory of pseudo holomorphic curves 
into such spaces with high symmetry.
Many mechanisms of the standard moduli theory over finite dimensional
spaces also work over these infinite dimensional spaces, 
which is based on a  simple functional analytic framework.
\end{abstract}
\begin{center}{ \bf  Introduction}\end{center}
Global  analysis of 
infinite dimensional spaces  is one of the
central theme to study in geometry.
There would immediately occur several difficulties when we try to follow and apply 
any established methods over finite dimensional spaces,
one of whose main reasons  come from  locally non compactness. 
Still the  extensive developments have been achieved
from the view points of the metric spaces, which 
essentially measure one dimensional objects. 

In this paper, we  study 
 infinite dimensional  geometry and analysis from 
 two dimensional view points, and 
 introduce  a class of 
 infinite dimensional geometric spaces 
consisted by sequences of embeddings by finite dimensional manifolds.
We develop a foundation of analytic tool to perform some functional analysis
under the conditions of  high symmetry over such spaces.

 Moduli theory of pseudo holomorphic curves
turned out to be   the very powerful tool and has become one of the central theme
in symplectic geometry. From the  view point of infinite dimensional geometry,
it has led us to study   
{\em  almost Kaehler sequences}  $[(M_i, \omega_i, J_i)] $  which consist of
 families of   embeddings by  almost Kaehler manifolds:
$$[(M_i, \omega_i, J_i)] =  (M_0 , \omega_0, J_0) \subset (M_1 , \omega_1, J_1) 
                \subset \dots \subset (M_i, \omega_i, J_i) \subset \dots $$
                In order to develop global analysis over these spaces,
                we introduce several set-ups of function spaces.
                Based on the analytic fields, we study 
                maps into these infinite dimensional spaces. In particular we 
                 develop the moduli theory of pseudo holomorphic curves into such  sequences
                 from two dimensional spheres.
  Our  construction provides with  two main ingredients.
One is Fredholm theory for the linearlized maps,
where it requires two conditions, 
 closedness of the maps and  finite dimensionality.
The other is non linear analysis where it also requies two conditions,
 regularity of holomorphic maps and
compactness of the moduli spaces. 
We discover that these properties also hold over 
the   infinite dimensional spaces
under the conditions of high symmetry.

As an application, these constructions
 allow us to calculate or estimate the capacity invariants 
 over such infinite dimensional spaces.
Capacity invariants have been
introduced by 
 Hofer and Zehnder   with the axiom of the invariants over  finite dimensional
 symplectic manifolds  ([HZ]). 
  Kuksin  extended the construction of the  capacity invariants over
   symplectic Hilbert spaces, and
 investigated  invariance under  the flow maps ([Ku], [B]).
It turned out that they play an important role in 
 relation with Hamiltonian PDE,
whose phase spaces are infinite dimensional.

This paper is organized by  three constructions:

(A) Formulation of a class of  infinite dimensional spaces
which allow us to analyze passing though infinite dimensional 
local charts. Infinite dimensional spaces behave quite flexibly
under embeddings from each other.
It allows us to introduce the {\em infinitesimal neighbourhoods} 
of such spaces, which are consisted by 
`convergent sequences' of infinite dimensional spaces.
These notions would provide us with rich classes of spaces
which are able to apply our  functional analytic methods.

(B)
Construction of the moduli theory
into such infinite dimensional spaces.

(C)
 Analysis of the capacity invariants over such spaces and their calculations.
Combining with these constructions, we study stability of the capacity invariants
under  deformation in the infinitesimal neighourhoods.
\vspace{3mm} \\
Let us explain more details of the contents.

We carefully  introduce  the differentiable functions over the infinite dimensional spaces,
which should be able to be extended over their completion to the Hilbert manifolds.

We  will verify  that 
many mechanisms of the moduli theory of holomorphic curves work over
almost Kaehler sequences with  some symmetric conditions.
Many of the  infinite complex homogeneous spaces satisfy such symmetries.


    $[(M_i, \omega_i, J_i)]$
is said to be a
  {\em symmetric almost Kaehler sequence}, if for any $k \geq 0$,
   there are some $l =l(k)$, 
  families of almost Kaehler submanifolds $M_k \subset W_i \subset M_i$ 
 for all $i \geq l$  with $W_l = M_l$ 
  and
  isomorphisms: 
$$P_i: \{(M, M_k), \omega, J\}  \cong \{(M,M_k), \omega, J\}$$
with $M \equiv \cup_j M_j$, 
which preserve $M_k$  and transform $W_i$ to $M_l$ as:
$$P_i: (W_i, \omega_i|W_i, J_i|W_i) \cong (M_l,\omega_l, J_l)$$
so that  at any $p \in M_k$, 
$D :   TM_k \oplus_{i \geq l}  N_{i,k}  \cong T M |M_k $
 give the uniform isomorphisms over $M_k$
with respect to the complete local charts, 
where $N_{i,k}$ are isomorphic to  the normal bundles of $TM_k \subset TM_l|M_k$:
$$ N_{i,k}= (P_i^{-1})_* [  \ (\text{ Ker } ( \pi_k)_* \cap TM_l) |M_k) \ ] , \quad
 D_p = \text{ id } \oplus (P_l)_* \oplus (P_{l+1})_* \oplus  \dots $$
(Precisely see def $1.4$.)
Symmetric Almost Kaehler sequences 
satisfy the uniform isomorpshims which is
 one  of the main property 
we need:
$$TM|M_k \cong TM_k \oplus (N_{k,l} \otimes {\bf R}^{\infty}).$$

Our infinite dimensional geometry relies  on the following functional analytic property:

\begin{lem}
Let $F: H \to H$ be a bounded operator with closed range and finite dimensional kernel.
Then for any closed linear subspace $L \subset H$ and the Hilbert space tensor product 
$L \otimes W$ with another Hilbert space $W$, 
the image of the induced operators: 
$$(F\otimes 1)(L \otimes W) \subset H \otimes W$$
also have closed range.
\end{lem}

Now let us introduce several functional spaces
 in order  to develop moduli theory of holomorphic curves.
 Later on we fix two points $p_0, p_{\infty} \in M_0 \subset M$.
 
 Let $E(J)_i, \ F_i \mapsto S^2 \times M_i$ be 
 vector bundles whose fibers are respectively:
\begin{equation} 
\begin{align*}
& E(J)_i(z,m) = \{ \phi: T_z S^2 \mapsto T_mM_i : 
   \text{ anti complex linear } \}, \\
& F_i(z,m) = \{ \phi: T_z S^2 \mapsto T_mM_i : \text{ linear} \}.
\end{align*}
\end{equation}

For fixed large $l \geq 1$, let $L^2_{l+1}(S^2, M_i)$ be the sets of
$L^2_{l+1}$ maps from $S^2$ to $M_i$, and 
define  the spaces of Sobolev maps:
\begin{equation}
\begin{align*}
 {\frak B}_i \equiv  {\frak B}_i (\alpha) & = \{  u   \in   L^2_{l+1}(S^2, M_i) :  
 [u] = \alpha ,   \\
 & \int_{D(1)} u^*(\omega) = \frac{1}{2}<\omega, \alpha>, \quad
  u(*) = p_* \in M_0  , * \in \{ 0, \infty \}  \}.
  \end{align*}
  \end{equation}
Then we have   two stratified Hilbert bundles over ${\frak B}_i$:
 \begin{equation}
 \begin{align*}
&  {\frak E}_i = L^2_l({\frak B}_i^*(E(J)_i)) = 
\cup_{u \in {\frak B}_i} \{u\} \times L^2_l(u^*(E(J)_i)), \\
&   {\frak F}_i = L^2_l({\frak B}_i^*(F_i)) 
= \cup_{u \in {\frak B}_i} \{u\} \times L^2_l(u^*(F_i)).
\end{align*}
\end{equation}
These spaces admit continuous $S^1$ actions.

 The non linear Cauchy-Riemann operator is  given  as  sections:
$$ \bar{\partial}_J \in C^{\infty}({\frak E}_i \mapsto {\frak B}_i),  
\quad
  \bar{\partial}_J (u) = Tu + J \circ Tu \circ i .$$
$u$ is called a {\em holomorphic curve} if it satisfies the equation 
$\bar{\partial}_J(u)=0$.
The moduli space of holomorphic curves is defined by:
$$     {\frak M}(\alpha, M_i,J_i)  = 
    \{ u \in C^{\infty}(S^2, M_i) \cap {\frak B}_i(\alpha) 
      : \bar{\partial}_{J} (u) =0   \}.$$
$J$ is called regular,
 if the linealizations
  are onto for all $u \in    {\frak M}(\alpha, M_i,J_i) $ and all $i \geq 0$:
 $$D \bar{\partial}_{J}(u): T_u {\frak B}_i \mapsto ({\frak E}_i)_u$$

Let $u \in {\frak B} \equiv \cup_i  {\frak B}_i$ and take an open neighbourhood
$U(u) \subset {\frak B}$. By introducing the Sobolev norms, one can make completion
$U(u) \subset \hat{U}(u)$. Notice that elements in $\hat{U}(u)$ cannot be realized 
by  maps into $M =\cup_i M_i$ in general.
Let: 
 $$ \bar{\partial}_{J}: \hat{U}(u) \mapsto \hat{\frak E}|\hat{U}(u)$$
 be the extension of the CR operator over the completion.
 The differential of the operator  is not necessarily onto even if it is regular, where
 the range may not be closed.
 An almost Kaehler sequence is said to be {\em strongly regular},
  if  the extensions are onto  at all:
 $$u \in   {\frak M}([(M_i, \omega_i, J_i)])   \equiv \cup_i \   {\frak M}(\alpha, M_i,J_i) .$$
 This is the key property of moduli theory we develop in the infinite dimensional setting.
 
 \begin{thm}
 Let $[(M_i, \omega_i, J_i)]$ be a symmetric  Kaehler sequence.
 
 (1)
Suppose it is  regular and dim $\cup_i $ Ker $ D_u \bar{\partial}_i =N$ is finite,
 then it is in fact strongly regular of index $N$.

 In particluar  
  $ {\frak M}([(M_i, \omega_i , J_i)]) $ is a regular $N$ dimensional $S^1$ free manifold.
  
  (2) Assume  moreover  it is  isotropic, and each connected component of
   $  {\frak M}([(M_i, \omega_i , J_i)])$ is bounded. Then 
   the equality holds:
   $$ {\frak M}([(M_i, \omega_i , J_i)]) =  {\frak M}( M_0,\omega_0, J_0) .$$
   
   In particular if it is minimal, then  $ {\frak M}([(M_i, \omega_i , J_i)])    $  is compact.
  \end{thm}
 The second condition in $ (2) $ is satisfied when $N=1$ ($2.D$).

\vspace{3mm}

Stability of geometric structures 
is one of the central  theme in geometry and analysis of dynamics.
Let us say that a  family of almost Kaehler sequences 
 {\em converges }: 
$$\{ [(M_i^l, \omega_i^l, J_i^l) ] \}_{l\geq 1}  \to [(M_i, \omega_i, J_i)]$$
as $l \to \infty$,
if there are positive $\epsilon >0$, 
  subindices $\{ k(i)\}_i$  and 
 $C^{\infty}$ and compatible embeddings:
 $$ I(i,l): M_i^l \hookrightarrow M_{k(i)}$$
 such that:
 
 (1) 
  $\{ I(i,l)(M_i^l) \}_{l \geq 1} \subset M_{k(i)}$ 
converges to $M_i$ in $M_{k(i)}$ in $C^{\infty}$ for each $i$,
 
 (2)
  there are uniform extensions of 
    almost complex $J(i,l)$ and  symplectic $\omega(i,l)$
 structures  over  the  open 
$\epsilon$  tubular  neighbourhoods: 
$$I(i,l) (M_i^l) \subset U(i,l) \subset M \equiv \cup_i M_i $$

(3)
 there are   compatible   families of holomorphic  maps:
$$\pi_k^l :  ( U(k,l) , J(k,l))  \mapsto (M_k^l, J_k^l).$$
(Precisely see def $3.1$.)
Analysis of perturbations of moduli theory is to study asymptotic
behaviour of structure of moduli spaces as $l \to \infty$.
In light of this aspect, we introduce the {\em infinitesimal neighbourhoods}
of $[(M_i, \omega_i, J_i)]$:
$${\frak N}( [(M_i, \omega_i, J_i) ] )$$
which consist of the equivalent classes
of the isomorphism classes of convergent families of almost Kaehler sequences,
where:
$$\{ [(M_i^l, \omega_i^l, J_i^l) ] \}_{l\geq 1} \sim 
\{ [(M_i^{l(h)}, \omega_i^{l(h)}, J_i^{l(h)}) ] \}_{  h \geq 1} $$
for all infinite subindices.

We say that the infinitesimal neighbourhood  is strongly regular, 
if  for any element $[ [(M_i^l, \omega_i^l, J_i^l) ] ]  \in {\frak N}( [(M_i, \omega_i, J_i) ] )$,
 there is a large $l_0$ so that for all 
$l \geq l_0$, 
$[(M_i^l, \omega_i^l, J_i^l) ] $ are all strongly regular.

In general regularity condition  in moduli theory
is not stable under convergence of
almost Kaehler sequences.
We verify that strong regularity condition overcomes  this difficulty and
controls analytic behaviour of  holomorphic maps.

Suppose  $[(M_i, \omega_i, J_i)]$ satisfies all  the conditions in theorem $0.1$:
  (1)  regular,  (2) minimal, (3) isotropically  symmetric 
  and   (4)   Kaehler sequence.
 \\
 (5) 
 Moreover assume  the moduli space ${\frak M}([(M_i, \omega_i, J_i)]) $ is $1$-dimensional.

Then by theorem $0.1$,  the followings hold:
\vspace{3mm}

(A) $ {\frak M}([(M_i, \omega_i, J_i)]) = {\frak M}(M_0, \omega_0, J_0)$, and both are compact.

(B)  $[(M_i, \omega_i, J_i)]$ is strongly regular.

\vspace{3mm}

Let us   introduce another geometric condition
which can be applied to analysis over infinitesimal neighbourhoods.
$[(M_i, \omega_i, J_i)]$ is {\em quasi transitive}, if
 for any  $N  >0$,
 there is  $k = k(N) $ so that automorphisms 
 $A: (M, \omega , J) \cong (M, \omega, J)$ exist 
with $A(p_i) \in M_k$  for   any 
$p_0, \dots, p_{N-1} \in M \equiv \cup_i M_i$ and 
all $0 \leq i \leq N-1$.

\begin{thm} Assume
moreover that  (6) $[(M_i, \omega_i, J_i)]$ is  quasi transitive. Then
the followings hold:
\vspace{2mm}

(C)  ${\frak N}^m([(M_i, \omega_i, J_i)])$ is strongly regular.
\vspace{3mm}

(D)There are homeomorphisms: 
$$ {\frak M}([(M_i, \omega_i, J_i)]) \cong {\frak M}([(M^l_i, \omega^l_i, J^l_i)]) $$
   for any $[[(M_i^l, \omega_i^l, J_i^l)] ]  
\in {\frak N}^m([(M_i, \omega_i, J_i)])$ and all large $l$.
\end{thm}
Here $m$ stands for minimality of maps in $\pi_2$.

In particular $(D)$ implies that roughly speaking 
there are no `divergent sequences'  $u^l \in 
 {\frak M}([(M^l_i, \omega^l_i, J^l_i)]) $ 
which approach to be  holomorphic  with respect to $J$, but
they do not converge to ant elements in 
$ {\frak M}([(M_i, \omega_i, J_i)]) $.

\vspace{3mm}

We apply the analysis of moduli theory over almost Kaehler sequences 
 to Hamiltonian dynamics. 
 
 Let $(M, \omega)$ be a finite dimensional symplectic manifold.
  A  Hamiltonian function  $f: M \mapsto [0, \infty)$ gives 
the Hamiltonian vector field $X_f$  by the relation
$df( \quad) = \omega(X_f, \quad)$.
A pre-admissible function  is called {\em admissible} if any non trivial periodic solutions 
have their periods larger than $1$.

Suprimum of the widths $\sup f - \inf f$ among all the admissible functions
is  called the capacity invariant of $(M, \omega)$, which 
contains deep information in Hamiltonian dynamics.

The capacity invariant cap$([(M_i, \omega_i, J_i)])$ over almost Kaehler
sequences are straightforwardly defined by use of bounded Hamiltonians
$f: M \to [0, \infty)$
over $M= \cup_i M_i$.

In order to study the capacity invariant  over
  perturbations of almost Kaehler sequences,
it turns out 
 that the  {\em asymptotic periodic solutions} of
Hamiltonian arise quite naturally, 
which consist of a family of loops 
$x_i : [0, T_i] \to M_i$ with:
  $$\sup_t |\dot{x_i} - X_{f_i}(x_i)| (t) \to 0, \quad i \to \infty.$$

Instead of periodic solutions, one can use asymptotic periodic solutions 
over Hamiltonians and define admissibility in a parallel way.
Then one obtains an analogue of capacity by use of such objects, 
which we call the {\em asymptotic capacity}:
$$\text{As-cap} ([(M_i, \omega_i, J_i)]).$$

Conversely if we regard asymptotic periodic solutions as though 
`periodic solutions' over  elements of an infinitesimal neighbourhood, 
then it leads us to formulate the capacity invariant
over the infinitesimal neighbourhoods:
\begin{align*}
 C^m([(M_i, \omega, J_i)]) =    
   \sup \{  \   \lim \sup_l  &  \text{ cap}([(M_i^l, \omega_i^l, J_i^l)]) : \\
&             [[(M_i^l, \omega_i^l, J_i^l)]] \in  {\frak N}^m([(M_i, \omega_i, J_i)]) \}
            \end{align*}
Notice that a priori estimates hold:
 $$C^m  \ \   \geq  \  \  \text{cap}  \  \ \geq    \ \  \text{As-cap}  \ \  \geq  \ \ 0.$$

\vspace{1cm}

 \begin{thm} Let $[(M_i, \omega, J_i)]$ be a minimal, 
  isotropically symmetric and quasi transitive  Kaehler sequence,
 with a  fixed minimal element 
 $ \alpha  \in \pi_2(M)$.
 
  If the 
moduli space of holomorphic curves is non empty, regular,
 $1$ dimensional
and $S^1$ freely cobordant to non zero,  then  the estimate:
$$C^m([(M_i,  \omega_i, J_i)]) \leq m$$
holds, where $m =<\omega, \alpha>$.
\end{thm}

Here we have the results of concrete calculations:

\begin{prop} 
(1) Let $ \text{Cap} $ be   $ \text{As-cap} $ or $ \text{cap}$. Then:
 $$C^m([{\bf CP}^i])= \text{Cap}([{\bf CP}^i])= \text{Cap}([D^{2i}]) = \pi.$$

(2) $C([(T^{4i}, \omega, J)]) = \infty$. 
\end{prop}
\quad
\vspace{3mm} \\
{\em Furthur research directions:}
Let us describe  possible developments for future, 
which  we partly study in this paper.

One of the important development is study of the
displacement energies, since our techniques in this paper
can be applied for $X \times  \lambda {\bf CP}^1$ if $X$ is the case with the additional condition of
projective or $\pi_2$ rank $1$,
where $\lambda {\bf CP}^1$ is ${\bf CP}^1$ equipped with the rescaled Fubini-Study form 
by $\lambda >0$ (see [LM]).
Let $M= \cup_i M_i$ and 
$\bar{M}$ be its completion to the Hilbert manifold.
Let us consider 
  a bounded Hamiltonian $H: M \to  [0, \infty)$ and
$D_H: \bar{M}\cong \bar{M}$  be the
 Hamiltonian vector field ($5.C$).

Let us say that an open subset $U \subset {M}$
is {\em displaceable}, if there is a bounded Hamiltonian $H$
as above so that
$D_H(\bar{U}) \cap \bar{U} = \phi $
hold in $ \bar{ M}$, where $\bar{U}$ is the closure of $U$ in $\bar{M}$.
We say that $H$ displaces $U$, and 
we  expect the energy estimates 
$\sup_{m \in M}  H(m)  \geq  \text{As-cap}(U)$.
For $X \subset {\bf CP}^{\infty}$, let $U(X) \subset {\bf CP}^{\infty}$ be the maximal 
open neighbourhoods of $X$ such that $U(X)$ are displaceable.
Estimates of capacity values of  $U(X)$ 
would be of particular interests for us,
 which might measure `symplectic  complexity' 
of $X$. 
 Below we introduce two categories:

{\em (A) Projective varieties:}
(Finite dimensional) projective varieties admit embeddings into ${\bf CP}^n$
for some $n$. 
Let us allow  to address two types of questions which are of interests for us.

 Let us consider  the moduli space of, say $K3$ surface, and take
 a non projective (but always Kaehler) variety $X$ which admit  projective
$K3$ surfaces $X_i \subset {\bf CP}^{n_i}$ converging to $X$ in the moduli spaces.
One can find some $X$ such that 
these  projective varieties must satisfy  $n_i \to \infty$.
Then we ask what are the asymptotic behaviour of values of $\{$As-cap$(U(X_i))\}_i$.

Let $X$ be a projective variety (possibly infinite dimension),
and ${\frak F}$ be the class of projective embeddings.
There are canonical embeddings such as Veronese 
and Pl\"ucker embeddings.
One may wander whether  positivity
$  \text{As-cap}(U(X)) \  >0$ might hold 
for both 
$X= Pl_r(Gr_r) , v_m({\bf CP}^{\infty}) \subset {\bf CP}^{\infty}$.

Let us fix $0 < \delta <1$ and consider $\delta$ balls:
$$B^{2n}(\delta) = \{ (z_1, \dots,z_n) : \Sigma_{i=1}^n |z_i|^2 < \delta^2\}
\subset { \bf C}^n.$$
equipped with the standard symplectic structure.
They are embedded into ${\bf CP}^n$ by:
$$(z_1, \dots,z_n) \in B^{2n}(\delta) \hookrightarrow [   \sqrt{1- |z|^2}  , z_1, \dots,z_n]$$
which preserve the symplectic structure compatibly with the Fubini-Study forms.
It induces the infinite dimensional embedding $I_2: B^{2\infty}(\delta) \subset {\bf CP}^{\infty}$.
Let 
$I_1: B^{2\infty}(\delta) \hookrightarrow B^{2\infty}(\delta)$
 be 
$(z_1, \dots,z_n, \dots) \to (0,z_1, \dots,z_n, \dots)$, and
 consider the composition:
$$S \equiv I_2 \circ I_1 :  B^{2\infty}(\delta)  \hookrightarrow {\bf CP}^{\infty}$$

\begin{lem} For $0 < \delta <1$, there are displaceable $U(S(B^{2\infty}(\delta))$
so that  positivity holds:
$$\text{As-cap }(U(S(B^{2\infty}(\delta)))) >0.$$
\end{lem}
{\em Proof:}
Let us choose $0 < \delta < \delta' <1$, and 
consider the embeddings 
$B^{2\infty}(\delta) \hookrightarrow B^{\infty}(\delta')$
 by
$(z_1, \dots,z_n, \dots) \to (0,z_1, \dots,z_n, \dots)$.

Let $f:[0, \infty) \to [0,1]$ be a smooth function 
with $f|[0, \delta^2]) \equiv \alpha >0$ and $f|[(\delta')^2, \infty) \equiv 0$,
where $\alpha >0$ is sufficiently small.
Let us put $z_0 =x_0+ \sqrt{-1}y_0$ and $\bar{z}= (z_1,z_2, \dots)$.
Then we have the bounded Hamiltonian:
$$F: B^{2\infty}(\delta') \to {\bf R}, \quad
F(z_0, \bar{z}) = f(|\bar{z}|^2)x_0.$$
Notice that it extends over ${\bf CP}^{\infty}$ by letting $0$ 
outside $B^{2\infty}(\delta') $.
The Hamiltonian vector field is given by the form
$X_F = (-\alpha, 0, *)$
on small neighbourhood of $S(B^{2\infty}(\delta) )$.
In particular the first coordinate of $D_F(S(B^{2\infty}(\delta)) )$ must be uniformly away from $0$.
This completes the proof.

\vspace{3mm}
Let us list two important embeddings into infinite dimensional projective spaces:

{\em Bergman (pseudo) metrics:}
Let $X$ be a finite dimensional complex manifold, and $W$ be the space of $L^2$-holomorphic $n$ forms over $X$.  The canonical inner product over $W$  gives  it a Hilbert space structure,
which gives rise to the Bergman kernel form. When it is pointwisely positive, then 
the map
$I: X \mapsto \bar{\bf CP}(W^*)$
is induced. When it is immersion then the Bergman metric is equipped over $X$ by pull-back of
the Fubini-Study.

{\em  Rational dynamics:}
In [K3], we have studied analytic behaviour of rational dynamics 
from the view point of scaling limits.
The simplest case of infinite algebraic varieties will be  the shift:
$$S : {\bf CP}^{\infty} \hookrightarrow {\bf CP}^{\infty}, \quad
[z_0,z_1, \dots ] \to [0,z_0,z_1, \dots]$$
whose image is  given by the homogeneous polynomial $F(z_0,z_1, \dots) =z_0$.

Let $f(z_0, \dots,z_{n-1})$ be a rational function and consider the iteration dynamics
given by the rule
$z_N = f(z_{N-n}, \dots, z_{N-1})$
with the initial values $(z_0, \dots,z_{n-1}) \in {\bf C}^n$.
Let us say that the orbits $\{z_N\}_N$ are affine, if 
$(z_0,z_1, \dots) \in {\bf C}^{\infty} \subset {\bf CP}^{\infty}$ 
by the embedding $(z_0,z_1, \dots) \to [1,z_0,z_1, \dots]$, and they 
 are rational  if $ [1,z_0,z_1, \dots]$ detemine points in ${\bf CP}^{\infty}$.
Now we have $X$ as all  the set of rational   orbits.

$f$ gives the {\em recursive} dynamics, if 
there exists some $M$ so that for any initial values, the corresponding orbits satisfy 
the equalities
$z_{N+M}=z_N$ hold for all $N \geq 0$.
When $f$ is recursive, non affine orbits are parametrized by 
algebraic sets in ${\bf C}^n$.

{\em Example:}
Let us consider
$f(z_0,z_1) = z_0^{-1} ( 1+z_1) $.
It is well-known that this gives the recursive dynamics of period $5$.
Notice that any affine orbits $\{z_N\}_N$ must satisfy $z_i \ne 0$ for all $i$,
which is equivalent to $z_i \ne 0$ for $0 \leq i \leq 4$.
By  straight forward computations, the non affine  set is given by
$V_s =\{(z_0, z_1):  z_0z_1=0\}$.

{\em (B) Complex vector bundles:}
Let  $[{\bf Gr}_{n,i}]$ be the complex Grassmannians and put
${\bf Gr}_n =  \cup_i  {\bf Gr}_{n,i}$.
Let $M$ be a complex manifold of dimension $0 \leq n <  \infty$, and 
 consider
the classifying maps $f: M \to {\bf Gr}_n$ for $TM$ with $X =f(M) \subset {\bf Gr}_n$.
 By modifying our functional spaces 
by change of dimension of fixing points as in defining the quantum cohomology,
one will be able to apply our techniques in this paper 
to the infinite Grassmannians.

\section{Almost Kaehler sequences}
{\bf 1.A  Function spaces over local charts:}
For positive $\epsilon >0$, 
let $D^{2k}( \epsilon)  \subset {\bf R}^{2k}$ be 
$\epsilon$ ball with the center $0$. We denote the $2i$ dimensional $\epsilon$ cube
$D(i) = D^2(\epsilon) \times \dots D^2(\epsilon)$ by multiplications of $D^2(\epsilon)$ by  $i$ times.
There are canonial embeddings 
$D_i =D(i) \times \{0\}  \subset D_{i+1}$ for all $i \geq 1$.
Let us put the infinite dimensional cube and disk by:
$$ D(\infty)   \equiv  \cup_{i \geq 1} \ D_i, \quad D(\epsilon ) \equiv  \cup_{k \geq 1} \ D^{2k}(\epsilon) 
 \quad  \subset \  {\bf R}^{\infty}$$
 respectively. Notice   diam $D(\infty)= \infty$.
\vspace{3mm} \\
  {\bf 1.A.1 Hilbert completion:}
 Let $H$ be the separable Hilbert space which is obtained by the 
 completion of ${\bf R}^{\infty}$ with the standard metric on it.
For $p= (p_0,p_1, \dots)  \in  {\bf R}^{\infty}$, 
let us denote by $D(\infty)(p)  \equiv D(\infty)+p$
and $D(\epsilon)(p) \equiv D(\epsilon) +p$  as
  the infinite dimensional cube and disk  with the center $p$ respectively.
We denote  the metric  completion by:
 $$\bar{D}(\infty)(p) , \quad \bar{D}(\epsilon)(p) \quad   \subset  \quad H.$$
 
  By a  neighbourhood or an 
open subset of $p \in B \subset {\bf R}^{\infty}$,
we mean that $B$ contains some $\delta >0$ disk with the center $p $,
with respect to the metric.

Let $ p \in B \subset {\bf R}^{\infty}$ be an open subset, and denote its closure by
 $\bar{B} \subset H$.
Let us consider a smooth and bounded function 
$f: B \mapsto {\bf R}$. We will regard the derivatives of $f$ at $p$
as the linear operators:
$$\nabla f:  T_p {\bf R}^{\infty} \equiv  \cup_{k \geq 1} T_p {\bf R}^{2k}
 \mapsto {\bf R}, \quad 
\nabla^2 f  : (T_p {\bf R}^{\infty})^{\otimes 2} \mapsto {\bf R},  \text{ etc}.$$
where
$\nabla^2 (f) (v,w) = \frac{\partial^2}{\partial s \partial t} f(p+sv+tw)|_{s=t=0}$
just for clarity.

Let us denote the operator norm by  $|\nabla^l f|(p)$, if
 it extends to the bounded  linear functional:
$$\nabla^l f : (T_p H)^{\otimes l} \mapsto {\bf R}.$$

\begin{defn}
We will say that $f$ is of {\em  completely  $C^k$ bounded geometry}
at $p=(p_0, p_1, \dots) \in B$, if:

\vspace{3mm} 

(1) 
 $f|B$ extends to a continuous function  on some neighbourhood of 
 $p $ in $\bar{B}(p)$ and 

(2) each  differential $\nabla^l f : T_pH \mapsto {\bf R}$
exists continuously on    some neighbourhood of 
 $p $ in $\bar{B}(p)$ for all $0 \leq l \leq k$
  (so  $|\nabla^l f|(p) < \infty$ hold) 
\end{defn}
 We say $f$ is of  completely  $C^k$ bounded geometry, if it is at any point
$p \in B$  satisfying uniformity:
$$||f||^2 C^k(B)  \equiv 
\text{Sup}_{p \in B} \Sigma_{0 \leq l \leq k} |\nabla^l f|^2(p) < \infty.$$  
Just completely bounded geometry implies $C^{\infty}$.

A pointwise operator $D$ on functions over $B$
is of completely bounded geometry, 
if $D$ extends to a smooth  operator  $D(p)$
over  functions of completely bounded geometry over $\bar{B}(p)$   at each $p \in B$.
Namely the followings are satisfied:

\vspace{3mm}

(1)  There is a constant $C$ with:
$$D: C^0(B) \mapsto C^0(B), \quad |Df|(p) \leq C | f|(p), \quad
p \in B.$$ 
We denote its extended and pointwise operator norm as $|D(p)|$.

(2) For all $k$, the following norms are all finite:
$$||D||^2 C^k(B)   \equiv 
\text{Sup}_{p \in B} \Sigma_{0 \leq l \leq k}
 |\nabla^l D|^2 (p) < \infty.$$ 
 
 $D$ gives a {\em complete isomorphism}, if it is of completely of bounded geometry.
Moreover 
  there are constants $0< c<c'$
 so that  the uniform bounds hold  for each $p \in B$:
 $$c  \leq ||D_p||  \leq c'.$$

For pointwisely bilinear forms, one has the parallel notion
of {\em complete nondegeneracy}.

Later we will always treat almost Kaehler  sequences whose almost complex structures,
symplectic structures or the induced Riemannian metrics are all
 completely nondegenerate  ($1.B$).
\vspace{3mm} \\
{\em Example 1.1:} Let $D^2  \subset {\bf R}^2$  be
 the standard  ball  with the center $0$, and
   consider   smooth functions 
$g, h : D^2 \mapsto [0,1]$
where:
\begin{equation}
\begin{align*} 
 & g(x) = exp(-\frac{|  x|^2}{1-|  x|^2}), 
 \quad  h(x) = exp(-\frac{|  x|}{1-|  x|}).
\end{align*}
\end{equation}
Let us prepare infinite numbers of 
the same $g$ and $h$, and in order to distinguish these, 
let us assign indices  as $g_i, h_i: D^2_i \mapsto [0,1]$.
Then we consider functions over 
$D(\infty) = D^2_0 \times D^2_1 \times  \dots $:
$$ G=g_0g_1g_2 \dots, \quad H= h_0h_1 h_2 \dots$$ 
by the pointwise multiplication.
Both $G$ and $H$ are  smooth on $D(\infty)$.
They have the following properties:

\vspace{3mm}

(1)  $G$ is of completely   bounded geometry on the unit balls with the center zero, and
 (2) $H$ is not at any points.
\vspace{3mm} \\
$H$ is not even continuous on $\bar{D}(\infty)$.
 In fact let us choose families of points  
$\{p(l)  =(p_0(l) ,  p_1(l) , \dots) \} \subset \bar{D}(\infty) \backslash D(\infty)$ with:
\begin{equation}
\begin{align*}
&   |p(l)|^2 L^2 \equiv \Sigma_{i=0}^{\infty} |p_i(l)|^2 \to 0,
    \quad l \to \infty, \\
& |p(l)| L^1 \equiv \Sigma_{i=0}^{\infty}  |p_i(l)| = \infty.
\end{align*}
\end{equation}
Then clearly $H(p(l))=0$, but $H(0) =1$.

\vspace{3mm}

Let $f:B \to {\bf R}$ be of completely  $C^1$ bounded geometry.
Then the one form $df= \Sigma_{i=1}^{\infty} \frac{\partial f}{\partial x_i} dx_i$ 
can be interpreted as a continuous map:
$$df : B \to H.$$
Then its higher derivatives give the functionals:
$$\nabla^l df: T\bar{B}^{\otimes l} \to H.$$
We define $df$ is of completely bounded geometry, 
if $f$ is  of completely  $C^1$ bounded geometry, and its higher derivatives 
$\nabla^l df$ give continuous maps with respect to the operator norms  for all $l \geq 0$.

 \begin{lem}
 Let  
 $D: C^0(B, H) \mapsto C^0(B,H)$ be  a  pointwise linear  functional,
 and assume it  gives a complete isomorphism. 
 Then its inverse  also gives a complete isomorphism.
  \end{lem}
  {\em Proof:}
  By the assumption, the inverse:
  $$D^{-1}: C^0(B, H) \mapsto C^0(B,H)$$ 
satisfy the equalities:
$$\nabla (D^{-1} ) = -D^{-1} \circ \nabla D \circ D^{-1}.$$
It is immediate to see that $D^{-1}$ is also of completely bounded geometry.

This completes the proof.
\vspace{3mm} \\
{\bf 1.A.2 Local charts:}
Let $(M_0,g_0)  \subset (M_1, g_1)  \subset \dots $ be embeddings of Riemannian
 manifolds with dim$M_i = 2d_i$, where we assume
 the compatibility condition $g_{i+1}|M_i = g_i$.
 We will denote such families by $[(M_i,g_i)]$.
 For $p,q \in M_i$, let us denote their distance in $M\equiv  \cup_i M_i $ by:
 $$d(p,q) \equiv \inf_{j \geq i} d_j(p,q) .$$
We denote $\epsilon$ tublar neighbourhood of $M_i$ by
 $U_{\epsilon}(M_i) \subset M $:
 $$U_{\epsilon}(M_i)= \{ m \in M :  d (m, M_i) < \epsilon \}$$

Let us recall 
$D^{2i}(\epsilon) \subset {\bf R}^{2i}$ 
and $ \bar{D}(p)$
be as in $1.A.1$.
Below we regard the Riemannian metrics $g = \{ g_i\}_i$ as the pointwise operator 
 over the local charts   $T_p D \equiv \cup_i T_p D^{2i}(\epsilon)$, $p \in D$.
If $g$ is of completely bounded geometry at $p \in D$,
then one can extend it to an operator  on the  Hilbert families
$ T_p \bar{D}(p)$.

\begin{defn}
The Riemannian family $\{ g_i \}_i$ is uniformly bounded,  if
 the following conditions are satisfied.  There exists a positive 
$\epsilon  >0$ such that:

\vspace{3mm}

(1) Every point $p \in M \equiv  \cup_i M_i$ admits
a  stratified local chart:
\begin{equation}
\begin{align*}
&  D^{2i}( \epsilon ) \subset D^{2(i+1)}( \epsilon ) \subset  
 \dots \subset  {\bf R}^{2\infty}, \quad
 \varphi(p) :D(\epsilon) \equiv \cup_i D^{2i}(\epsilon)
     \hookrightarrow  \cup_i M_i, \\
& \varphi(p)_i \equiv \varphi(p)|D^{2d_i}(\epsilon)  \mapsto M_i, 
\quad \varphi(p)(0) =p.
\end{align*}
\end{equation}
 
\vspace{3mm}

(2) With respect to $\varphi(p)$ as above,
the induced  Riemannian metrics $\{ g_i\}_i$ are uniformly 
 of completely bounded geometry.
   This means that for any $l >0$, there are constants $C(l) \geq 0$
   independent of $p$
so that  the estimates hold:
$$\text{Sup}_{p \in M} \text{Sup}_{m \in  \cup_i D^{2i}(\epsilon)}
\Sigma_{0 \leq k \leq l}
|\nabla^k ( \varphi(p)^*g) |(m) \leq C(l) \quad (*)$$

(3) There is an increasing and proper function
$h:(0, \infty) \to (0, \infty)$
so that for any $i$, any  pairs of points $p,q \in M_i$ satisfy
the uniformly bounded distance property:
$$d(p,q) \geq h(d_i(p,q))$$
 where $d_i$ and $d$ are the distances  on $ M_i$ and $M = \cup_i M_i$ 
respectively.
\end{defn}

\vspace{3mm}

We will say that the stratified local charts as above  are
{\em complete local charts}.
Also the above family $\{ (p, \varphi(p))\}$ 
will be called a {\em uniformly bounded covering}.
Later on uniform implies independence of choice of points as above.

Let $f: M=  \cup_i M_i \mapsto {\bf R}$ be a bounded smooth function and:
 $$ \varphi(p)^*(f) : D(\epsilon) \to {\bf R}$$
  be the induced functions with respect to the uniformly bounded covering.
We say that
$f$ is of {\em completely bounded geometry}, if 
all $\varphi(p)^*(f)$ are uniformly of completely bounded geometry as:
$$||f||C^k(M) \equiv
\sup_{p \in M} ||\varphi(p)^*(f)||C^k(D(\epsilon)) \leq C_k$$
for
 all $k =0,1,2 \dots$
 and  some constants $C_k$ 
  independent of $p \in M$.

\begin{lem}
Let $[(M_i,g_i)]$ be   Riemannian embeddings as above with
a uniformly bounded covering $\{(p, \varphi(p))\}$ with $\epsilon >0$.
Then exp$_p: \bar{D}(\epsilon') \mapsto \bar{D}(\epsilon)$ exists and is smooth,
 with respect to   the induced Riamannian metrics $\varphi(p)^*(g)$,
where  $\epsilon >  \epsilon' >0   $ and we regard $  D(\epsilon') \subset T_0\bar{D}(\epsilon) $.
\end{lem}
For a proof, see [Kl] ($p 57$, $p72$).
Notice that the geodesic coordinate does not preserve the stratifications.

\vspace{3mm}

   Let
       $f_n, g : M  \to { \bf R}$ be a smooth and bounded family of functions
       for  $n=0,1,2,  \dots$.
      We say that $\{f_n\}_n$ {\em converges weakly} to $g$ in $C^l$,
      if  the restrictions converge in $C^l$ for all $k=0,1,2, \dots$:
      $$f_n|M_k \to g|M_k.$$

   \begin{lem}
   Let $[(M_i, g_i)]$ be a uniformly bounded Riemannian family
   such that each $M_k$ is compact.
   Let
       $f_n : M  \to { \bf R}$ be a family of smooth functions 
       of completely bounded geometry
       for  $n=0,1,2,  \dots$
       
       Suppose $C^{l+1}$ norms are 
       uniformly bounded:
       $$||f_n|| C^{l+1}(M)  \leq Const(l).$$
    Then a subsequence $f_{n_j}$ weakly converges in $C^l$ to a smooth function
    $g:  M \to { \bf R}$ of completely bounded geometry.
       \end{lem}
{\em Proof:}
By the condition, the restrictions $\{f_n|M_k\}_n$ satisfy uniformity
of $C^{l+1}$ norms $||f_n||C^{l+1}(M_k) \leq C(l)$.

It follows from  Rellich lemma that 
there is some $C^l$ function $g_1: M_1 \to {\bf R}$ 
so that  a subsequence $\{f_{n(i)}|M_1\}_i$ converges
to $g$ in $C^l(M_1)$.

By the same argument, 
there is some $C^l$ function $g_2: M_2 \to {\bf R}$ 
so that  a subsequence $\{f_{n(i,2)}|M_2\}_i$ converges
to $g_2$ in $C^l(M_2)$ for 
another subsequence $\{n(i,2) \}_i \subset \{n(i)\}_i$.
Clearly $g_2|M_1=g_1$ holds.

By choosing subsequences successively, 
$\{f_{n(i,k)}|M_k\}_i$ converge to some  $C^l$ function $g_k: M_k \to {\bf R}$,
with $g_k|M_{k-1}=g_{k-1}$. 
These satisfy uniformity of $C^l$ norms $||g_k||C^l(M_k) \leq c < \infty$.

Let  $g : M \to {\bf R}$ be a smooth and bounded function
defined by $g|M_k \equiv g_k$.
Then the subsequence $\{f_{n(i,i)}\}_i$ converge weakly to $g$ in $C^l$.

 This completes the proof.
  \vspace{3mm} \\
{\bf 1.B Almost Kaehler sequence:}
Let $(M, \omega , J)$ be a finite dimensional symplectic manifold
 equipped with a compatible almost
complex structure.
 Namely $g( \quad , \quad ) = \omega( \quad , J \quad )$ gives
a Riemannian metric on $M$.
Such a manifold is called  an {\em almost Kaehler manifold}.

Let $ (M_0,\omega_0,J_0) \subset (M_1, \omega_1, J_1) 
      \subset \dots \subset  (M_i, \omega_i, J_i) \subset \dots$
be  infinite embeddings of almost Kaehler manifolds.
  If one denotes the inclusion $I(i): M_i \hookrightarrow M_{i+1}$,
 then the above implies
$\{ I(i)\}_i $ gives a family of holomorphic embeddings,
$J_{i+1} \circ I(i)_* = I(i)_* \circ J_i$, and
 the symplectic forms  $I(i)^*(\omega_{i+1})=\omega_i$ are 
the restrictions.

 Suppose  dim $M_k=2d_k$, and
    $U_{\epsilon}(M_i) \subset M \equiv  \cup_j M_j $
be   $\epsilon$ tublar neighbourhoods of $M_i$. Let:
 $$\pi_k' : \cup_i D^{2i}(\epsilon) \mapsto D^{2d_k}(\epsilon)$$
be the standard projections.
\begin{defn}
An almost   Kaehler sequence
  consists of a family of   embeddings  by  almost Kaehler manifolds:
$$[(M_i,\omega_i,J_i)] = (M_0,\omega_0,J_0) \subset (M_1, \omega_1, J_1) 
      \subset \dots \subset  (M_i, \omega_i, J_i) \subset \dots$$
   and a positive  $\epsilon >0$
    so that    any points $p \in M \equiv \cup_i M_i$  admit  $\epsilon$
    uniformly bounded coverings $\{(p, \varphi(p))\}$
    which satisfy the followings:

\vspace{3mm}

(1)  $\omega|\cup_i D^{2i}(\epsilon)$ and 
  $J|\cup_i D^{2i}(\epsilon)$ are
 of completely bounded geometry.

(2) The induced symplectic form satisfies:
$$\varphi(p)^*(\omega) = 
    \frac{\sqrt{-1}}{2}
    \Sigma_{i=0}^{\infty} \   dw_i \wedge d\bar{w}_i \quad  \text{ at } p$$
where $(w_1, \dots, w_i)$ are the coordinates
 on $D^{2i}(\epsilon) \subset {\bf C}^i$.

 (3)
There are     families of holomorphic  maps:
$$\pi_k : U_{\epsilon}(M_k) \mapsto M_k$$
which satisfy the following properties:
$$ \pi_k |  M_k  = \text{ id }, \quad
  \pi_k(\varphi(p)(m) ) = \varphi(p)(\pi_k'(m)) \quad \text{ for all } 
m \in U_{\epsilon}(M_k).$$
\end{defn}
Uniformly bounded coverings which satisfy the properties  (1) (2) (3)  above, 
are called $\epsilon$ {\em complete almost Kaehler charts}.

An almost Kaehler data $\{(\omega_i, J_i)\}$ gives a uniformly bounded and
 compatible family of  Riemannian metrics on $\{ M_i\}_i$.
Notice that the equalities $< v,u> = <\pi_k(v),u>$ hold
 for $u \in T_p M_k$ and $v \in T_p U_{\epsilon}(M_k)$
 with respect to the induced Riemannian metric.

\vspace{3mm}
        Later on, we fix   a uniformly bounded covering
by $\epsilon$ complete almost Kaehler charts.

We will say that $[(M_i, \omega_i, J_i)]$ is  a {\em Kaehler sequence},
if it is an  almost Kaehler sequence  consisted by
a  uniformly bounded covering 
by  holomorphic
complete  Kaehler charts  $ \varphi(p) $ at  all points $p$,
where we equip with  the standard complex structure on $\cup_i D^{2i}(\epsilon)$
(see [GH] $p107$).

Let  $f: M=\cup_k M_k \to {\bf R}$ be a smooth bounded function
over an almost Kaehler sequence.
 We will say  that $f$  is a {\em bounded Hamiltonian function},
 if it is of completely bounded geometry.

\vspace{3mm}

Let $(M, \omega)$ be a finite dimensional symplectic manifold. 
The following  facts are  well known ([G1]):

\vspace{3mm}

(1) there exist compatible almost complex structures, and 

(2) the space of compatible almost complex structures
is contractible.
\vspace{3mm} \\
In the infinite dimensional situation, the condition (1) depends on the spaces,
but the same  thing holds for (2) for a fixed uniformly bounded covering.

\begin{lem}
Let $[(M_i, \omega_i)]$ be   a  symplectic sequence.
Suppose there exists a family of  compatible almost complex structures
$\{ J_i\}_i$ so that $[(M_i, \omega_i, J_i)]$ consists of an almost Kaehler sequence
with respect to a uniformly bounded covering $\{(p, \varphi(p))\}$.
Then the space of such family:
\begin{align*}
{\frak J}([(M_i, \omega_i)]) =&  \{  \ \{ J_i \}_i :  [(M_i, \omega_i, J_i)] : \\
& \text{ almost Kaehler  sequence with respect to } \{(p, \varphi(p))\}  \   \}
\end{align*}
is contractible.
\end{lem}
$Proof:$ 
We follow a well known argument in the  finite dimensional case.

Let us choose a reference family of almost complex structures $\{ J_i^0 \}_i$.
Take another one $\{ J_i^1\}_i$. Let us connect these by a compatible family of almost complex structures
$\{ J_i^t\}_i$, $t \in [0,1]$. For $\alpha = 0$ or $1$,
 let us put $h_i^{\alpha}( \quad , \quad ) = \omega_i(\quad , J_i^{\alpha} \quad)$.
Then $\{h_i^{\alpha}  \}_i$ gives a family of Riemannian metrics. 
Moreover each $J^{\alpha}_i$
is uniquely determined by $h_i^{\alpha}$. Now let us consider a smooth family of Riemannian metrics:
$$h_i^t = h_i^0 + t( h_i^1 - h_i^0).$$
Then for each $i$, there exists a unique and smooth family of compatible almost complex structures
$J_i^t$, $t \in [0,1]$ satisfying $h_i^t (\quad , \quad) = \omega_i(\quad,  J_i^t \quad)$.
 
Let us choose 
 a complete almost Kaehler  chart at  $p \in M_i \subset M_{i+1}$:
 $$\omega_i = \Sigma_{j \leq i}  \ dp_j \wedge dq_j, 
        \quad \omega_{i+1} = \Sigma_{j \leq i+1}  \ dp_j \wedge dq_j
\quad  \text{ at } p$$
  and denote the  local projections by $\pi'_i: D^{2d_{i+1}}(\epsilon) \mapsto D^{2d_i}(\epsilon)$.
  Let us   check  the compatibility condition
$J_{i+1}^t \circ I(i)_* = I(i)_* \circ J_i^t$  at $p$ and for each $t$.
Let us take  $v_i \in T_p M_i$. Then:
\begin{equation}
\begin{align*} 
 \omega_{i+1}(\quad ,   J_{i+1}^t v_i)  
       &    = h_{i+1}^0(\quad , v_i) + t( h_{i+1}^1 (\quad , v_i) - h_{i+1}^0(\quad , v_i)) \\
&     = \omega_{i+1}(\quad , J_{i+1}^0  v_i) + t\{ \omega_{i+1}(\quad , J_{i+1}^1  v_i) -
                                \omega_{i+1}(\quad , J_{i+1}^0  v_i) \} \\
 & = \omega_{i+1}(\quad, J_i^0v_i) + t \{    \omega_{i+1}(\quad, J_i^1 v_i) -
                                                            \omega_{i+1}(\quad, J_i^0 v_i)   \} \\
                                                            & =  \omega_{i+1}(\quad ,   J_i^t v_i)  \\
  & = \omega_i(\pi'_i  \quad , J_i^0 v_i) + t \{    \omega_i( \pi'_i \quad, J_i^1 v_i) -
                                                            \omega_i( \pi'_i \quad, J_i^0 v_i)   \} \\
  & = \omega_i( \pi'_i \quad , J_i^t v_i).
\end{align*}
\end{equation}     
The fourth equality implies the 
the  compatibility condition.

Moreover  the following equalities hold
 from the equality between the first and the last above:
\begin{align*}
&    \omega_i( \pi'_i  J_{i+1}^t(w), J_i^t \pi'_i(v)) =
\omega_{i+1}(J_{i+1}^t(w) ,   J_{i+1}^t \pi'_i(v)) \\
& = \omega_{i+1}(w ,    \pi'_i(v)) =  \omega_i(\pi_i'(w) ,    \pi'_i(v))
=  \omega_i(J_i^t(\pi_i'(w)) ,    J_i^t(\pi'_i(v))).
\end{align*}
This implies the equality: 
$$ \pi'_i  J_{i+1}^t =    J_i^t \pi_i'$$
and so $\pi_i$ are holomorphic with respect to $J^t$.
This completes the proof.
\vspace{3mm} \\
{\em Remark 1.1:} 
It is not clear whether the conclusion might still 
 hold when we do not fix a uniformly bounded covering.
\vspace{3mm} \\
{\bf 1.B.2 Embeddings of almost Kaehler sequences:}
Let  $[(M_i, \omega_i, J_i) ]$ 
be an almost Kaehler sequence equipped with 
 complete local charts $\varphi(p) :  \cup_{s \geq 1} D^{2s}(\epsilon)
\hookrightarrow \text{ im } \varphi(p) \subset  M$ for all $p \in M$.

Let us say that $[(M_i', \omega_i', J_i') ]$ is 
 {\em embeddable}   into   $[(M_i, \omega_i, J_i) ]$, 
if there are positive $\epsilon >0$, 
  subindices $\{ k(i)\}_i$ with $k(i) \geq i$ and 
  compatible embeddings between almost Kaehler manifolds:
 $$ I_i: (M_i' , \omega_i', J_i') \hookrightarrow (M_{k(i)}, \omega_{k(i)}, J_{k(i)})$$
 so that
  there 
 are     families of holomorphic  maps:
$$\pi_i :  U_i  \mapsto M_i'$$
from    the  open 
$\epsilon$  tublar  neighbourhoods 
$I_i (M_i') \subset U_i \subset M \equiv \cup_i M_i$, 
 which satisfy the properties:
$$ \pi_i  |  M_i'  = \text{ id }, \quad
  \pi_i(\varphi(p)(m) ) = \varphi(p)(\tilde{\pi}_i(m))$$
  for all  
$m \in U_i$, where $\tilde{\pi}_i :   \cup_{s \geq 1} D^{2s}(\epsilon)  \to D^{2d_i}(\epsilon)$
are the projections with $d_i = \dim M_i'$.
  \vspace{3mm} \\
{\em Example 1.2:}
Let us fix $p \geq 1$ and consider the canonical embeddings of the Grassmannians
$Gr_{p,q} \hookrightarrow Gr_{p,q+1}$
which embed  each $p$ plane $L \subset {\bf C}^{p+q} \subset {\bf C}^{p+q+1}$.
These admit the canonical and compatible Kaehler forms, and 
the direct limits $Gr_p \equiv \lim_{q \to \infty} Gr_{p,q}$ consiste of the Kaehler sequences.

Let us consider the Pl\"ucker embedding
$Gr_{p,q} \hookrightarrow {\bf CP}^N$,
where $N=N(p,q)=
\begin{pmatrix}
&p+q \\
& p
\end{pmatrix} -1$, which associate each $p$ plane $L \subset {\bf C}^{p+q}$
and its basis $\{v_1, \dots, v_p\} $ to the complex line
$[v_1 \wedge \dots \wedge v_p] \in {\bf CP}^N$.

It is well known that these embeddings preserve the canonical Kaehler forms, and so 
they give the embedding of the Kaehler sequences:
$$I: [Gr_{p,q}] \hookrightarrow [{\bf CP}^n]$$
where $(M_i,\omega_i, J_i) = Gr_{p,i}$ and 
$(M_i',\omega_i', J_i')={\bf CP}^i$ with $k(i) =N(p,i)$.

Moreover the Schubert calculus verifies the isomorphisms:
$$I_* : H_2(Gr_{p,q}: {\bf Z}) \cong H_2({\bf CP}^N: {\bf Z}) \cong {\bf Z}$$
between simply connected spaces.
\quad
\vspace{3mm} \\
{\em Example 1.3:}
Let us consider the {\em Veronese maps} defined as follows.
Let us introduce the lexicographic order for
 two indices $(i_0, \dots, i_n)$ and $(j_0, \dots, j_p)$.

Let us fix $m \in\{1,2, \dots\}$, and 
take ${\bf CP}^n$ with the homogeneous coordinate $[z_0, \dots,z_n]$.
For $N = \begin{pmatrix}
& n+m \\
& m
\end{pmatrix}
-1$, we define the Veronese map:
\begin{align*}
& v_m : {\bf CP}^n \hookrightarrow {\bf CP}^N, \\
& v_m([z_0, \dots,z_n]) = \{ z_0^{i_0}, \dots z_n^{i_n} : \Sigma_{l=0}^n i_l =m\} .
\end{align*}

With $n_1=1$, let us  define  numbers inductively by
$n_{i+1} = 
\begin{pmatrix}
& n_i+m \\
& m
\end{pmatrix}
-1$.

Now we have two different embeddings:
$${\bf CP}^{n_i} \subset_v {\bf CP}^{n_{i+1}}, \quad {\bf CP}^{n_i} \subset {\bf CP}^{n_{i+1}}$$
where the first is the given by the Veronese map and the second is by the canonical embedding.

\begin{lem} The following diagram commutes:
$$
\begin{matrix}
{\bf CP}^{n_i} & \subset_v & {\bf CP}^{n_{i+1}} \\
\cap && \cap \\
{\bf CP}^{n_{i+1}} & \subset_v & {\bf CP}^{n_{i+2}}
\end{matrix}$$
\end{lem}
{\em Proof:}
This follows since we have used  the lexicographic order
for the coordinates. 
This completes the proof.

\vspace{3mm}

\begin{cor}
There is a canonical embeddings of ${\bf CP}^{\infty}$ to itself:
$$v_m : {\bf CP}^{\infty} \subset_v {\bf CP}^{\infty}$$
of degree $m$, so that the restrictions are given by the Veronese maps.
\end{cor} 
{\em  Remark 1.2:}
We have the Veronese  sequence by the
embeddings by the iterations of the Veronese maps:
$${\bf CP}^{n_1} \subset_v {\bf CP}^{n_2} \subset_v  \dots \subset_v {\bf CP}^{n_l} \subset 
\dots \subset  V \equiv  \cup_i {\bf CP}^{n_i}  .$$
This is not Kaehler sequence, since the degree grows unboundedly in the total space.
Study  of this embeddings will require much harder analysis.
\vspace{3mm} \\
        {\bf 1.B.3 Symmetric almost Kaehler sequence:}
Let us  introduce geometric classes of almost Kaehler sequences.
    Their symmetric properties  will allow us to analyze structure of  holomorphic maps.

\begin{defn}
$[(M_i, \omega_i, J_i)]$
 is  a {\em symmetric almost Kaehler sequence}, if
there are   $\epsilon >0$, 
uniformly bounded coverings   $\{ \varphi(p) \}$ at any $p \in M$ 
and some $l =l(k)  >k$  for any $ k \geq 0$
so that 
there are
  families of almost Kaehler submanifolds $M_k \subset W_i \subset M_i$ 
  with $W_l = M_l$ for all $i \geq l$ and
  isomorphisms: 
$$P_i: \{(M, M_k), \omega, J\}  \cong \{(M,M_k), \omega, J\}$$
which preserve $M_k$  and transform $W_i$ to $M_l$ as:
$$P_i: (W_i, \omega_i|W_i, J_i|W_i) \cong (M_l,\omega_l, J_l)$$
such that  at any $p \in M_k$:
$$D :   TM_k \oplus_{i \geq l}  N_{i,k}  \cong T M |M_k $$
 give the uniform isomorphisms over $M_k$
with respect to the complete local charts, 
where:
\begin{align*}
& N_{i,k}= (P_i^{-1})_* [  \ (\text{ Ker } ( \pi_k)_* \cap TM_l) |M_k) \ ] \\
& D_p = \text{ id } \oplus (P_l)_* \oplus (P_{l+1})_* \oplus  \dots 
\end{align*}
\end{defn}
We say that the family of maps $\{(P_i, \pi_k)\}_{i,k}$ give symmetry
of the almost Kaehler sequence with respect to $(M_k,M_l)$.

If all these properties hold by use of complex structure,
then we say that it is a symmetric Kaehler sequence.

\vspace{3mm}

Suppose  $[(M_i,\omega_i, J_i)]$ is symmetric.
It is {\em isotropic}, 
if there are families of parametrized  isomorphisms  for each  $0 \leq t \leq 1$:
$$P^t_i: \{(M, M_k), \omega, J\}  \cong \{(M,M_k), \omega, J\}$$
 with:
$$P_i^0\equiv  \text{ id }, \quad P^1_i =P_i.$$
\quad
  \vspace{3mm}  \\
{\em Examples 1.4:}
(1) Let $(X,\omega, J)$ and $(Y, \tau, I)$ be two almost Kaehler manifolds, 
and choose a base point $y_0 \in Y$. 
Let us consider the products:
$$(X \times Y_1 \times Y_2 \times \dots, \ \omega+\tau_1 +\tau_2 + \dots, \ J \oplus I_1 \oplus I_2 \oplus  \dots)$$
where all $(Y_i,\tau_i,I_i)$ are the same $(Y,\tau,I)$, and 
 we embed $X \times Y_1 \subset X \times Y_1 \times Y_2$ 
by identifying $ X \times  Y = X \times Y \times \{y_0\}$ and similar for others.

The  infinite product sequence admits symmetric structure  by
 choosing:
  $$M_k = X \times Y_1 \times \dots \times Y_k , \quad 
  W_i =M_k  \times \{ y_0 \} \times  \dots \times \{y_0\} \times Y_i.$$
 $P_i$ are given by the obvious exchange of the coordinates.
 
 \vspace{3mm}

(2) Let $ M$ be a complex manifold,
and take any   holomorphic curve $u: {\bf C}P^1 \mapsto M$.
Then the holomorphic vector bundle $u^*(TM) \mapsto {\bf C}P^1$
splits as the direct sum of  holomorphic line bundles.
This fact can be regarded as
 `infinitesimal  symmetric  property' (see [OSS]).
  
\vspace{3mm}

(3)
Let us consider the projective spaces with the Fubini Study forms
$[({\bf CP}^i, \omega_i)] = {\bf CP}^1 \subset {\bf CP}^2 \subset \dots \subset {\bf CP}^n \subset \dots {\bf CP}^{\infty}$.  
This is an isotropic  symmetric Kaehler sequence, and we 
 denote it  by ${\bf C}P^{\infty}$.
There are standard charts ${\bf C}^i \subset {\bf CP}^i$ and $\omega_i$ can be expressed as:
$$\omega_i |{\bf C}^i =  \frac{ \sqrt{-1}}{2}
    [ \frac{ \Sigma_l dw_l \wedge d\bar{w}_l }{(1+ w\bar{w})}
     - \frac{ ( \Sigma_l \bar{w}_ldw_l ) \wedge 
 (\Sigma_l w_l d\bar{w}_l)}{(1+ w\bar{w})^{2}}]$$
where $w=(w_1, \dots,w_i)$ are the coordinates on ${\bf C}^i$.
Then $[\omega_i \equiv \omega|D^{2i}]$
is completely of bounded geometry, where $D^{2i} \subset {\bf C}^i$
are the unit balls.
In order to obtain another charts at  any $p \in {\bf C}P^i$,
one may use any constant unitary matrix $U \in Mat_{i+1}({\bf C})$
with $U([1,0, \dots, 0]) = p \in {\bf C}P^i $.

Let $U_{\epsilon}({\bf C}P^k) \subset \cup_i  {\bf C}P^i$
be $\epsilon$ tublar neighbourhood. Then there are natural projections,
$\pi_k : U_{\epsilon}({\bf C}P^k) \mapsto {\bf C}P^k$:
$$ \pi_k([z_0, \dots, z_k,z_{k+1},  \dots]) 
= [z_0, \dots, z_k,0 , \dots].  $$

Let us put
$M_k ={\bf CP}^k$ and $W_i$ by:
$$W_i =\{[z_0: \dots: z_k: 0  \dots : 0 : z_i: 0 : 0 : \dots ] \in {\bf CP}^{\infty} \}$$
for all $i \geq l= k+1$ with  $M_l = {\bf CP}^{k+1}$.
$P_i :  W_i \cong {\bf CP}^{k+1}$ are given just by  exchange of  the coordinates:
$$[z_0: \dots: z_k: 0: \dots : 0:  z_i: 0 \dots] \to [z_0: \dots: z_k : z_i: 0: \dots].$$
This is isotropic, by putting:
\begin{align*}
 P_i^t( [z_0:  \dots: & z_k:  \dots]) 
=  [z_0: \dots:z_k  
 : \cos  \frac{\pi t}{2}  z_{k+1}+\sin \frac{\pi t }{2} z_i : \\ 
 & z_{k+2}:  \dots : z_{i-1}: 
 - \sin \frac{\pi t}{2}  z_{k+1}+\cos \frac{\pi t}{2} z_i : z_{i+1}: \dots].
\end{align*}

(4)
There are many variants. For example one can change ${\bf C}$ by ${\bf H}$.
For others, let us consider the Grassmannians:
$$Gr_{r,n}({\bf C}) = \{ H \subset {\bf C}^{r+n} \ ; \ H : r \text{ dimensional } {\bf C} \text{ vector subspaces } \}.$$
One can canonically embed as $H \subset {\bf C}^{r+n+1}$, and by taking the direct limit,
one obtains the Kaehler sequence
$Gr_r({\bf C})= \lim_{n \to \infty} Gr_{r,n}({\bf C})$
equipped with the standard Kaehler structure.

This space also admits isotropic and symmetric structure. Let us put:
$${\bf C}^{k,i}=\{ (z_1, \dots, z_k, 0, \dots, 0, z_{k+i}) : z_j  \in {\bf C} \} \subset  {\bf C}^{k+i}$$
and choose
$M_k =Gr_{r,k}$ and $W_i\equiv W^r_{k,i}$ are consisted by all elements of the form:
$$W_{k,i}^r =\{H  \subset {\bf C}^{k+r,i}   \ ; \ H : r \text{ dimensional } {\bf C} \text{ vector subspaces } \}.$$
 The required isomorphisms and isotropies
can be obtained by the same way as the above $P_i$ and $P_i^t$.

\begin{lem}
Let $[(M_i, \omega_i, J_i)]$ be a symmetric   almost Kaehler sequence, 
and choose any   pair $(M_k,M_l)$ as above.

Then there is a bundle $N_{k,l} \to M_k$ 
with  the uniform  isomorphisms:
$$TM |M_k \cong TM_k \oplus (N_{k,l} \otimes {\bf R}^{\infty})$$
with respect to the complete local charts over $M$.
\end{lem}
{\em Proof:}
Let us put:
$$N_{k,l} =( \text{ Ker } (\pi_k)_* \cap TM_l  )|M_k.$$
There are holomorphic isomorphisms $TM_l|M_k \cong TM_k \oplus N_{k,l}$.
Then the conclusion follows by use of the  isomorphisms of  the tangent bundles:
$$ TM_k \oplus N_{k,l} \cong TW_i|M_k , \quad 
 (v,w) \to (v, (P_i^{-1})_*(w)).$$

This completes the proof.
\vspace{3mm} \\
{\bf 1.B.4 Quasi transitivity:}
Let $[(M_i, \omega_i, J_i)]$ be an almost Kaehler sequence.
We  say  $[(M_i, \omega_i, J_i)]$ is {\em quasi transitive}, if
 for any  $N  >0$,
 there is  $k = k(N) $ so that for any mutually different points 
$p_0, \dots, p_{N-1} \in M \equiv \cup_i M_i$,
 there are automorphisms of the almost Kaehler
sequence $A: (M, \omega , J) \cong (M, \omega, J)$
such that:
$$A(p_i) \in M_k$$
hold for all $0 \leq i \leq N-1$.

\begin{lem}
The infinite  projective space
$[({\bf C}P^i, \omega_i, J_i)]$ is
quasi transitive. 
\end{lem}
{\em Proof:}
Let us construct automorphisms $A^i: (M,M_{l_i}) \cong (M,M_{l_i})$
inductively so that they satisfy the followings:
$$  A^i(p_i) \in M_{l_i}, \quad A^i|M_{l_j} = id$$
for all $j < i$. We put $A \equiv A^{N-1} \circ A^{N-2} \circ \dots \circ A^0$
and $k  =l_{N-1}$.

Let us choose any $p_0=[z_0,z_1, \dots] \in  {\bf C}P^L \subset {\bf C}P^{\infty}$.
Firstly let us move $p_0$ to $ [1,0,0, \dots]$ by
a unitary  automorphism $A^0 \in U(L+1) \subset $ Aut ${\bf C}P^{\infty}$.

Let us consider $u_1= A^0(p_1) \in  {\bf C}P^{\infty}$.
 We put $A^1=$ id, if $u_1 \in {\bf C}P^1$.
Suppose $u_1 =[u_1^0,u_1^1, \dots] \not\in {\bf C}P^1$.
Then   $(u_1^1 , u_1^2, \dots)$ is non zero and so defines an element   in $ {\bf CP}^{\infty}$.
Let us choose another unitary automorphism $T_1$ with
$T_1([u_1^1,u_1^2, \dots]) =[1, 0, \dots]$. 
Then we put $A^1 =$ diag $(1, T_1)$.

Let us consider $u_2= A^1 \circ A^0(p_2) \in  {\bf C}P^{\infty}$.
 We put $A^2=$ id, if $u_2 \in {\bf C}P^2$.
Suppose $u_2 =[u_2^0,u_2^1, \dots] \not\in {\bf C}P^2$.
Then   $(u_2^2 , u_2^3, \dots)$  defines an element   in $ {\bf CP}^{\infty}$.
By  another unitary automorphism $T_2$ with
$T_2([u_2^2,u_2^3, \dots]) =[1, 0, \dots]$ 
Then we put $A^2 =$ diag $(1,1,  T_2)$.

By the same way one can inductively construct $A^3 , \dots, A^{N-1}$.

 This completes the proof.
 \vspace{3mm} \\
{\em Remark 1.3:}
(1) A similar argument can be used to verify that 
the infinite Grassmannians $Gr_N({\bf C}) =\lim_{L \to \infty} Gr_{N,L}$
also satisfy  quasi transitivity.

(2)
  For  our later applications, it is enough to assume 
{\em weakly quasi transitivity}, in the sense that  the above  $k=k(N,L)$ can also depend
on $L = \max_{i,j} d(p_i,p_j)$.

Notice that if the diameter of $M$ are bounded, then it is quasi transitive 
whenever weakly quasi transitive.
\vspace{3mm} \\
 {\bf 1.B.5 Minimality:}
  Let $[(M_i, \omega_i,  J_i)]$ be  an almost Kaehler sequence. 
Let us introduce its  invariant  ([HV]):
 \begin{equation}
 \begin{align*}
 m([(M_i, & \omega_i, J_i)]) 
       = \text{ inf}_u \{ <\omega, u> ; \\
       & u: S^2 \mapsto M \equiv \cup_i M_i 
                                                   : \text{ non constant holomorphic  curve} \}.
 \end{align*}
 \end{equation}
 By restriction to  the symplectic  sequence,
 one obtains the  invariant:
$$m([(M_i, \omega_i)]) = \text{ inf}_{\alpha} 
         \{ <\omega, \alpha> \ >0 : \alpha: S^2 \mapsto \cup_i M_i\}.$$
 We  say      $[(M_i, \omega_i,  J_i)]$  is {\em minimal}, if both  the equality 
 and positivity hold:
$$m([(M_i, \omega_i)]) = m([(M_i, \omega_i , J_i)]) >0 .$$
The minimal  homotopy class 
 $\alpha \in $ Homot$\{ S^2 \mapsto M\}$ is given by the equality
 $\int_{\alpha} \omega = m([(M_i,  \omega_i, J_i)]) $, which 
 plays an important role in section $3$.
\vspace{3mm} \\
{\em Examples 1.5:}
(1)
 Notice that if  $[(M_i, \omega_i, J_i)]$ is of  $\pi_2$ rank $1$,  
 $\pi_2(\cup_i M_i) / $ Tor $\cong {\bf Z}$,
 then minimality is equivalent to existence of
non constant holomorphic curves representing $1 \in  \pi_2 / $ Tor.

  The Fubini Study form on ${CP}^n$ with the standard complex structure
 gives $\pi_2$ rank one minimal data $(\omega, J)$ with $m =  \pi$.

 (2) Let $({\bf CP}^1, \omega, J)$ be the standard curve and
 $[(M_i, \omega_i, J_i)]$ be minimal. Then the product
 $[((M_i \times {\bf CP}^1, \omega_i+ \omega, J_i \oplus J)]$
 is also minimal.

(3) Suppose $[(M_i, \omega_i, J_i)]$ is algebraic with 
each $\omega_i \in H^2(M_i : {\bf Z})$.
Then it is  minimal, if any generating elements  in $H_2(M: {\bf Z})$
can be represented by  some holomorphic curves.
In particular it is the case when it  is simply connected, algebraic,
and any generating elements  in $\pi_2(M)$
can be represented by  some holomorphic curves.
\vspace{3mm} \\
{\bf 1.C Transition functions:}
Let $[(M_i, \omega_i, J_i)]$ be an almost Kaehler sequence, and choose
 a uniformly bounded covering $\{(p, \varphi(p))\}$.

For any $p,p' \in M=  \cup_i M_i$, let us put:
$$\Phi(p,p') \equiv \varphi(p')^{-1} \circ \varphi(p):
B(p,p') \cong B(p',p)$$ where
 $B(p,p') \equiv 
\varphi(p)^{-1}( $ im $\varphi(p)  \cap $ im $\varphi(p')) 
\subset \cup_i D^{2i}(\epsilon)$.

\begin{lem}  Let $[(M_i, \omega_i, J_i)]$ 
be an almost Kaehler sequence.
Then 
$\Phi(p,p')$ give the complete isomorphisms:
$$\sup_{p,p'}| |\nabla^l \Phi(p,p')||C^0 \leq c_l$$
for constants $c_l$ independently of $p,p'$ for all $l \geq 0$.
\end{lem}
{\em Proof:} 
See  def $1.1$ and below it  for the terminologies.

This can be seen by use of the exponential mappings 
([Kl] $p74$). For convenience we check the uniform estimates for $ l \leq 2$,
which will be  used to construct diffemorphisms in lemma $5.3$ below.

The estimates on the first derivative come from    uniformity of
compatible Riemannian metrics.

Let us put $g(p) = \varphi(p)^*(g)$ on $\cup_i D^{2i}(\epsilon')$.
Thus $\Phi(p,p')^*(g(p')) = g(p)$ holds.
Since both of $g(p)$ and $g(p')$ are uniformly bounded, 
the conclusion  holds for $l= 1$.
Let us put $\Phi=  \Phi(p,p')$.

Suppose $\nabla \Phi$ could be unbounded at $m$, 
  and take  
a $C^1$ curve $\gamma : [-1,1] \to B(p,p')$
with $\gamma(0) =m \in B(p,p') $.

 Let  $X_0 \in T_m M$ be a unit tangent vector with $|X_0|=1$,
and  extend $X_0$ as the vector fields  $X$ along $\gamma$
by constant.
Then we  put
$W_0=  \frac{d \Phi_*(X_s)}{ds}|_{s=0}  \in T_{\Phi(m)} M$
and  extend  $W_0 $ to
   $W$  along $\Phi(\gamma)$ similarly.
 We also have another vector field $Z= \Phi^{-1}(W)$ along $\gamma$.
 Notice $\frac{d X}{ds} = \frac{dW}{ds} =0$ and 
 $ \frac{dZ}{ds}|_{s=0} = X_0$.

Let us  consider the equalities:
$$g(p)(X, Z) = \Phi^*(g(p'))(X,Z) = g(p')(\Phi_*(X),W) .$$
By differentiating both sides at $s=0$, 
one has the equality:
$$g(p)' (X_0, \Phi^{-1}(W_0) ) + g(p)(X_0,X_0) 
= g(p')'(\Phi_*(X_0),W_0) +
g(p')(W_0,W_0) $$
where $g(p)' = d g(p)(\gamma(s))/ ds |_{s=0}$.
So there is a constant $C$ with the point wise estimates:
$$||\nabla \Phi_*(X_0)||^2 \leq C \{  ||\nabla \Phi_*(X_0)|| + 1\}$$
which contradicts to the assumption.

\section{ Moduli spaces  of   holomorphic curves}
In this section we study theory of holomorphic curves into almost Kaehler sequences.
In particular we develop the analytic tools to 
construct finite dimensional moduli spaces over   sequences which satisfy 
some symmetric properties.
\vspace{3mm} \\
{\bf 2.A Finite dimensional preliminaries:}
 We recall basics of moduli theory of holomorphic curves into finite dimensional 
  symplectic manifolds.
  Most  of the contents in $2.A$  are already in  [HV], which 
  are preliminaries for  $2.B$ where
 we  formulate Sobolev spaces over the infinite dimensional spaces  $M= \cup_i M_i$.

 ${\bf CP}^1$ has particular points $0, \infty \in {\bf CP}^1$,
 and let $ 0 \in D(1) \subset S^2= {\bf C}P^1$ be the hemi sphere. 
 We choose and fix the following data:
 \vspace{2mm} 
 
(1) a large $l \geq 1$,  

(2) a non trivial homotopy class  $\alpha \in  \pi_2(\cup_i M_i)$, and

(3)  different fixed points
  $p_0, p_{\infty} \in M_0 \subset M \equiv  \cup_i M_i$.
\vspace{2mm} \\
Let $L^2_{l+1}(S^2, M_i)$ be the sets of
$L^2_{l+1}$ maps from $S^2$ to $M_i$.
(In $2.B$, we will define them in detail).
  Let us put the spaces of Sobolev maps:
\begin{equation}
\begin{align*}
 {\frak B}_i \equiv  {\frak B}_i (\alpha) & = \{  u   \in   L^2_{l+1}(S^2, M_i) :  
 [u] = \alpha ,   \\
 & \int_{D(1)} u^*(\omega) = \frac{1}{2}<\omega, \alpha>, \quad
  u(*) = p_* \in M_0  , * \in \{ 0, \infty \}  \}.
  \end{align*}
  \end{equation}

 Let $E(J)_i, \ F_i \mapsto S^2 \times M_i$ be 
 vector bundles whose fibers are respectively:
\begin{equation} 
\begin{align*}
& E(J)_i(z,m) = \{ \phi: T_z S^2 \mapsto T_mM_i : 
   \text{ anti complex linear } \}, \\
& F_i(z,m) = \{ \phi: T_z S^2 \mapsto T_mM_i : \text{ linear} \}.
\end{align*}
\end{equation}
Then we have   two stratified Hilbert bundles over ${\frak B}_i$:
 \begin{equation}
 \begin{align*}
&  {\frak E}_i = L^2_l({\frak B}_i^*(E(J)_i)) = 
\cup_{u \in {\frak B}_i} \{u\} \times L^2_l(u^*(E(J)_i)), \\
&   {\frak F}_i = L^2_l({\frak B}_i^*(F_i)) 
= \cup_{u \in {\frak B}_i} \{u\} \times L^2_l(u^*(F_i)).
\end{align*}
\end{equation}
On all of these Hilbert manifolds ${\frak B}_i, {\frak E}_i , {\frak F}_i$,
there exist compatible, free and continuous  $S^1$ actions 
which are induced from the one on ${\bf C} \subset {\bf C}P^1$.
\vspace{3mm} \\
{\em Remark 2.1:}
One may regard that 
 $E(J) = \cup_i E(J)_i$  and $F = \cup_i F_i$ are stratified vector bundles
 over $S^2 \times M $,  $M= \cup_i M_i$, and so 
  we have   stratified Hilbert bundles
$ {\frak E} = \cup_i {\frak E}_i$ and  $     {\frak F} = \cup_i {\frak F}_i $
over $ {\frak B} \equiv  \cup_i {\frak B}_i$.

\vspace{3mm}

 The non linear {\em Cauchy-Riemann operators} and their tangent maps
 are defined respectively  as  sections:
 \begin{equation}
 \begin{align*}
 & \bar{\partial}_J \in C^{\infty}({\frak E}_i \mapsto {\frak B}_i),  
\quad
  \bar{\partial}_J (u) = Tu + J \circ Tu \circ i , \\
& T \in C^{\infty}({\frak F}_i \mapsto {\frak B}_i),  
\quad  T (u) = Tu, \quad   u \in {\frak B}_i
 \end{align*}
 \end{equation}
 where $i$ is the complex conjugation on $S^2 = {\bf CP}^1$.
If $u$ satisfies $\bar{\partial}_J(u)=0$, then we say that 
$u$ is a {\em holomorphic curve} or $J$-curve.

Now let us   define the moduli space of holomorphic curves by:
$$   {\frak M}(\alpha, M_i,J_i)  = 
    \{ u \in C^{\infty}(S^2, M_i) \cap {\frak B}_i(\alpha) 
      : \bar{\partial}_{J} (u) =0   \}.$$
        There is an induced $S^1$ free action on ${\frak M}_i$.

\vspace{3mm}

We  say that $J$ is {\em regular},
 if for any $u \in {\frak M}_i$, the linealizations:
 $$D \bar{\partial}_{J}(u): T_u {\frak B}_i \mapsto ({\frak E}_i)_u$$
 are onto for all $i$.

\vspace{3mm}

The following estimates follow from
 the inverse function theorem and the Riemann Roch theorem:
\begin{prop}
Let $[(M_i,  \omega_i, J_i)]$ be a regular
almost Kaehler sequence.
Then the moduli spaces are $S^1$ free manifolds with dimension:
$$    \text{ dim } {\frak M}(\alpha, M_i, J_i)  =  2  <c_1(T^{1,0} M_i) , [u]>   -2 \text{ dim } M_i -1 .$$

If moreover it is  minimal, then  each ${\frak M}(\alpha, M_i,J_i )$ is compact. 
\end{prop}
Later on  we will omit to denote $\alpha$  explicitly.

\begin{defn}
The moduli spaces of holomorphic curves for  almost Kaehler sequences are given by:
$$  {\frak M}([(M_i,  \omega_i, J_i)]) = \cup_{i \geq 1}  \ 
        {\frak M}(\alpha, M_i,J_i).$$
\end{defn}

\begin{lem}
Let  $[(M_i,  \omega_i, J_i)]$ be a regular
almost Kaehler sequence.
Then  ${\frak M}([(M_i, \omega_i, J_i)])$     is a  $S^1$ freely
stratified manifold.
\end{lem}
{\em Example 2.1:}
Let us consider the standard holomorphic embedding
${\bf CP}^1 \hookrightarrow {\bf CP}^n$
with fixed two points. Modulo $S^1$ action, this is the unique element  
in the moduli space which is regular in the minimal class.
\vspace{3mm} \\
{\bf 2.B Sacks-Uhlenbeck's estimates:}
\begin{lem}
Let $[(M_i, \omega_i, J_i)]$ be a  minimal almost Kaehler sequence.
  For each $\alpha$, there is a constant $c(\alpha) \geq 0$ so that
any $u \in {\frak M}([(M_i, \omega_i, J_i)])$ satisfy  uniform estimates:
$$|\nabla^{\alpha} u|C^0(S^2)   \leq  c(\alpha).$$
\end{lem}
{\em Proof:}
We will only verify  the uniform estimate $|\nabla u|C^0 \leq c$.
The estimates for higher devrivatives 
follow from the elliptic regularity.

  There is a biholomorphic isomorphism:
$$\Phi: Z   = {\bf R} \times S^1    
 \cong {\bf CP}^1 \backslash  \{ 0, \infty \}, \quad (r,t) 
\to exp(r + 2\pi  it)$$
where we equip the standard complex structure on $Z$. 
For any holomorphic curve 
$u \in {\frak M}$, let us  regard $u: {\bf R} \times S^1 \mapsto M$ with
$u(- \infty) = p_0$ and $u(\infty) = p_{\infty} \in M_0$.

It follows from  $\frac{\partial }{\partial s} u + J(u) \frac{\partial }{\partial t}u=0$ 
that  the equalities hold:
$$||du||^2= \omega(\frac{\partial }{\partial s} u, J(u) \frac{\partial }{\partial s}u)
+ \omega(\frac{\partial }{\partial t} u, J(u) \frac{\partial }{\partial t}u)
=2 \omega(\frac{\partial }{\partial s} u , \frac{\partial }{\partial t}u)
= 2 ||u^*(\omega)||^2.$$

\begin{sublem}[SU]
There  are  constants $C$  and $\epsilon >0$
determined by $[(M_i, \omega_i, J_i)]$  so that 
for any holomorphic disk $u : D^2 \mapsto M = \cup_i M_i$ and
$E = \int_{D^2} u^*(\omega) \leq \epsilon$, 
the estimate holds:
$$ \varphi (x) \leq C E, \quad \varphi = |du|^2$$
for all $x \in D^2(\frac{1}{2})$.
 \end{sublem}

{\em Proof of the lemma:}
Let us fix a small positive constant $\delta >0$.
Then for any $u \in {\frak M}([(M_i, \omega_i, J_i)])$, we put
 $s(u) \equiv s_{\infty}(u)  - s_0(u)>0$, where:
\begin{align*}
& s_0(u) =  \sup \{ s \in {\bf R}:
\quad d(u(( - \infty , s) \times S^1) , p_0) \leq \delta \} , \\
& s_{\infty}(u) = \inf  \{ s \in {\bf R} :\quad
 d(u((s, \infty) \times S^1) , p_{\infty}) \leq \delta \}.
 \end{align*}

{\bf Step 1:}
We claim that for $0< \mu \leq    \frac{s(u)}{3}$, 
there is a positive $\epsilon >0$ determined by 
$[(M_i, \omega_i, J_i)]$   and $\mu$ with the estimates:
$$ \int_{ (- \infty, s_0(u)+ \mu] \times S^1} u^*(\omega) ,  \ \ 
 \int_{[s_{\infty}(u)-  \mu, \infty) \times S^1} u^*(\omega) \ \  \geq \epsilon.$$
We will only verify the first estimate. The latter follows by  the same argument.
Notice that  the translation on $Z$ is an automorphism
(which  does not preserve  the required condition 
$\int_{D(1)} u^*(\omega) = \frac{1}{2}<\omega, \alpha>$ on ${\frak B}$).

Let us choose a translation $T$ on $Z$ 
so that $s_0(u \circ T) = 0$ holds. Notice $s(u \circ T) =s(u) \geq   3 \mu$.
One may assume $s_0(u)=0$, since  the equality: 
$$ \int_{ (- \infty, s_0(u \circ T)+ \mu] \times S^1} (u \circ T)^*(\omega) =
 \int_{ (- \infty, s_0(u)+ \mu] \times S^1} u^*(\omega)$$
holds.
Let $D^2(b) \subset S^2$ be the disk with the radius $b$.
Let   $  a >0$ be as:
 $$ (- \infty, s_0(u)+ \mu] \times S^1 =  D^2(1+a) 
\backslash 0  \subset S^2$$ where we identify 
 $ (- \infty, s_0(u)] \times S^1 =  D^2(1) \backslash 0$. 
 We put  $D=D^2(1)$ and $D' = \ D^2(1+a)$.
 
   Let us put $B_{\delta}(0) \equiv \{ m \in M: d(p_0, m) < \delta) \} \subset M$
   as $\delta$ neighbourhood of $p_0$.
   Then 
$u(s,t) \in  \partial   B_{\delta}(0) $ and so the equality 
$d(u(s,t), u(-\infty)) = \delta $ holds 
  at $s= s_0(u)$ and some $t \in S^1$
   with $(s,t) \in \partial D$.

Suppose $\int_{D'} u^*(\omega) < \epsilon$ for 
sufficiently  small $\epsilon  =\epsilon(\mu)>0$.
 Then by sublemma $2.1$, 
  the uniform estimates of the derivative: 
 $$ |du|   \leq C(\mu) \sqrt{\epsilon} $$ hold on 
 all points of $D$. 
 This is a contradiction if $\epsilon >0$ are too small, since  $u(s,t) \in \partial B_{\delta}(0)$ 
 and $d(p_0, u(s,t)) = \delta$ as above.
 
This completes the proof of the claim.

\vspace{3mm}

{\bf Step 2:}
Let us  proceed the proof of the lemma by the contradiction argument.
Suppose contrary. Then there are families
$\{ u_i\}_i \subset {\frak M}([(M_i, \omega_i, J_i)])$
and $\{p_i\}_i \subset S^2$ with 
$ |\nabla u_i|(p_i) \to \infty$.
As [HV] $p611$, one may assume the inequalities:
\begin{equation}
\begin{align*}
&  |\nabla u_i|(x) \leq 2 |\nabla u_i|(p_i) \text{  for  all  }
x \text{ with }  d(x, p_i) \leq \epsilon_i,  \\
& \epsilon_i |\nabla u_i|(p_i) \to \infty, \quad\epsilon_i \to 0.
\end{align*}
\end{equation}
Let $D_i = D_i(p_i)$ be $r$ balls with the center $p_i$ for some small $r >0$.
 We put   rescaled balls with the center $p_i$ as
 $B_i = |\nabla u_i|(p_i) D_i(p_i)$
 by multiplying the numbers $|\nabla u_i|(p_i)$, 
 where one regards $B_i \subset {\bf C}$. 
 By conformal invariance, one gets a family of 
holomorphic maps $v_i : B_i \mapsto M = \cup_i M_i$.
This family satisfies uniform bounds:
 $$  |d v_i | (0) =1,     \quad 
|d v_i|(x)    \leq 2 \quad \text{ for }      |x| \leq   \epsilon_i |\nabla u_i|(p_i)   .$$
In particular by choosing  small $1 >> a, \epsilon' >0$, 
 $| |dv_i|^2(x)  - |dv_i|^2(0) | \leq  \epsilon'$
holds for all $x \in D(a)$ by the elliptic regularity,
where $D(a) \subset B_i$ is $a$ ball with the center $0$,
and $a$ is independent of $i$.
 This implies  the estimates
$|dv_i|(x) \geq \sqrt{1 -  \epsilon'}$.
In particular the uniform estimates hold from below:
$$\int_{D(b_i)} v_i^*(\omega) \geq \int_{D(a)} v_i^*(\omega) \geq C >0$$ 
for all $b_i \geq a$ with $D(b_i) \subset B_i$.

{\bf Step 3:}
On the other hand
 the uniform bounds 
$\int_{B_i} v_i^*(\omega) \leq m$ hold from above
where  $m$ is the minimal invariant.
We claim that there is 
 some family 
$ \{  R_i  \leq  \epsilon_i   |\nabla u_i|(p_i) \}$ with  $R_i \to \infty$
 such that the length $\delta_i$ of
$x_i \equiv v_i(exp(2 \pi i R_it)): S^1  \mapsto M$
must decay   $\delta_i \to 0$.

Notice that $[a,b_i] \times S^1 \subset B_i$ hold. Then there are some $R_i$ so that
the decay:
$$\int_{[R_i-1, R_i+1] \times S^1} v_i^*(\omega) \to 0$$
must hold. So the decay
$\sup_{x \in R_i \times S^1} |dv_i|(x) \to 0$ hold
by sublemma $2.1$, which verifies the claim.

{\bf Step 4:}
Thus  there is a family of small disks $\{ d_i \}_i \subset M$
which span $x_i$, and $\int_{d_i} \omega \to 0$.
Let $B_i' \subset B_i$ be $R_i$ balls with the center $p_i$,
whose  boundaries  are $x_i$.
Let us put two `almost'  holomorphic
spheres:
$$  u_i' =  
\begin{cases}
 & u_i
  \text{ on } S^2  \backslash B_i'  \\
  &   d_i
    \end{cases} ,
    \quad \quad 
  v_i' = B_i' \cup d_i.$$
By the condition, these must satisfy:
\begin{equation}
\begin{align*}
& <\omega, u_i'> + <\omega, v_i'> \to m >0, \\
& \text{ lim}_i <\omega, u_i'> , \quad 
 \text{ lim}_i <\omega, v_i'> \ \  \geq \ 0
\end{align*}
\end{equation} 
By minimality, one of $<\omega, u_i'>$ or $<\omega, v_i'>$ 
must be zero for all large $i$.
By step $2$ and  $3$,
  $<\omega, v_i'>$ must be positive and equal to $m$. So
  $<\omega, u_i'> =0$ must hold.

First of all, suppose  there is a uniform lower bound
$s(u_i) \geq 3\mu_0 >0$.
There are three cases;

(1) an infinite subset of $\{ p_i\}_i$ is contained in $(-\infty , s_0(u_i)] \times S^1$
or

(2)  is contained in  $[s_{\infty}(u_i), \infty) \times S^1$ or

(3) in $[s_0(u_i), s_{\infty}(u_i)] \times S^1$.
\vspace{3mm} \\
Suppose the case (1). Then by step $1$,
there is a positive $\epsilon >0$
with
$ \int_{[s_{\infty}(u_i) - \mu_0, \infty) \times S^1} (u_i')^*(\omega) \geq \epsilon$.
This implies the asymptotic bounds:
$$\text{ lim}_i <\omega, u_i'> \  \geq \  \epsilon$$
which give a contradiction.
The other cases can be considered similarly.

\vspace{3mm} 

{\bf Step 5:}
Let us verify that 
  $s(u_i) \to 0$ cannot happen.
  This will complete the proof of the lemma. 
Suppose contrary.
Let us take $p =0, q=\frac{1}{2} \in S^1$.
Then since $u_i(o \times s_*(u_i) ) \in B_{\delta}(*)$,
 $* = 0, \infty$ and $o = p,q$,
and since
 $d(B_{\delta}(0), B_{\delta}(\infty))  >0$ is positive,
 there are families
 $\{t_i \}$ and $\{r_i \}$, $t_i, r_i  \in [s_0(u_i), s_{\infty}(u_i)] $,
such that $|\nabla u_i|(p \times t_i ), |\nabla u_i|(q \times r_i) \to \infty$.
On the other hand one has a lower bound
 $d(p \times t_i, q \times r_i) \geq \frac{1}{2}$
in ${\bf R} \times S^1$.
By the same arguments as step $2$ and $3$, 
one obtains two non trivial almost holomorphic spheres,
which also cannot happen by  minimality of the homotopy class. 

 This completes the proof.
\vspace{3mm} \\
{\bf 2.C Hilbert completion of function spaces:}
Here  we introduce the basic function spaces 
for the infinite dimensional analysis.

Let   us take an element:
$$u \in {\frak B}_i(\alpha) \subset {\frak B}(\alpha)  \equiv \cup_{i \geq 1}  {\frak B}_i(\alpha)$$
 and let $U(u) \subset  {\frak B}(\alpha)$
 be  a small  neighbourhood of $u$ in the set of $L^2_{l+1}$
maps from $S^2$ to $M$.
Below we will describe its completion to a Hilbert manifold $\hat{U}(u)$.

Let us check the Sobolev embeddings for maps into Hilbert space.
\begin{lem}
There are constants $c_l$ with the uniform estimates:
$$||u||C^{l-1}(S^2)  \leq c_l ||u|||L^2_{l+1}(S^2).$$
\end{lem}
{\em Proof:}
By uniformity of complete local charts, 
it is enough to verify 
 the uniform estimates:
$$||u||C^{l-1}_c  \leq c_l ||u| |L^2_{l+1}$$
for $u \in C_c( U ; H ) $ with open subset $U \subset {\bf R}^2$.

The Sobolev estimate
$||\tilde{u}||C^{l-1}_c(U) \leq c_l ||\tilde{u}||L^2_{l+1}(U)$ hold 
for $\tilde{u} \in C_c(U)$.
Let $H $ be the closure of ${\bf R}^{\infty}$ with the standard norm,
and express $u =(\tilde{u}_0, \tilde{u}_1, \dots) \in C_c^{l-1}(U; H)$.
Then we have the estimates:
\begin{align*}
\Sigma_{k=0}^{l-1} 
|\nabla^k u|^2 (m) & = \Sigma_{k=0}^{l-1}
 \Sigma_{j \geq 0} |\nabla^k \tilde{u}_j|^2(m) \\
& \leq c_l  \Sigma_{j \geq 0} || \tilde{u}_j||^2L^2_{l+1}(U) = c_l ||u||^2 L^2_{l+1}(U)
\end{align*}
for any $m \in U$. By taking sup of the values of the left hand side,
we obtain the desired estimates.
This completes the proof.
\vspace{3mm} \\
{\em Remark 2.2:}
(1) Later 
we will find a reason why to use such completion of function  spaces,
rather than stratified spaces in some  step by step ways,
where  we will use  some automorphisms on almost Kaehler sequences
which do not preserve stratifications.

(2)  All functional spaces as Hilbert manifolds admit the
free and continuous $S^1$ actions. In precise there is a constant $C>0$
with the inequalities:
$$ C^{-1} ||u|| \leq \text{ Sup}_{t \in S^1} ||tu|| \leq C ||u||$$ for all elements $u$ 
in such  Hilbert spaces.

\vspace{3mm}

Let  $ \varphi(p) :D \equiv \cup_i D^{2i}(\epsilon)
     \hookrightarrow  \cup_i M_i$ be a
  complete almost Kaehler chart at $p$.
Sometimes we will identify $D$ and $D(p)$ where: 
$$D \equiv \cup_i D^{2i}(\epsilon)
\subset {\bf R}^{\infty} \subset H , \quad
D(p) \equiv   \varphi(p)( \cup_i D^{2i}(\epsilon)) \subset M.$$
By definition 
$D$ is equipped with the induced metric which is 
uniformly equivalent to the standard one on $H$.
Let us take the following data:
\vspace{2mm} 

(1)  finite set
of points $s_0, \dots, s_k \in S^2$, 

(2) an open cover $U_0, \dots, U_k$ with
$s_i \in U_i \subset S^2$, and 

(3) a partition of unity $f_0, \dots, f_k$
 over $S^2$.
 \vspace{3mm} \\
 For  $u \in {\frak B}(\alpha)$, 
one can choose  large $k$ so that 
each  image  $u(U_j)$ is contained in
a complete almost Kaehler chart at $\varphi(p_j)$ with $p_j = u(s_j)$.
Then one can express its restrictions as
$u|U_j : (U_j,s_j) \mapsto (D(p_j), p_j)$.
Identifying $D(p_j)$ with  $D$ as above, one may regard these maps as:
$$u|U_j: (U_j,s_j) \mapsto (D, 0) \subset (H,0).$$
Let us introduce precisely 
the Hilbert norm on the set of sections of $u^*(E(J))$ as
follows;
let us take any $\varphi \in \Gamma(u^*(E(J)))$. 
Then one may express
 the restriction  as:
$$\varphi|U_j:  TU_j \mapsto TD = D \times {\bf R}^{\infty}$$
which is anti linear with respect to $(i, J_{u( m)})$ at $(m,u(m))$.
Then one can define:
$$||\varphi||^2 L^2_l \equiv  \Sigma_{0 \leq j \leq k}  \Sigma_{ 0 \leq a \leq l}
\int_{U_j}  |\nabla^a(f_j \varphi|U_j)|^2 (m) dm.$$
By taking completion with respect to the above norm, 
one obtains the Hilbert bundles:
$$ \hat{\frak E} = \cup_{u \in {\frak B}} \hat{\frak E}_u \equiv
 \cup_{u \in {\frak B}} \{u\} \times L^2_l(u^*(E(J))).$$
One can  also make completion of  the functional spaces  ${\frak F}$  and  obtain
$\hat{\frak F}$.

Let $u \in {\frak B} = \cup_i {\frak B}_i$, and $U(u) \subset {\frak B}$
be a small  open subset.  
Let us complete $U(u)$ so that 
one obtains  a Hilbert manifold $\hat{U}(u)$ as below.
Let us write 
$u|U_j: (U_j,s_j) \mapsto (D, 0)$ for any $u \in {\frak B}_i$.
Then locally any element $v \in U(u) $
can be expressed as $v|U_l : U_l \mapsto {\bf R}^{\infty} \subset H$.
Then  we introduce Sobolev norms on $ U(u)$ by:
$$||v||^2 L^2_{l+1} = \Sigma_{0 \leq j \leq k} \Sigma_{0 \leq a \leq l+1}  
        \int_{U_j} |\nabla^a ( f_j v)|^2 (m) dm \quad (*)$$
By completion, one obtains  the Hilbert manifolds:
$$ \hat{U}(u) \quad \text{ over } \cup_i M_i \qquad 
 (u \in {\frak B}= \cup_{i \geq 1} {\frak B}_i)$$
 on   neighbourhoods  of $u$,
where  the local   Hilbert-structures are obtained 
passing through the exponential map.

Notice that if $u$ is holomorphic, then $k$ above can 
be chosen uniformly by lemma $2.2$.

\vspace{3mm}

Let us introduce: 
$$ \hat{\frak M}( [(M_i,  \omega_i, J_i)] ) = 
\cup_{u \in {\frak M}( [(M_i, \omega_i, J_i)] )}
   \ \  \{ v \in  \hat{U}(u)
      : \bar{\partial}_{J} (v) =0   \}.$$
Apriori this space is bigger than the moduli space
${\frak M}([(M_i,\omega_i,J_i)])$.
Later we study on their coincidence each other.
\vspace{3mm} \\
{\bf 2.C.2 Functional framework:}
Let $H$ be a Hilbert space, and $L \subset H$ be a closed subspace.
\begin{lem}
Let $F: H \to H$ be a bounded operator with closed range,
whose kernel consists of finite dimensional subspace.
Then $F(L) \subset H$ is also closed.
In particular if $F$ is injective, then $F(L)$ is closed.
\end{lem}
{\em Proof:}
If kernel $F=0$, then $F: H \cong F(H)$ gives an isomorphism.
In particular $F(L)$ is closed.

Suppose ker$(F) = K \subset H$ is of finite dimension.
Then $F$ induces an isomorphism
$F: H/ K \cong F(H)$, where we equip with the metric 
on $H/K$ by use of orthogonal decomposition $H = K^{\perp} \oplus K$.
Then it is enough to see that the image of the projection 
$pr(L) \subset H/K$ is still closed.

One may assume that $L \cap K= 0$ by replacing $L$ by $(L \cap K)^{\perp}$  in $L$,
when it  has positive dimension.

Suppose a sequence $\{ \bar{v}_i\}_i \subset pr(L)$ converge to some element
$\bar{v} \in H/K$. 
By the assumption, the representatives $v_i \in L$ of $\bar{v}_i$ are unique.
Let us represent  $v_i  = v_i^1+v_i^2 \in L$ with 
 respect to the decomposition $H=K^{\perp}\oplus K$.

We claim that $||v_i||$ are uniformly bounded. 
Suppose contrary and assume $||v_i|| \to \infty$.
Then by normalizing as $w_i = ||v_i||^{-1}v_i = w_i^1+w_i^2$,
both $||w_i^1|| \to 0$ and $||w_i^2|| \to 1$ hold.
Since $K$ is finite dimensional and $L$ is closed,
a subsequence $w_i$ converges to some element $w \in L \cap K$ with $||w||=1$.
This contradicts to  our assumption, which verifies the claim.

Now since $\{v_i^2\}_i \subset K$ is a bounded sequence,
a subsequence converges to some element $v^2 \in K$.
Since $v_i^1$ converges to $v$, it follows from these 
that a subsequence of  $\{v_i\}_i$ converges to $v+ v^2 \in L$.
This implies $\bar{v} \in pr(L)$.

This completes the proof.
\vspace{3mm} \\
{\em Remark 2.3:}
The assumption of finite dimensionality is necessary.
Let $H$  be a separable  infinite dimensional Hilbert space,
and choose an orthonormal basis $\{v_i\}_i$.
Let $0< a_i  \to 0$
be a decreasing family of numbers.

Let us consider the surjective bounded map:
$$F = \text{ id } \oplus 0 : H \oplus H \to H$$
and the closed subspaces $L $
spanned by the basis:
$$L = \text{ span } \{ w_i = (a_i v_i, v_i): i =0,1,2, \dots\} \subset   H \oplus H.$$
We claim that the image of the restriction $F|L$ is not closed.
Suppose  contrary. Then since $F|L$ is injective, the restriction must be
an isomorphism by the open mapping theorem.
So there must exist some $C>0$ with the uniform estimates:
$$|a_iv_i| = |F(w_i)| \geq C|w_i| = C \sqrt{a_i^2+1}. $$
But the left hand side converge to $0$,
which contradict to the right hand side.

\vspace{3mm}

The following abstract property is the key to our Fredholm theory 
we develop later:
\begin{cor}
Suppose the above situation, and choose another Hilbert space $W$.
Then the image of the Hilbert space tensor product $L \otimes W$
over  the induced operator $F \otimes 1: H \otimes W \to H \otimes W$,
 still has closed range.

In particular if $F$ is an isomorphism, then $F \otimes 1$ is also the same.
\end{cor}
{\em Proof:}
 Let us put $E= L \cap \text{ Ker } (F)$, and decompose 
 $L \cong L' \oplus E$. Then  $F(L)=F(L')$ holds.
 Since the restriction $F|L'$ is injective, 
 it gives the isomorphism onto $F(L)$
 by the open mapping theorem.
 
Since  the restriction $F \otimes 1|L' \otimes W$ gives the isomorphism 
onto $F \otimes 1(L' \otimes W)= F\otimes 1(L \otimes W)$,
the conclusion follows.
This completes the proof.
\vspace{3mm} \\
{\bf 2.D Geometric conditions:}
We study  functional analytic properties of the Cauchy-Riemann operators 
over almost Kaehler  sequences which satisfy 
the  geometric conditions we have introduced in $1.B$.

Let us say that a subset
$S \subset {\frak B}([(M_i, \omega_i, J_i)]) = \cup_i {\frak B}(M_i, \omega_i, J_i)$
is {\em bounded}, if there is some $i_0$ so that the images of $u$ are contained in $M_{i_0}$ 
for any elements  $u \in S$. 

Our aim in $2.D$  is to verify the following:

 \begin{thm}
 Let $[(M_i, \omega_i, J_i)]$ be a symmetric  Kaehler sequence.
 
 (1)
Suppose it is  regular and dim $\cup_i $ Ker $ D_u \bar{\partial}_i =N$ is finite,
 then it is in fact strongly regular of index $N$.

 In particluar  
  $ {\frak M}([(M_i, \omega_i , J_i)]) $ is a regular $N$ dimensional $S^1$ free manifold.
  
  (2) Assume  moreover  it is  isotropic, and each connected component of
   $  {\frak M}([(M_i, \omega_i , J_i)])$ is bounded. Then 
   the equality holds:
   $$ {\frak M}([(M_i, \omega_i , J_i)]) =  {\frak M}( M_0,\omega_0, J_0) .$$
   
   In particular if it is minimal, then  $ {\frak M}([(M_i, \omega_i , J_i)])    $  is compact.
  \end{thm}
 The second condition in $ (2) $ is satisfied when $N=1$.
 
 \vspace{2mm}

The former  follows from combination of lemma $2.7$ and 
proposition $2.2$ below.
The latter is verified in lemma $2.8$.
\vspace{3mm} \\
{\bf 2.D.1 Strong regularity over symmetric Kaehler  sequences:}
Let $[(M_i, \omega_i, J_i)]$ be a symmetric
  almost Kaehler sequence, and 
  choose symmetric  data  $\{ (\pi_k, P_i)\}_{i,k}$ with respect to $(M_k, M_{l'})$
  for some $l'=l(k)$.

For any $u \in {\frak B}_k \subset {\frak B}_{l'} \subset {\frak B}$,
let $\hat{U}(u) \subset {\frak B}$ be as in $2.C$.
Let us take small neighbourhoods 
  $U(u)_{l'}  \subset {\frak B}_{l'} \cap \hat{U}(u)$.
  
  There are extended projections:
  $$\bar{\pi}_k : \hat{U}(u) \to U(u)_k$$
  with $\bar{\pi}_k|U(u)_k =$ id.
  Then
 the isomorphism:
 $$T_u U(u)_{l'} \cong  T_u {\frak B}_k \oplus V_u(k,l')$$
hold, where:
$$V_u(k,l')= \text{ Ker } (\bar{\pi}_k)_*  \cap T_uU(u)_{l'} .$$

\begin{lem}
The complete isomorphisms hold:
$$T_u \hat{U}(u) \cong T_u {\frak B}_k \oplus V_u(k,l') \otimes H$$
where $H$ is a  separable Hilbert space.
\end{lem}
{\em Proof:}
This follows from  lemma $1.6$.
This completes the  proof.

\vspace{3mm}

The Cauchy-Riemann operator $\bar{\partial}_J$ and  the tangent map $T$
 give smooth sections respectively:
$$  \bar{\partial}_J : \hat{U}(u) \mapsto \hat{\frak E}|\hat{U}(u), \quad
  T:  \hat{U}(u) \mapsto  \hat{\frak F}|\hat{U}(u).$$

\begin{defn} Let  $[(M_i,  \omega_i, J_i)]$ be
 a regular  almost Kaehler sequence.
It is strongly regular,  if the differential:
$$D\bar{\partial}_u: 
T_u \hat{U}(u) \mapsto \hat{ \frak E}_u | \hat{U}(u)$$ is surjective
 for  any $ u \in {\frak M}([(M_i, \omega_i, J_i)])$.
\end{defn}

\begin{lem}  Let $[(M_i, \omega_i, J_i)]$ be a symmetric 
   Kaehler sequence. 
   
   Then
  $D\bar{\partial}_J  : T_u \hat{U}(u) \mapsto \hat{\frak E}_u$ 
 has  closed range. 
\end{lem}
{\em Proof:} 
Firstly we verify that 
 $DT:  T_u \hat{U}(u) \mapsto \hat{\frak F}_u$  has closed range, and then
 we verify the conclusion.

{\bf Step 1:}
 Let us take some $k$ so that 
$u \in   {\frak M}(M_k, \omega_k, J_k)$ with $u : S^2 \to M_k$.
Let $P_i: \{(M, M_k), \omega, J\}  \cong \{(M,M_k), \omega, J\}$ be the locally homogeneous  data
 for all $i \geq  l' =l(k)$.
 
 Let us consider the bundles over $S^2 \times M_k$:
 $$N_{k,l'}(m,z)= \{ \phi: T_z S^2 \to \text{ Ker }( \pi_k )_*
 \cap T_m M_{l'}: \text{ linear } \} \subset F_l(z,m)$$
  and put the Hilbert subbundles over ${\frak B}_k$:
 $${\frak F}_{k,l'}= L^2_l ({\frak B}_k^*(N_{k,l'}))= \cup_{u \in {\frak B}_k} \{u\} \times L^2_l(u^*(N_{k,l'})) \subset {\frak F}_{l'}| {\frak B}_k.$$
 
 There are
   bundle decompositions  
   ${\frak F}_{l'} |{\frak B}_k \cong {\frak F}_k \oplus {\frak F}_{k,l'}$ 
  over ${\frak B}_k$ given by:
  $$\bar{\phi} \to ((\pi_k)_*(\bar{\phi}),\bar{\phi} -  (\pi_k)_*(\bar{\phi})).$$
 It follows from  symmetric property that  
 the bundle decompositions:
 $$\hat{\frak F}|{\frak B}_k \cong {\frak F}_k \oplus {\frak F}_{k,l'} \otimes H$$ 
 hold  as lemma $2.5$.
 Let:
   $$DT_{l'}: T_u {\frak B}_{l'} = T_u {\frak B}_k \oplus V_u(k,l')
 \to {\frak F}_{l'}|_u = {\frak F}_k|_u \oplus {\frak F}_{k,l'}|_u$$ be the tangent map.
 Clearly this is diagonal $DT_{l'}= DT_k \oplus DT_{k,l'}$ with respect to these decompositions.
Then the  total tangent map:
 $$DT: T_u \hat{U}(u) \cong T_u {\frak B}_k \oplus V_u(k,l') \otimes H
\to  \hat{\frak F}|_u \cong {\frak F}_k|_u  \oplus {\frak F}_{k,l'}|_u  \otimes H$$ 
is also  diagonal:
$$DT= DT_k  \oplus (DT_{k,l'} \otimes \text{ id }).$$

Now 
 $DT_l: T_u {\frak B}_l \to {\frak F}_l|_u$ has closed range with finite dimensional  kernel,
 which follow from the well known analysis of holomorphic curves 
 into finite dimensional symplectic manifolds (see [HV]).
 Since 
 $V_u(k,l) \subset T_u{\frak B}_l$ are closed subspaces, 
it follows from   corollary $2.1$ that 
 $DT_{k,l'}  \otimes \text{ id }$ has closed range.
Since $DT_k$ has closed range, the direct sum
 $DT_k \oplus (DT_{k,l'} \otimes \text{ id })$ also has closed range.

{\bf Step 2:} 
 One can follow the above proof by changing $T$ by $\bar{\partial}$ by use of the 
complete Kaehler charts. 
Notice the formula:
$$ D\bar{\partial}_J(v) = (DT + J \circ DT \circ i)(v) +  N(v) $$
where $N$ involves $\nabla J$, 
and  $ N \equiv 0$ when $J$ is integrable.
So if it  is Kaehler, then
$   D\bar{\partial}_J   = DT + J \circ DT \circ i  $ holds.
In particular if we decompose the holomorphic local charts as in step $1$,
then $D\bar{\partial}_J $ can be also expressed 
as the form $K_1 \oplus (K_2 \otimes \text{ id })$.
The rest of the argument is parallel to  step $1$.

This completes the proof.
\vspace{3mm} \\
{\bf 2.D.2 Index computations:}
Let $[(M_i, \omega_i, J_i)]$ be a symmetric almost Kaehler sequence.
Thus
there are holomorphic projections
$\pi_k: U_{\epsilon}(M_k) \to M_k$ with $\pi_k|M_k=$ id from 
  small neigbourhoods  in $ M= \cup_i M_i$.

For  $u \in {\frak B}_k$, 
let $\bar{\pi}_j: \hat{U}(u) \mapsto {\frak B}_j$ be 
the  induced   projections for all $j \geq k$.

\begin{lem} Let
 $[(M_i, \omega_i, J_i)]$ be an almost Kaehler sequence.
If    $\cup_i $ Ker$D_u\bar{\partial}_i$ is of finite dimension,
   then the equality holds:
$$   \text{ Ker } D_u\bar{\partial}_J    = \cup_i   \text{ Ker } D_u\bar{\partial}_i.$$

In particular the left hand side is of finite dimension.
\end{lem}
$Proof:$ 
The condition implies    
$\cup_i $ Ker$D_u\bar{\partial}_i = $
 Ker$D_u\bar{\partial}_{i_0}$ 
for some $i_0$.

Suppose contrary and assume  Ker $D_u\bar{\partial}_J    \ne \cup_i   $ Ker $D_u\bar{\partial}_i$.
Let $u_t \subset \hat{U}(u)$ be a smooth curve with $u_0 =u$ and 
$u_t'|_{t=0} \equiv v \in $ Ker  $D_u\bar{\partial}_J  $
but $v \not\in  \cup_i   $ Ker $D_u\bar{\partial}_i$.
  Then
 $D_u\bar{\partial}_j((\pi_j)_*(v)) = (\pi_j)_*(D_u\bar{\partial}(v))  =0$ vanish 
  for  all $j \geq k$.
So $(\pi_j)_*(v) $ lies in Ker$D_u \bar{\partial}_j$. 
It  must be contained in  Ker$D_u\bar{\partial}_{i_0}$
by the assumption.
 Since $j$ is arbitrary, this implies $v \in $ Ker$D_u\bar{\partial}_{i_0}$. 
 
This completes the proof.

\vspace{3mm}

Let us denote 
$\bar{\partial}_J : \hat{U}(u) \mapsto \hat{\frak E}|\hat{U}(u)$ and
$\bar{\partial}_i  :      {\frak B}_i   \mapsto {\frak E}_i $ respectively.

\begin{prop}
Let $[(M_i, \omega_i, J_i)]$ be
 a symmetric  Kaehler sequence.
  Let us  choose any $u \in {\frak M}(M_k, \omega_k, J_k)$.
  
If the uniform bounds  dim Coker $D_u\bar{\partial}_i \leq M$ hold  for all $i\geq k$,
then  dim Coker $D_u \bar{\partial}_J \leq M$ also holds.

In particular if it is regular, then it is in fact strongly regular.
\end{prop}
$Proof:$
$D\bar{\partial}_J$ has closed range by  lemma $2.6$.
Suppose the estimates dim Coker $D_u\bar{\partial}_J \geq M+1$ could hold, 
and take orthogonal unit elements  $u_1, \dots, u_{M+1}$
in Coker $D_u\bar{\partial}_J$.

There are large $l >>k$ so that  $u_i^l = (\pi_l)_*(u_i)$ are  defined for all $1 \leq i \leq M+1$.
For  small $\epsilon >0$, let us choose sufficiently large  $l$ so that  the estimates below
hold, where $B \subset$ im$D_u \bar{\partial}_l \subset  ({\frak E}_l)_u$ 
are the unit balls:
$$|u_i^l|^2 \geq 1 - \epsilon, \quad |< u_i^l, u_j^l>| \leq \epsilon,
\quad | < B , u_i^l>| \leq \epsilon.$$

There are numbers $a_1, \dots, a_{M+1} \in {\bf R}$ with   $\Sigma_j |a_j|^2 =1$
so that  $v \equiv \Sigma_i a_i u_i^l $ lie in  im $D_u\bar{\partial}_l$, 
since  dim Coker $D_u\bar{\partial}_l \leq M$ hold.
Let us pick up $i$ with $|a_i| = $ Sup$_{1 \leq j \leq M+1} |a_j| \geq \frac{1}{\sqrt{M+1}}$.
Then one should have the estimates:
$$\epsilon \geq |<v,   u_i^l>| \geq |a_i| (1- \epsilon) -  \epsilon \Sigma_{i \ne j} |a_j|
\geq |a_i|(1- \epsilon) - \sqrt{M} \epsilon.$$
Since $\epsilon$ are arbitrarily small, this  is a contradiction.

This completes the proof.
\vspace{3mm} \\
 {\em Example 2.2:}
  ${\bf C}P^{\infty}$  is strongly regular of index $1$ by propopsition $2.2$.
  \vspace{3mm}

So for regular and symmetric  Kaehler sequences, 
the moduli spaces of  holomorphic curves are strongly regular
 with the expected indices.
      With respect to these obsevations, we would like to propose the following:
  \vspace{2mm} \\
{\em Conjecture 2.1:} Let $[(M_i, \omega_i, J_i)]$ be 
a symmetric Kaehler sequence.

(1) 
 One can perturbe the complex structure (to be almost Kaehler)
 so that the result  could become strongly regular.
 
 (2) For irregular case,  index $ D\bar{\partial}_J = M$ hold  when
 index $D\bar{\partial}_i =M $  for all $i$  and 
  $\cup_i $ Ker$D\bar{\partial}_i$  is of finite dimension.

 (3) Suppose moreover it is regular (and hence strongly regular).
 Then 
the embedding:
$$ {\frak M}( [(M_i,  \omega_i, J_i)] )  = \cup_i {\frak M}_i(M_i,  \omega_i, J_i)
\subset \hat{\frak M}( [(M_i,  \omega_i, J_i)] )$$
is in fact equality.
\vspace{3mm} \\
 {\bf 2.D.3 Compactness of moduli spaces:}
 Let  $[(M_i,  \omega_i, J_i)]$ be an isotropic  symmetric  almost Kaehler sequence.

   \begin{lem}
Suppose that 
each connected component of 
$ {\frak M}( [(M_i,  \omega_i, J_i)] )  $ is bounded.
 Then  the equality holds:
 $$ {\frak M}( [(M_i,  \omega_i, J_i)] )   = {\frak M}( M_0,  \omega_0, J_0) .$$
 
In particular if it
 is minimal, then  ${\frak M}( [(M_i,  \omega_i, J_i)] )  $ is compact.
 \end{lem}
 Notice that the condition is satisfied for regular almost Kaehler sequences 
 whose moduli spaces have $1$ dimensional.
 
 \vspace{3mm}

 {\em Proof:}
 Let us choose an element $[u] \in  {\frak M}( [(M_i,  \omega_i, J_i)] )  $.
  By the assumption there is some $l_0$ so that
    the connected component ${\frak M}(u)$
  containing $u$ has all their images  in $M_{l_0}$.
  
Suppose  there could exist  some   $k$ with  its symmetry over $(M_k, M_{l'})$ for  $l'=l(k)>k$,
  such that  ${\frak M}(u)$  has all their images  in $M_{l'}$ but not all  in $M_k$.
  
  Let   $P_i^t: (M, W_i,M_k) \cong (M, M_{l'},M_k)$ be the isotropies for $i \geq l'$,
  where $M=\cup_i M_i$. There is some $u' \in {\frak M}(u)$ so that
   the images of $P^0_i(u') =u'$ are contained in $M_{l'}$, 
  while $(P^1_i)^{-1}(u')$ are not the case for all  $i > l'$.
  This implies that the images of  ${\frak M}(u)$ cannot be contained in $M_{l'}$,
  since $(P^1_i)^{-1}(u')$ must be contained in ${\frak M}(u)$.
  This  contradicts to the assumption.
  So ${\frak M}(u)$ must be contained in $M_k$.
  
  Next let us replace the pair $(k,l'(k))$ by $(k-1, l'(k-1))$.
  Because of the relation $l'(k-1) >k-1$, the inequality
  $l'(k-1) \geq k$ must hold.
  In particular the images of  ${\frak M}(u)$ are contained in $M_{l'(k-1)}$,
  and proceed the same argument. Then we find that 
  ${\frak M}(u)$ is contained in $M_{k-1}$.
  
  Let us continue this process. Then finally we find that the images of
  ${\frak M}(u)$ must be contained in $M_0$.
   This completes the proof.

\vspace{3mm}

It would be interesting to study 
what happens for the cases of  positive dimension.
We would like to  propose the following:
\vspace{3mm} \\
 {\em Conjecture 2.2:}
 Suppose $ {\frak M}( [(M_i,  \omega_i, J_i)] )  $ is a smooth manifold of finite dimension.
 Then the equality:
 $$ {\frak M}( [(M_i,  \omega_i, J_i)] )   = {\frak M}( M_k,  \omega_k, J_k) $$
 holds for some $k$.

\vspace{3mm}

 Notice that the strong regularity condition is   stable under small perturbations,
  and we expect that such property can be studied over deformations
  of these sequences. This is the topic at the next section.

 \vspace{3mm}

\section{Infinitesimal neighbourhoods of almost Kaehler sequences}
One of the quite characteristic properties which 
infinite dimensional spaces possess, is that they can contain 
many spaces as their proper subsets  which can  include even themselves.
In our formulation, infinite dimensional spaces are consisted by 
sequences of finite dimensional spaces, which leads us canonically 
to introduce {\em neighbourhoods} of such sequences as some sets of 
another sequences.
In the  finite dimensional case, notion of neighbourhoods will require some dimension restrictions.
However such limitations are free for our infinite dimensional sequences.

In later sections, we apply such notions to study {\em stability} of
the invariants of almost Kaehler sequences
under ``very small perturbations'', which measure
continuity of these invariants  in the framework of the neighbourhoods.
\vspace{3mm} \\
{\bf 3.A Convergence:}
Let us start from the  finite dimensional case.
Let $\{ X_i \}_{i=1}^{\infty}$ be a family of smooth  manifolds of the same dimension
embedded into another finite dimensional smooth manifold $M$.

We say that the set $\{ X_i \}$ {\em converges } to
 $X$ in $M$, if
there exist  coverings of $X = \cup_l U_l$ and $X_i = \cup_l U_l^i$
so that for all sufficiently large $i >>0$, there are  diffeomorphisms
$F_l^i : U_l \cong U_l^i (\subset M)$ which  converge to the identity in $C^{\infty}$ topology in $M$.
Notice that if $X_i$ are sufficiently near $X$ in $C^{\infty}$ topology, then 
they are  isotopic to $X$.
\vspace{3mm} \\
{\bf 3.A.1 Convergence of almost Kaehler sequences:}
Let us introduce the following:
\begin{defn}
Let $[(M_i, \omega_i, J_i)  ]$ be an almost Kaehler sequence.
Let us say that a  family of almost Kaehler sequences 
$\{ [(M_i^l, \omega_i^l, J_i^l) ] \}_{l=1}^{\infty}$ {\em converges } to
$[(M_i, \omega_i, J_i)  ]$, 
if there are positive $\epsilon >0$, 
  subindices $\{ k(i)\}_i$ with $k(i) \geq i$ and compatible
 $C^{\infty}$ embeddings of 
 $[(M_i^l, \omega_i^l, J_i^l) ] $ into $M = \cup_k M_k$ for all $l$:
 $$ I(i,l): M_i^l \hookrightarrow M_{k(i)}, \quad I(i+1,l)|M_i^l =I(i,l)$$
 with   almost complex $J(i,l)$ and  symplectic $\omega(i,l)$
 structures  on the  open 
$\epsilon$  tublar  neighbourhoods 
$I(i,l) (M_i^l) \subset U(i,l) \subset M $, which extend the given ones:
$$J(i,l)|I(i,l)(M_i^l) = J_i^l, \quad \omega(i,l)|I(i,l)(M_i^l)=\omega_i^l$$
 so that the following three conditions hold:
 
(1)  For each $l$, 
there is a uniformly bounded covering
 $\{(p, \varphi(p))\}$ on $[(M_i^l, \omega_i^l, J_i^l)]$,
and its extension $\{(p, \psi(p))\}$ 
by 
 $\epsilon$
complete almost Kaehler charts  over $( U(i,l), \omega(i,l), J(i,l))$
for all
 $p \in I(i,l)(M_i^l)$:
\begin{equation}
\begin{align*}
&   \varphi(p) :  \cup_{s \geq 1} D^{2s}(\epsilon)
\hookrightarrow \text{ im } \varphi(p) \subset  \cup_{r \geq 1} M_r^l \\ 
& \qquad \qquad  \qquad  \cap   \qquad \quad  \cap   \\
& \psi(p):   \cup_{s \geq 1} D^{2s}(\epsilon)  
  \hookrightarrow  U(i,l) \subset \quad  \cup_{r \geq 1} M_r
\end{align*}
\end{equation}

 (2)
There are     families of holomorphic  maps:
$$\pi_i^l :  U(i,l)  \mapsto M_i^l$$
with respect to $\omega(i,l)$, which satisfy the following properties:
$$ \pi_i^l  |  M_i^l  = \text{ id }, \quad
  \pi_i^l(\psi(p)(m) ) = \psi(p)(\pi_i^l(m))$$
  for all  
$m \in U(i,l)$, where $\pi_i^l :   \cup_{s \geq 1} D^{2s}(\epsilon)  \to D^{2d^l_i}(\epsilon)$
are the projections with $d^l_i =$ dim $M_i^l$.

\vspace{2mm}

(3)  All derivatives of the operators converge to zero as $l \to \infty$:
$$ \sup_i||\nabla^{\alpha}( J(i,l) - J)||_g U(i,l), \quad
\sup_i||\nabla^{\alpha}(\omega(i,l) - \omega)||_g U(i,l) \to 0.$$

\vspace{2mm}

 (4)  For each $i$, $\{ I(i,l)(M_i^l) \}_{l \geq 1} \subset M_{k(i)}$ 
converges to $M_i$ in $M_{k(i)}$ in $C^{\infty}$.
\end{defn}
Notice that in general $I(i,l)(M^l_i)$ may be far away from $M_i$
as $i,l \to \infty$ (convergence in (3) is not assumed  uniform with respect to $i$).

In short, we will denote convergence by the notation:
$$\{ [(M_i^l, \omega_i^l, J_i^l) ] \}_l \to [(M_i, \omega_i, J_i)].$$
We say that
$\{ [(M_i^l, \omega_i^l, J_i^l) ] \}_l \to [(M_i, \omega_i, J_i)]$
is {\em minimal convergence}, if all almost Kaehler sequences
are minimal.
\vspace{3mm} \\
{\em Examples 3.1:}  (1)
Let us choose a uniform family of almost Kaehler manifolds
  $(X_i, \omega_i, J_i)$, $i= 1, 2, \dots$, and 
  take another family of almost Kaehler manifolds
  $(X_0, \omega_0^l , J_0^l)$ so that both of $\omega_0^l $ and $  J_0^l$
  converge to $\omega_0$ and $J_0$ in $C^{\infty}$ as $l \to \infty$ respectively. 
    Let us put:
  $$(M_i^l , \omega_i^l , J_i^l) \equiv (X_0 \times \dots \times 
  X_i ,   \omega_0^l + \omega_1 + \dots + \omega_i,
   J_0^l \oplus J_1 \oplus \dots , \oplus J_i).$$
  Then $\{[(M_i^l , \omega_i^l , J_i^l)]\}$ is  convergent to the product  with $I(i,l)=$ id.

(2) Let $(X, \omega, J)$ be an almost Kaehler manifold with a fixed point $x_0 \in X$.
Let $\omega(i) = \omega + \dots + \omega$ and $ J(i) = J \oplus
\dots  \oplus J$ be almost Kaehler data on       $X \times \dots \times X$.
We put $M_i = \times^i X$ where $M_i = \times^i X \times \{ x_0 \} \subset M_{i+1}$.

Let $a: {\bf N} \mapsto {\bf N}$ be any proper function.
Let us choose all the same  $M_i^l =(M_i, \omega(i), J(i))$,
but take embeddings $I(i,l): M_i^l \hookrightarrow M_{i+1}$ as:
$$I(i,l)((m_0, \dots, m_i))= (m_1, \dots, m_{a(l)}, x_0, m_{a(l)+1}, \dots, m_i).$$
Then 
   $\{[(M_i^l , \omega(i) , J(i))]\}$ is a convergent family.
  \vspace{3mm} \\
 {\bf 3.B Infinitesimal neighbourhoods:}
Let us take an  almost Kaehler sequence
$ [(M_i, \omega_i, J_i)]$, and 
two  familes of  almost Kaehler sequences:
$$ \{[(M_i^l, \omega_i^l, J_i^l)]\}_{l \geq 1} , \  \{[(N_i^l, \tau^l_i, I^l_i)]\}_{l \geq 1}
\ \ \to \ \    [(M_i, \omega_i, J_i)]$$
which both converge to the same almost Kaehler sequence.

Let us say that the two families 
are {\em equivalent}, if 
there are infinite subindices $\{k(l)\}_l$  and $\{k'(l)\}_l$ 
such that there are isomorphisms between almost Kaehler sequences   for all $l$:
$$ [(M_i^{k(l)}, \omega_i^{k(l)}, J_i^{k(l)} ]_{i \geq 1} \cong  [(N_i^{k'(l)}, \tau^{k'(l)}_i, I^{k'(l)}_i)]_{i \geq 1}.$$
We denote 
  the equivalence class  of $ \{[(M_i^l, \omega_i^l, J_i^l)]\}_{l \geq 1}$ by
  $[[(M_i^l, \omega_i^l, J_i^l)]]$.

\begin{defn} An infinitesimal neighbourhood of  $[(M_i, \omega_i, J_i)]$
 is given by  the set of equivalece classes by  convergent families:
$$ {\frak N}([(M_i, \omega_i, J_i)]) = \{ [[(M_i^l, \omega_i^l, J_i^l)]] : 
 \{ [(M_i^l, \omega_i^l, J_i^l)] \}_{l \geq 1} \to [(M_i, \omega_i, J_i)] \} .$$
\end{defn}

Let  $[(M_i, \omega_i, J_i)]$ and 
$\{[(M_i^l, \omega_i^l, J_i^l)]\}_{l \geq 1}  $ be 
 minimal  almost Kaehler sequences, and suppose the family converges to 
  $ [(M_i, \omega_i, J_i)]$ with the data $\{I(i,l): M_i^l \hookrightarrow M_{k(i)} \}_{i,l}$.
 We  say the convergence {\em keeps} minimal classes, 
 if $I(l)_*$ maps the minimal classes to the one over 
   $ [(M_i, \omega_i, J_i)]$
   on $\pi_2$.

      Now by restricting on minimal almost Kaehler sequences, we define
      the infinitesimal neighbourhoods of minimal almost Kaehler sequences:
\begin{align*}
  {\frak N  }^m   (&[(M_i,   \omega_i, J_i)]) = \{ [[(M_i^l, \omega_i^l, J_i^l)]] :  
 \{ [(M_i^l, \omega_i^l, J_i^l)] \}_{l \geq 1} \to [(M_i, \omega_i, J_i)]  \\
&    \{ [(M_i^l, \omega_i^l, J_i^l)] \}_{l \geq 1} : \text{ minimal families 
which   keep minimal classes  }  \} .
\end{align*}
{\vspace{3mm} \\
{\bf 3.C Moduli theory over perturbations of spaces:}
 So far we have studied global analysis  of the moduli spaces of holomorphic curves into 
 infinite dimensional spaces. 
 It turned out that their behaviours are well controlled if we assume
 integrability with high symmetry, where
 they are  stable
 under small perturbations of their structures
   under such situations.
  This is the key aspect which allows us to study 
  moduli theory over  infinitesimal neighbourhoods.

 Let $[(M_i,  \omega_i, J_i)]$ be an almost Kaehler sequence.
 Let us consider its infinitesimal neighbourhood and  take an element:
$$ [[(M_i^l,  \omega_i^l, J_i^l)]]  \in  {\frak N}^m([(M_i,  \omega_i, J_i)]).$$
The following holds since regularity condition is open:

\begin{lem}
Let    $[(M_i,  \omega_i, J_i)]$  be a  minimal and regular
almost Kaehler sequence.
Suppose ${\frak M}([(M_i,  \omega_i, J_i)])$ is bounded
and $S^1$ freely $1$ dimensional manifold.
Then there is  $l_0>>0$ so that 
there are $S^1$ compatible embeddings for all $l \geq l_0$:
$${\frak M}([(M_i,  \omega_i, J_i)]) \subset {\frak M}([(M_i^l,  \omega_i^l, J_i^l)]).$$
\end{lem}
{\em Proof:}
It is enough to use well known analysis which are applied for
finite dimensional almost Kaehler manifolds.
By the assumption, there is $i$ so that 
${\frak M}([(M_i,  \omega_i, J_i)]) = {\frak M}(M_i,  \omega_i, J_i)$
are both compact.
 Let us consider a family of embeddings
$I(i,l): M_i^l \hookrightarrow M_k$ in definition  $3.1$.

By compactness there is a large $l_0$ so that for all $l \geq l_0$, 
the images of any $u \in {\frak M}(M_i,  \omega_i, J_i)$ are contained in
$U(i,l) \subset \cup_j M_j$. 
Then by use of the projections $\pi_i^l : U(i,l) \to M_i^l$,
let us consider a smooth maps $u' \equiv \pi_i^l \circ u : {\bf CP}^1 \to M_i^l$.
The differentials of the Cauchy-Riemann operators must be surjective
at $u'$ since $(M_i^l, \omega_i^l, J_i^l)$ converge to $(M_i, \omega_i, J_i)$  smoothly.
Let us apply the infinite dimensional implicit function theorem to $u'$ to 
obtain a holomorphic curve $u''$ with respect to $(M_i^l,\omega_i^l, J_i^l)$.
This assignment $u' \to u''$ extends to the $S^1$ freely equivariant ones.
This completes the proof.

\vspace{3mm}

Notice that we have used compactness of the moduli spaces 
so that the  inverse  of the differentials of C-R operators 
satisfy     uniform estimates from above.

 In order to obtain the converse embedding, we use
the strong regularity condition below.

For any map $u: {\bf C}P^1 \to M_i^l$, let us
consider the compositions:
 $$v \equiv I(i,l) \circ u : {\bf C}P^1 \to M_{k(i)}$$
with the embeddings $I(i, l): M_i^l \hookrightarrow M_{k(i)}$.
There is a unique  $a>0  $ so that
 the translation 
 $v'(m) = v(a\ m)$ on ${\bf C} \subset {\bf C}P^1$
 satisfies the condition
 $\int_{D(1)} (v')^*(\omega_{k(i)}) = \frac{1}{2}<\omega, \alpha>$
 as in $2.A$.  This gives the embedding:
 $$ {\frak B}(M_i^l,  \omega_i^l, J_i^l) \subset  {\frak B}(M_{k(i)},  \omega_{k(i)}, J_{k(i)}) .$$
Notice that if $l$ is sufficiently large, then this is very near just the 
induced  maps by the composition with $I(i,l)$,
namely $a $ above is near $1$.

Let us  consider a  holomorphic curve $u \in {\frak M}([(M_i^l,  \omega_i^l, J_i^l)])$,
and  regard:
 $$u \in {\frak B}([(M_i,  \omega_i, J_i)]).$$
As in $2.C$, 
let  $\hat{U}(u)$ be the  Hilbert completion of the open subset $U(u)$  in 
${\frak B} ( [(M_i,  \omega_i, J_i)] )$, and 
consider 
$\bar{\partial}_J : 
\hat{U}(u) \mapsto 
\hat{\frak E}|\hat{U}(u)$ with  the differential $D\bar{\partial}_J: 
T_u \hat{U}(u) \mapsto \hat{ \frak E}_u$.

\begin{defn}   ${\frak N}^m([(M_i, \omega_i, J_i)])$
is strongly regular, if 
for any  element $[[(M_i^l, \omega_i^l, J_i^l)]]$,
 there is $l_0 $
and  $C >0$  so that  for  all $l \geq l_0$
and for any 
 $u \in {\frak M}([(M_i^l,  \omega_i^l, J_i^l)])$:
$$D\bar{\partial}_J: 
T_u \hat{U}(u) \mapsto \hat{ \frak E}_u$$ is surjective
over $[(M_i,  \omega_i, J_i)]$. Moreover 
the linear inverse  $( D\bar{\partial}_u)^{-1}$
satisfies  the    uniform estimate:
$$ C \geq    |(D \bar{\partial}_u)^{-1}|.$$
 \end{defn}
  Notice that
    $[(M_i, \omega_i, J_i)]$ must be  strongly regular, while
$[(M_i^l,\omega_i^l, J_i^l)]$ might not be strongly regular.

\begin{prop} Let $[(M_i, \omega_i, J_i)]$ be a minimal
  almost Kaehler sequence.
  
Suppose:

\vspace{2mm}

(1)   the  moduli space
${\frak M}([(M_i, \omega_i, J_i)])$ is $S^1$ freely $1$  dimensional, and 

\vspace{2mm} 

(2) ${\frak N}^m([(M_i, \omega_i, J_i)])$ is  strongly regular.
\vspace{2mm} \\
Then
for any  element
$ [[(M_i^l, \omega_i^l, J_i^l)]] \in {\frak N}^m([(M_i, \omega_i, J_i)])$, 
there is  $l_0>>0$ 
so that  the $S^1$ equivariant homeomorphisms:
$$ {\frak M}([(M_i^l, \omega_i^l, J_i^l) ])  \cong
  {\frak M}([(M_i, \omega_i, J_i)])$$
 hold  for all $l \geq l_0$.
   \end{prop}
$Proof:$ 
Combining with lemma $3.1$, it is enough to construct 
the $S^1$ equivariant embeddings
$ {\frak M}([(M_i^l, \omega_i^l, J_i^l) ])  \hookrightarrow
  {\frak M}([(M_i, \omega_i, J_i)])$
   for all $l \geq l_0$.

Let us take  $u \in {\frak M}(M_i^l, \omega_i^l, J_i^l)$.
By lemma $2.2$, there are constants $c_{\alpha}$ independent of $u$
so that the pointwise estimates 
$|\nabla^{\alpha} u| \leq c_{\alpha}$ hold.

Let us  regard $u \in   {\frak B}(M_{k(i)}, \omega_{k(i)}, J_{k(i)})$, and
let $\hat{U}(u)$ be the Hilbert manifolds over $\cup_i M_i$ as above.
 It follows from the  uniform estimates that 
  there is  small  $\delta >0$ independent of $u$ so that
    $\delta$ ball $B_{\delta}(u) \subset T_u \hat{U}(u)$
can be regarded as an open subset in $\hat{U}(u)$.

It follows from the  strong regularity condition that for any  element
$ [[(M_i^l, \omega_i^l, J_i^l)]] \in {\frak N}^m([(M_i, \omega_i, J_i)])$, 
there is  $l_0>>0$
so that  the  $S^1$ freely equivariant embeddings:
$$\Phi : \hat{\frak M}([(M_i^l, \omega_i^l, J_i^l)] ) \hookrightarrow
\hat{\frak M}( [(M_i, \omega_i, J_i)])$$
are given into the smooth manifolds
 for all $l \geq l_0$ (see above $2.C.2$).

 Let us   denote the holomorphic projection 
 $\pi_k : \hat{U}(u) \mapsto {\frak B}(M_k, \omega_k, J_k)$.
Then $v \equiv \pi_k( \Phi(u))$ must satisfy  $\bar{\partial}_J(v)=0$
for all large $k$. It  implies
  the equality $v = \Phi(u)$ for some  large $k$,
  since ${\frak M}([(M_i, \omega_i, J_i)]) $ are $1$ dimensional smooth manifolds,
and    
    the injectivity radii on 
  ${\frak M}([(M_i, \omega_i, J_i)]) \subset {\frak B}([(M_i, \omega_i, J_i)])$
  are uniformly bounded from below
   by  uniform surjectivity of    $D \bar{\partial}_J$.
   This completes the proof. 
  \vspace{3mm} \\
  {\bf 3.C.2 Strong regularity on infinitesimal  neighbourhoods:}
We have introduced a notion of  quasi transitivity in $1.B.4$.
Such property will be   satisfied if the automorphism group is sufficiently large.
For example  ${\bf C}P^{\infty}$ is the case.
In general,  regularity on an almost Kaehler sequence
does not imply the strong  one on its infinitesimal neighbourhoods.
In order to guarantee such property,
we will require high symmetry and integrability.

\vspace{3mm}

\begin{prop} 
Let    $[(M_i, \omega_i, J_i)]$  be a  
minimal and quasi-transitive
 almost   Kaehler sequence. 
Suppose   it is strongly regular.

Then 
$ {\frak N}^m([(M_i, \omega_i, J_i)])$ 
is  also strongly regular.
\end{prop}
{\em Proof:}
Let us use the notations in $3.A.1$.
 Let $U(k,l) \subset M \equiv \cup_i M_i$ be the open neighbourhoods of 
$I(k,l)(M_k^l)$, and extend $\omega_k^l$ and $J_k^l$ as $\omega(k,l)$
and $ J(k,l)$ over $U(k,l)$ respectively.

{\bf Step 1:}
 We  verify   that for any $i >0$,
there is a large $l_0 =l(i)$ so that
 $(U(i,l), \omega(i,l), J(i,l))$
is strongly  regular at $u$
for   any holomorphic curve
$u \in {\frak M}(M_i^l, \omega_i^l, J_i^l)$  for all $l \geq l_0$.

We claim that 
 there exists
 $v \in {\frak M}(M_{k(i)}, \omega_{k(i)}, J_{k(i)})$ which is sufficiently near $u$
in ${\frak B}([(M_i, \omega_i, J_i)])$.
Suppose contrary. Then there exist $\epsilon >0$
and $u_l \in {\frak M}(M_i^l, \omega_i^l, J_i^l)$ so that 
$\epsilon$ neighbourhoods $B_{\epsilon}(u_l) \subset {\frak B}([(M_i, \omega_i, J_i)])$
 contain no elements in  ${\frak M}(M_{k(i)}, \omega_{k(i)}, J_{k(i)})$.
 However by lemma $2.2$, a subsequence must converge to 
 a solution  in $ {\frak M}(M_{k(i)}, \omega_{k(i)}, J_{k(i)})$.
  This verifies the claim.
  
Since 
 $(U(k,l), \omega, J)$ and $ (U(k,l), \omega(k,l), J(k,l))$ are
uniformly near and since   strong  regularity is an open condition,
this implies that the latter is strongly  regular at $u$.

{\bf Step 2:}
Let us fix   sufficiently large  $N >0$ and
 $N$ points $m_0, \dots, m_{N-1}  \in S^2$ with
 $m_0 =0,   m_1= \infty$.
$d(u(m_i),u(m_j))$ are uniformly bounded 
for any $u \in {\frak M}([(M_i^l, \omega_i^l, J_i^l)])$
by   lemma $2.2$.

Let us use the condition of quasi transitivity, and
 choose an automorphism $A$ on $M= \cup_i M_i$ 
with respect to $\{ p_i = u(m_i) \}_{i=0}^{N-1}$
as in $1.B.4$, so that $A(p_i) \in M_{k_0}$ for some $k_0$.

Strong  regularity at $u$ is equivalent to 
that on $A \circ u$, which follows from uniformity of 
complete local charts.
Thus by replacing $u$ by $ A \circ u$, 
one may assume $u(m_i) \in M_{k_0} \subset M= \cup_i M_i$.

 There is  small $\epsilon >0$
independent of $u$ so that the bounds
 $d(M_{k_0}, u(m)) < \epsilon$ hold at any $m \in  S^2$
 again by lemma $2.2$.
Let $U_{\epsilon}(M_{k_0}) \subset M = \cup_i M_i$
be the $\epsilon$  neighbourhood, and 
   $\pi_{k_0}:U_{\epsilon}(M_{k_0}) \mapsto M_{k_0}$ be 
 the  holomorphic projection.
Then by composition, the projection 
 $u' = \pi_{k_0}(u) :S^2 \mapsto M_{k_0}$ is a small deformation 
 of $u$, where 
the uniform estimates:
   $$||u - u'||C^{\alpha} < \epsilon_{\alpha}$$ 
   hold by  lemma $2.2$. 
Here $\epsilon_{\alpha}$ are independent of $i, l$ of ${\frak M}([(M_i^l, \omega_i^l, J_i^l)])$.

Let $\bar{\partial}_J$ be the C-R operator with respect to 
$[(M_i, \omega_i, J_i)]$. Then 
there is a small constant $\delta >0$ with 
$||\bar{\partial}_J(u')| | < \delta$.

{\bf Step 3:}
We show that 
 $l_0$ in step $1$ can be chosen independently of choice of $i$.
This is enough to verify the proposition.

Let us choose 
$u_l \in  {\frak M}([(M_i^l, \omega_i^l, J_i^l)])$ and put
$u'_l =\pi_{k_0}(u_l)$.
 By step $2$, these family  admit uniformly bounded derivatives,
 and $||\bar{\partial}_J(u'_l)||$ converges to zero uniformly as $l \to \infty$.
 So a subsequence 
  of $u'_l$ converges to some $u$ as $l \to \infty$
and $\bar{\partial}_J(u)=0$ holds.
 This implies that
  $u'_l$ are
strongly  regular on $[(M_i, \omega_i, J_i)]$
 for all large $l \geq l_0$ where $l_0$ are independent of $i$.
Then 
$u_l$  are also the same by the above uniform estimates.

  This  completes the proof.
  \vspace{3mm} \\
{\em Example 3.1:}
${\frak N}^m([({\bf C}P^i, \omega_i, J_i)])$ is strongly regular.

\vspace{3mm}

Combining with proposition $3.1$ and $3.2$, we obtain the following:
\begin{cor} Let $[(M_i, \omega_i, J_i)]$ be a minimal,
quasi-transitive and strongly regular 
  almost Kaehler sequence.
  
  If    the  moduli space
${\frak M}([(M_i, \omega_i, J_i)])$ is $S^1$ freely $1$  dimensional, 
then
for any  element
$ [[(M_i^l, \omega_i^l, J_i^l)]] \in {\frak N}^m([(M_i, \omega_i, J_i)])$, 
there is  $l_0>>0$ 
so that  the $S^1$ equivariant homeomorphisms:
$$ {\frak M}([(M_i^l, \omega_i^l, J_i^l) ])  \cong
  {\frak M}([(M_i, \omega_i, J_i)])$$
 are given  for all $l \geq l_0$.
   \end{cor}

\vspace{3mm}

Now we collect the previous results which we have induced so far:
\begin{thm} Let  $[(M_i, \omega_i, J_i)]$
 be a regular and minimal almost Kaehler sequence.
 
    (1)  Suppose  it is isotropic and symmetric. 
     If each connected component of the moduli space is bounded, 
     then the equality holds:
     $$ {\frak M}([(M_i, \omega_i, J_i)])= {\frak M}(M_0, \omega_0, J_0).$$
   
   In particular if     
    the moduli space is $1$ dimensional, then the above equality holds,
     which are both compact.

   (2)  If it is  symmetric 
    Kaehler, then it is strongly regular.
    
    If moreover it is quasi-transitive, then 
    $ {\frak N}^m([(M_i, \omega_i, J_i)])$ 
is   strongly regular.

(3) Under all the conditions in (1) and (2), 
it follows that for any  element
$ [[(M_i^l, \omega_i^l, J_i^l)]] \in {\frak N}^m([(M_i, \omega_i, J_i)])$, 
there is  $l_0>>0$ 
so that  the $S^1$ equivariant homeomorphisms:
$$ {\frak M}([(M_i^l, \omega_i^l, J_i^l) ])  \cong
  {\frak M}([(M_i, \omega_i, J_i)])$$
 are given  for all $l \geq l_0$.
\end{thm}

\section{Application to Hamiltonian dynamics}
In this section we apply theory of moduli spaces we developed 
so far to  study of Hamiltonian dynamics 
defined by smooth and bounded functions over almost Kaehler sequences.
\vspace{3mm} \\
{\bf 4.A Capacity invariant:}
In Hamiltonian dynamics,
capacity invariant contains deep information on the periodic solutions
and has been playing  one of the central roles 
 over finite dimensional symplectic manifolds.

Let $(M, \omega)$ be a closed symplectic manifold of finite dimension.
  We say  that a smooth function  $f: M \mapsto [0, \infty)$ is Hamiltonian.
The {\em Hamiltonian vector field} $X_f$ is uniquely defined by the relation:
$$df( \quad) = \omega(\quad, X_f).$$
A {\em periodic solution} $x:[0,T ]\to M$ with $x(0)=x(T)$ satisfies  the equation:
 $$\dot{x} = X_f(x)$$
and $T\geq 0$ is called {\em period}. 
\vspace{3mm} \\
{\em Remark 4.1:}
Let $a >0$ be a positive number and $\tilde{f}(x)=af(x)$ be the function multiplied by $a$.
If $x:[0, T] \to M$ is a periodic solution to $f$
with period $T$, 
then $\tilde{x}:[0,a^{-1}T] \to M$ given by $\tilde{x}(t) =x(at)$
is also a periodic solution to $\tilde{f}$,
since the equality  $X_{\tilde{f}} =aX_f$ holds.
So its dynamical properties are scaling invariant under multiplication by positive numbers.
Moreover the   `height' of functions is in inverse proportion to  periods.

\vspace{3mm}

Let us introduce an invariant over symplectic manifolds 
with respect to periodic solutions ([HZ]).
A Hamiltonian function  is {\em pre admissible} if
there are open sets $U , V \subset M$ with $f|U \equiv c =  \sup f$
and $f|V \equiv 0$. Moreover
it is {\em  admissible}  if in addition,  any periodic solution
  is either constant or period $T= T(x) >1$, where
$x: [0, T]  \to M$ with $x(0)=x(T)$,  and $X_f$
is   the Hamiltonian vector field.

Let us denote the set of admissible functions by $H_a(M, \omega)$. Then we
 define the {\em capacity} of $(M, \omega)$ by the following:
$$c(M, \omega) =     \sup \{ m(f)   \equiv \sup f -  \inf  f \geq 0
: f \in H_a(M, \omega)\}.$$

Let $N$ be a symplectic manifold with boundary. $f: N \mapsto [0, \infty)$ 
is pre admissible,  if there is an open set $U$  with $f|U \equiv c = \sup f$ and
 it vanishes on some neighbourhood of boundary.
It is  admissible if in addition, 
any periodic solution $\dot{x} = X_f(x)$ is either constant or period $>1$.
We define the capacity by the same way over symplectic manifolds with boundary.

This numerical invariant satisfies some axioms of capacity. 
 Notice that  admissible functions always exist,
since the Darboux's  chart exists at any point, 
 and $c(D, \omega) = \pi$  where
$(D, \omega)$ is the standard symplectic  disk.

Let us consider the upper bounds of the invariants.
 Let   $(M, \omega, J)$  be a minimal
almost Kahler manifold of finite dimension,
and let us fix a minimal element $\alpha \in \pi_2(M)$.

Let us state a basic  relation between periodic solutions and holomorphic curves:
 \begin{lem}[HV]  
 Let $f: M \to [0, \infty)$ be a pre-admissible function.
 Then there are non trivial periodic solutions 
whose  periods satisfy the estimates 
$$T < m(  f)^{-1}m$$
where $m$ is the minimal invariant, 
if the moduli space of holomorphic curves with respect to $\alpha$
is non empty, regular, 
$1$ dimensional and $S^1$ freely cobordant to non zero.
\end{lem}
 In particular under the above conditions, the estimates:
 $$c(M, \omega) \leq m$$
hold.
This is used to estimate capacity invariants over almost Kaehler sequences defined below.
\vspace{3mm} \\
{\bf 4.B Asymptotic periodic solutions:}
Let $[(M_i, \omega_i, J_i)]$ be an almost Kaehler sequence.
 Later on we  fix a uniformly bounded covering 
 by $\epsilon >0$ complete almost Kaehler charts.

Let $f: M \mapsto {\bf R}_+$ be a bounded Hamiltonian function
and put $f_i = f|M_i$.
A family of smooth loops 
 $x_i: [0, T_i] \mapsto M_i$ with
$x_i(0)=x_i(T_i)$ 
is a periodic solution
over $[(M_i, \omega_i, J_i)]$, if these satisfy the equations:
$$\dot{x_i} (t)= X_{f_i}(x_i(t))$$
for all $0 \leq t \leq T_i$ and  $i=0,1, \dots$
We will denote such a family by $[x_i]$.

    \vspace{3mm}
    
 Let    $ [[(M_i^l,  \omega_i^l, J_i^l)]] \in {\frak N}([(M_i,  \omega_i, J_i)])$
 be an element of the infinitesimal neighbourhood.
Let us consider a family of bounded Hamiltonians:
 $$f_l : M^l \equiv \cup_i M_i^l \to  {\bf R}$$
  and take a family of 
 periodic solutions 
 $[x_i^l]$ over $[(M_i^l,\omega_i^l,J_i^l)]$ with their periods $T^l_i$.
 Passing through the embeddings $I(i,l): M_i^l \hookrightarrow M_{k(i)}$, 
 one may regard these loops as:
 $$x_i^l : [0, T_i^l] \to M_{k(i)}.$$
 
 Let $\pi_i : U_{\epsilon}(M_i) \to M_i$ be the holomorphic projections.
    
   \begin{prop} Suppose  $  [(M_i,  \omega_i, J_i)]$ is  quasi transitive, and 
   both the uniform bounds hold for all $\alpha  \geq 0$:
$$  \sup_l  ||f_l|| C^{\alpha}(M^l) \leq C_{\alpha}, \quad \sup_{i,l} T^l_i  \leq T < \infty.$$
   Then  there is a bounded Hamiltonian
$g: M \equiv \cup_i M_i \to {\bf R}$ and a family of automorphisms $A_i $ 
over $ (M, \omega, J)$
so that for some $l_i$, 
the family of loops   $z_i \equiv \pi_i(A_i \circ x_i^{l_i}): [0, T_i^{l_i}] \to M_i$ 
satisfy the following asymptotics:
  $$\lim_{i \to \infty} \  \sup_t  \ |\dot{z_i} - X_{g_i}(z_i)| (t) \ = \  0.$$
   \end{prop}
       {\em Proof:}     Let us recall the notations   in $3.A.1$.

   {\bf Step 1:}
   Let us choose a sufficiently large $N$ and $N$ points
   $0 \leq  t_i = \frac{T}{N-1}i  \leq T$ for
   $i =0,1, \dots, N-1$.
   By quasi transitivity, there are some $k$ and automorphisms $A_i^l $    over 
   $[(M_i, \omega_i, J_i)]$ so that:
   $$y_i^l(t_i) \equiv A_i^l(x_i^l(t_i)) \in M_k$$
   hold for all $i $.
   By the assumption of the uniform bounds, 
   there is small $\delta>0$ with the bounds in $M$  for all $0 \leq t \leq T_i^l$:
   $$d(y_i^l(t), M_k) \leq \delta.$$

{\bf Step 2:}
    $M^l_i$
   admit the embeddings $I(i, l): M^l_i    \hookrightarrow M$,
    and  there are some open subsets:
       $$I(i, l)(M^l_i) \subset U(i,l)     \subset M$$ 
    where    $U(i,l)$ contain
      $   \epsilon$ neighbourhooods of $I(i,l)(M^l_i)$ for some $   \epsilon >0$.

Let us fix $i$, and consider the restrictions $f_l|M_i^l$ and regard them as:
    $$f_l: I(i, l)(M^l_i) \to {\bf R}.$$ 
 We  extend them as
    $f_i^l: M \to {\bf R}$  as follows.
   Let 
    $\pi_i^l : U(i,l) \to I(i,l)(M^l_i)$ be the holomorphic projections. 
   Let $\mu: [0,\epsilon) \to [0,1]$ be a cut off function 
   with $\mu(m)=0$ for all $m \geq \frac{1}{2}\epsilon$
  and $\mu \equiv 1$ near $0$.
   Then we put:
   $$f_i^l(m)\equiv \mu(d(m, M^l_i)) \ f_l(\pi_i^l(m))$$
   for $m \in U(i,l) \subset M$, and $f_i^l| M \backslash U(i, l) \equiv 0$.
   Let us put
   $h_i^l \equiv f_i^l \circ (A_i^l)^{-1}$.
   These families satisfy the following properties:
   
(1)     $\{h_i^l\}_{i,l}$ are   uniformly bounded as
   $||h_i^l||C^{\alpha}(M) \leq C_{\alpha}$ for all $\alpha$.
   
   (2) For $\bar{h}_i^l \equiv h_i^l|M_i$,
      $x_i^l :[0, T_i^l] \to M_{k(i)}$ satisfy the asymptotics:
   $$\lim_{l \to \infty} \sup_t |\dot{y}_i^l (t) - X_{\bar{h}_i^l}(y_i^l(t))|=0.$$
   
      (3) For the projections $\pi_i: U_{\epsilon}(M_i) \to M_i$,
      $$\lim_{i \to \infty} ||y_i^l - \pi_i(y_i^l)||C^1 =0.$$

   {\bf Step 3:}
 A  subsequence of the family $\{h_i^{l_i}\}_i$ converges weakly to   
a bounded Hamiltonian $g$ 
so that
for $z_i \equiv \pi_i(y_i^{l_i}) $, the asymptotics hold:
$$\lim_{l \to \infty} \sup_t |\dot{z}_i (t) - X_{g_i}(z_i(t))|=0$$
by  step $2$ and lemma $1.3$.
      This completes the proof.
\vspace{3mm}

In the context of this paper,  asymptotic analysis
play one of the central roles in studying increasing sequences of manifolds.
Proposition $4.1$   presents a motivation quite naturally
to introduce some  families of loops  
which approach  periodic solutions asymptotically.
Conversely  asymptotic periodic solutions defined below
can be regarded as though  a kind of `periodic loops'  over the infinitesimal neighbourhoods.

Let $[(M_i, \omega_i, J_i)]$ be an almost Kaehler sequence.
Let $ [x_i ]$ be a family of smooth loops:
 $$x_i: [0, T_i] \mapsto M_i$$ with
$x_i(0)=x_i(T_i)$ and $i=0,1, \dots$

We say that $[x_i]$ is an {\em asymptotic loop}, 
if     $ T= \limsup_i T_i < \infty$ is finite.
We call $T$  as the {\em periods} of the family.

An asymptotic loop $[ x_i]$ is {\em small}, if
 $0 \leq   \limsup_i  T_i \leq 1$ holds.

 We say that
    an asymptotic  loop $[x_i]$  is {\em non trivial}, if 
    uniform estimates:
$$ \text{ length } x_i \equiv \int_0^{T_i} |\dot{x}_i(t) |dt \geq \delta >0$$  hold
    for all  large $i \geq i_0 = i_0([x_i])$ and 
   some positive $\delta  = \delta([x_i])>0$.

\vspace{3mm}

 \begin{defn}
Let $f: M \mapsto {\bf R}_+$ be a bounded Hamiltonian function
and put $f_i = f|M_i$.
An  asymptotic periodic solution  is  an asymptotic loop $[x_i ]$
 satisfying:
  $$\sup_t |\dot{x_i} - X_{f_i}(x_i)| (t) \to 0, \quad i \to \infty.$$
   \end{defn} 
  \quad 
   \vspace{3mm}  \\
  { \em Remark 4.2:}
  (1)
   Trivial asymptotic periodic solutions $[x_i]$ with length $x_i \to 0$,   always exist 
with  any periods, if $f$ attain  maximum or minimum values.

(2) Any asymptotic periodic solutions $[x_i]$ must satisfy uniform bounds:
$$ \limsup_i  \text{ length } x_i <  \infty.$$

(3)
Suppose a family of bounded familtonians $\{f^l\}$
over an almost Kaehler sequence
converges as $||f^l -f||C^{\alpha}(M) \to 0$
for all $\alpha$.
Let $\{x_i^l\}_i$ be families of periodic solutions with respect to $f^l$.
Then $y_i \equiv x_i^i$ is an asymptotic periodic solution
with respect to $f$.
\vspace{3mm} \\
     {  \bf 4.B.2 Comparison with periodic solutions:}
  Let us compare asymptotic periodic solutions with exact solutions.

  Let us say that an asymptotic periodic solution $[x_i]$ is {\em exact}, 
  if $x_i$ are periodic solutions:
   $$ \dot{x}_i - X_{f_i}(x_i)=0$$
   to the restrictions $f|M_i$ 
  for infinitely many $i$.

  There are almost Kaehler sequences and bounded Hamiltonians over them
  so that there are no exact priodic solutions but do exist asymptotic periodic solutions.
    \vspace{3mm} \\
   { \em Example 4.1:}
   Let $N_i=(T^4,  \omega_i)$ be a family of 
   symplectic tori, where:
    $$ \omega= dx_0 \wedge dy_0 +dx_1 \wedge dy_1 + a_i dx_0 \wedge dx_1$$ with $a_i  \to 0$.
   If we choose $a_i =  \frac{n_i}{m_i}$ with $n_i, m_i  \to  \infty$,
   then any non trivial periodic solutions with respect to, say $f(x,y)= \cos(2 \pi x_0)$
   must have divergent periods as $i  \to  \infty$.
   
   In particular the almost Kaehler sequences:
    $$M_i=(N_0  \times  \dots \times N_i,  \omega_0 + \dots + \omega_i)$$
    certainly admits asymptotic periodic solutions with the period $1$,
    with respect to the bounded Hamiltonian: 
    $$f(x_0,y_0, x_1,y_1, \dots) =\cos(2 \pi x_0)$$
   but they do not admit any periodic solutions.
   Here we embed $M_i 
    \hookrightarrow M_{i+1}$ by $(z_0,   \dots, z_i) 
       \to (z_0,   \dots, z_i,0)$.
       \vspace{3mm}  \\
      { \bf 4.C  Capacity functions over almost Kaehler  sequences:}
      In this section we introduce $3$ variants  of capacity invariants:
    $$   \text{cap},  \ \ \text{As-cap} , \ \  \text{C} $$
  over almost Kaehler sequences. cap is the most standard one and 
  extends  the finite dimensional case directly.
  As-cap and C both  arose by looking at 
  `infinitesimal stabilizations' of cap.
      We verify that in some case they coincide each other.
         \vspace{3mm}  \\
       {   \bf 4.C.1 Capacity on almost Kaehler sequences:} 
      Let $f: M  \to [0, \infty)$ be a bounded Hamiltonian
       over
an almost Kaehler sequence $[(M_i, \omega_i, J_i)]$.
We say that $f$ is  {\em pre admissible}, if
there exist open sets $U, V \subset M$ with $f|U \equiv \sup f$
and $f|V \equiv 0$.

 Let us say that $f$ is {\em admissible}, if it is pre admissible, and
 for any non trivial asymptotic loops $[x_i]$, 
 if it consists of 
  periodic solutions  $x_i:[0,T_i]  \to M_i$ with respect to the restrictions
  $f_i \equiv f|M_i: M_i \to {\bf R}$ for all $i \geq i_0 \geq 0$, 
  then thier periods 
satisfy the bounds:
   $$ T(x_i)= T_i> 1.$$

Let $[(X_i, \omega_i, J_i)]$ be an almost Kaehler sequence  with boundary 
$\partial X = \cup_i \partial X_i$ with
 $\partial X_i \subset  \partial X_{i+1}$.
A  bounded Hamiltonian $f: \cup_i X_i \mapsto [0, \infty)$
is {\em pre admissible}, if there is an open set $U \subset  \cup_i X_i$ with $f|U \equiv $ sup$f$,
and  $f|V \equiv 0$ holds on $\delta >0$ neighbourhood $V$ of the boundary
$\partial X  \subset V$ for some positive $\delta =\delta(f) >0$.
We say that $f$  is {\em admissible}, if  it is pre admissible and
 satisfies the above property.

Let us put the set of admissible functions by
$H_a([(M_i,\omega_i, J_i) ]) $. 
Recall $m(f)= \sup f - \inf f \geq 0$.

 \begin{defn} The capacity function cap is  
    defined by:
    $$\text{cap}([(M_i, \omega_i, J_i)]) = \text{ Sup}     \{ m(f) :    f \in H_a([(M_i, \omega_i, J_i)]) \}.$$
  \end{defn}
$cap$ takes values in ${\bf R}_+ \cup \{ + \infty \}$.

\begin{prop}
 Let   $[(M_i, \omega_i, J_i)]$  be a minimal
almost Kahler sequence with a fixed minimal element $\alpha \in \pi_2(M)$
with $m =<\omega, \alpha>$.

Then 
the estimates:
$$\text{cap}([(M_i, \omega_i, J_i)]) \leq \ m$$
hold, if the 
moduli space ${\frak M}([(M_i, \omega_i, J_i)])$
 is non empty, regular,
has $1$ dimensional and $S^1$ freely cobordant to non zero.
\end{prop}
{\em Proof:}
Let $f : M \to {\bf R}$ be pre admissible, and $f_i : M_i \to {\bf R}$ be its restrictions.
It is enough to check that there is a family of periodic solutions
$x_i : S^1 \to M_i$ with the properties:

(1) their periods satisfy uniform lower bounds $T(x_i) > \delta >0$ from below, 

(2) they consist of non trivial asymptotic loops with length $x_i > \delta' >0$.
\vspace{2mm} \\
By lemma $4.1$, periodic solutions $x_i: S^1 \to M_i$ exist certainly.
In fact the construction by [HV] with the proof of lemma $2.2$
verifies the above two properties for these periodic solutions.
We outline the proof below.

Let  $B_{\mu}(p_0), B_{\mu}(p_{\infty}) \subset M = \cup_i M_i$ be $\mu >0$ balls
so that $f|B_{\mu}(p_0) \equiv 0$ and $f|B_{\mu}(p_{\infty}) \equiv \sup f$.

{\em Step 1:}
For $c >0$, there exist $\epsilon >0$ so that for any  $s_0$ and solution 
$(\lambda, u) $  to the equation $(f_i)_{\lambda}(u)=0$ as in 
[HV] $(2.12)$ on page $  599$ with:
\begin{align*}
& |\nabla u(s,t)| \leq c  \quad \text{ for } s \in ( - \infty , s_0] \times S^1, \\
& u((- \infty, s_0] \times S^1) \cap B_{\mu}(p_0)^c \ne \phi
\end{align*}
then the lower estimates hold:
$$\int_{-\infty}^{s_0} \int_0^1 u^*(\omega) \geq \epsilon.$$

The proof goes by the same way as [HV] in lemma $3.1$,
where we introduce minor modifications of the arguments 
 in order to induce the uniform bounds as above.

Let us consider smooth maps $\gamma : S^1 \to B_{\mu}(p_0)$
with the $L^2$ bounds  length $\int_{S^1} |\dot{\gamma} |^2 \leq c $
for a constant $c >0$. 
 Let us regard $\gamma : S^1 \to  B_{\mu}(p_0) \subset H$ as maps into the Hilbert space,
 and consider  $p_0 - \gamma : S^1 \to H$.
We claim that for any $\delta >0$, 
 there is a constant $\tau=\tau(\delta,c)$
so that if $ ||p_0-  \gamma ||L^2 \leq  \tau$, then
$\text{ diam } \gamma \leq  \delta $ hold.
The argument in [HV] $p 608$ almost works except that we are concerning maps into 
infinite dimensional spaces, and so we cannot use the maps themselves directly 
to apply the Ascoli-Arzela Theorem.
We just modify the argument by replacing $p_0 - \gamma$ by 
  the  functions $|p_0 - \gamma| : S^1 \to [0, \infty)$ by taking the pointwise norms.

Let $\gamma_k : S^1 \to  B_{\mu}(p_0)$ ba a family of smooth maps 
with the uniform bounds $\int_{S^1} |\dot{\gamma}_k |^2 \leq c $.
By the compact embedding $L^2_1(S^1)  \hookrightarrow C^0(S^1)$
between the real valued  functional spaces, 
if $|p_0 - \gamma_k| :S^1 \to [0, \infty)$ converges to zero in $L^2$, 
then they must also converge to zero in $C^0$.
 In partucluar the diameters of $\{\gamma_k \}_k$ 
also converge to zero. This is enough to conclude the claim.

Let us introduce the inifinite dimensional Poincar\'e lemma,
which states if a closed form on a ball in the Hilbert space satisfy $C^0$ bounds, then
one can find  its primitive form which lies in the dual space of $L^2$ vectors.
Let:
$$\omega = \Sigma_{i> j \geq 0} \ a_{ij} dx_i \wedge dx_j$$
be a closed form on $B_{\mu}(p_0)$ with the uniform bounds
 $||a_{i,j}||C^0(B_{\mu}(p_0)) \leq C$.
Then the primitive $\theta$ with $d \theta = \omega$ is given by the one form:
$$\theta = \Sigma_{i>j \geq 0}  \ (  \  \int_0^1 a_{ij} ( t  \bar{x})  \ t \  dt \  )  \ x_i \  dx_j$$
which satisfies  $d \theta =\omega$ and is clearly  $L^2$ one form.

As a result  there is a constant $c$ with the estimates 
$|\int \gamma^*(\theta)| \leq c \int |\dot{\gamma}|^2$
for any loops $\gamma : S^1 \to B_{\mu}(p_0)$ by the Cauchy-Schwartz estimate.

{\em Step 2:}
There are uniform constants $c_{\alpha}$ so that 
any solutions 
$(\lambda, u) $  to the equations $(f_i)_{\lambda}(u)=0$ satisfy 
uniform bounds $||u||C^{\alpha}(S^2) \leq c_{\alpha}$.
This is verified by following the argument in [HV] except taking the limit
to obtain holomorphic spheres, where instead we use a similar method
as the proof of lemma $2.2$ particularly step $4$.

The rest of the arguments are also parallel again by use of a similar method as 
step $4$ as above.
This completes the proof.
 \vspace{3mm}  \\
{   \bf 4.C.2 Asymptotic  capacity:}   
      Below we introduce two stable versions of the capacity function.
             Let $f: M  \to [0, \infty)$ be a bounded Hamiltonian.
 Let us say that $f$ is { \em as-admissible}, if it is pre admissible, and
any  non trivial  asymptotic periodic solution has
 its period:    $$\liminf_i T(x_i) > 1.$$

Similarly we define as-admissibility for bounded Hamiltonians on  almost Kaehler sequences  with boundary.

Let us put the set of as-admissible functions by
$H_{\text{as}}([(M_i, \omega_i, J_i)]) $, and 
 we define the {\em asymptotic  capacity}:
$$\text{As-cap}([(M_i, \omega_i, J_i)]) = \text{ Sup}     \{ m(f) :   f \in H_{\text {as}}([(M_i, \omega_i, J_i)]) \}.$$
As-cap also takes values in ${\bf R}_+ \cup \{ + \infty \}$.
Notice that a-priori estimate: 
 $$\text{As-cap}([(M_i, \omega_i, J_i)])    \leq \text{cap}( [(M_i, \omega_i, J_i)])$$
  holds.
    \vspace{3mm} \\
  {   \bf 4.C.3 Capacity over infinitesimal neighbourhoods:}   
        Let us introduce the geometric version of 
       stabilized capacity function
      which is formulated by use of the infinitesimal neighbourhoods.
      
     Let $[(M_i,  \omega_i, J_i)]$ be an almost Kaehler sequence.
     The {\em capacity invariant over   infinitesimal neighbourhoods}
   are defined by:
\begin{align*}
 C^m([(M_i, \omega, J_i)]) =    
   \sup \{  \   \lim \sup_l  &  \text{ cap}([(M_i^l, \omega_i^l, J_i^l)]) : \\
&             [[(M_i^l, \omega_i^l, J_i^l)]] \in  {\frak N}^m([(M_i, \omega_i, J_i)]) \}
            \end{align*}

Notice that a priori estimates hold:
 $$\text{As-cap}  \ \leq    \text{cap} \  \leq  \ C^m$$

 \begin{thm} Let $[(M_i, \omega, J_i)]$ be a minimal, 
  isotropically symmetric and quasi transitive  Kaehler sequence,
 with a  fixed minimal element 
 $ \alpha  \in \pi_2(M)$.
 
  If the 
moduli space of holomorphic curves is non empty, regular,
 $1$ dimensional
and $S^1$ freely cobordant to non zero,  then  the estimate:
$$C^m([(M_i,  \omega_i, J_i)]) \leq m$$
holds, where $m =<\omega, \alpha>$.
\end{thm}
{ \em Proof:} 
Strong regularity holds since it is a  regular and symmetric Kaehler sequence.
Then   ${\frak N}^m([(M_i,  \omega_i, J_i)])$ is 
strongly regular since it is minimal and quasi-transitive.
For any  $ [[(M_i^l, \omega_i^l, J_i^l)]] \in  {\frak N}^m([(M_i, \omega_i, J_i)]) $,
there is large $l_0 $ so that for all $l \geq l_0$, 
the moduli spaces over $[(M_i^l, \omega_i^l, J_i^l)]$
are all homeomorphic to the one over $[(M_i, \omega_i, J_i)]$
by  proposition $3.1$.

Then the estimates
$\text{cap} ([(M_i^l, \omega_i^l, J_i^l)]) \leq m$ hold
by combining with  these with lemma $4.1$.
This completes the proof.

\vspace{3mm}

So far we have introduced three versions of the capacity invariants.
We would like to propose:
\vspace{3mm} \\
  {\em Conjecture 4.1:}
   Let $[(M_i, \omega_i, J_i)]$ be a Kaehler sequence.
   Under the conditions of theorem $4.1$, 
        the equalities hold:
  $$   \text{As-cap}([(M_i, \omega_i, J_i)])   = 
  \text{cap}([(M_i, \omega_i, J_i)])  = C^m( [(M_i, \omega_i, J_i)]).$$ 
We verify the equalities   for both the infinite unit disk and ${\bf CP}^{\infty}$ below.
\vspace{3mm}  \\
      {      \bf 4.D Examples:}
Here  we estimate the values of capacities over
some concrete cases. These include infinite dimensional disks, 
infinite projective space and infinite tori.
 \vspace{3mm}  \\
{   \bf 4.D.1 Infinite disks:}   
Let
 $D^{2i} \subset {\bf R}^{2i}$  be the standard unit disks, and $[D^{2i}]$ be
 the standard  Kaehler sequence with
$  \cup_i D^{2i} \subset {\bf R}^{2 \infty}$.

Let us fix $N \geq 1$ and
 put the infinite products  of $D^{2N}$:
 \begin{align*}
&  \tilde{D} = D^{2N} \times D ^{2N} \times \dots, \\
 & D(i) = D ^{2N} \times \dots \times D^{2N}  \subset \tilde{D} \qquad  (i \text{ times})
 \end{align*}
 where we embed $D(i) \subset D(i+1)$ by $(m_1, \dots, m_i) \to (m_1, \dots, m_i,0)$.
 
 Here we verify the following equalities:
 \begin{prop} 
 Let $ \text{Cap} $ be   $ \text{As-cap} $ or $ \text{cap}$. Then
  the equalities hold:  
 $$C^m([{\bf CP}^i])= \text{Cap}([{\bf CP}^i])= \text{Cap}([D^{2i}]) = \pi.$$
\end{prop}
{\em Proof:}
We split the proof into three cases:
$$(1) \  \pi \geq C^m([{\bf CP}^i]) , \quad 
(2)  \   \   \text{As-cap}({\bf CP}^{\infty})  \geq  \text{As-cap}([D^{2i}]) \geq \pi.$$

The proof of (2) is  postponed to the following lemmas below.

Let us consider (1).
$ {\bf CP}^{\infty}$ is minimal with $m =\pi$, 
 since it has  $\pi_2$-rank $1$
 and  the embedding ${\bf CP}^1 \subset {\bf CP}^{\infty}$ is holomorphic.
It is isotropically symmetric by Examples $1.4$ and quasi transitive by lemma $1.7$.

It is strongly regular by example $2.2$, whose moduli space is homeomorphic to $S^1$
which is $S^1$ freely  cobordant to non zero.

So all the conditions in theorem $4.1$ are satisfied, and so (1) follows.
 
\begin{lem} 
$  \text{As-cap}([D(i)])  \geq  \text{As-cap}([D^{2i}]) \geq \pi$.
\end{lem}
{\em Proof:}
The first inequality is clear, since there are compatible embeddings 
$D^{2iN} \hookrightarrow D(i)$.  So we will consider the second one.

We follow a parallel argument  to [HZ] $p72$.
Here one needs to be careful about treating asymptotic solutions.

In [HZ],  $c(D^{2i}) = \pi$ is verified for any $i$. In fact the following is shown;
let $f: [0,1] \mapsto [0, \pi]$ be any function satisfying:

\vspace{2mm}

(1)
$0 \leq -f' < \pi - \delta$, $\delta >0$,

(2) $f(t) = \pi  - \delta'$,  $t \sim 0$, $\delta' >0$
and $f(t) \equiv 0$, $t \sim 1$.
\vspace{2mm} \\
Then  $F: D^{2i} \mapsto [0, \pi]$ 
 gives an admissible function on $D^{2i}$ by $F(x) = f(|x|^2)$.
Let us put $F: \cup_i D^{2i} \mapsto [0, \pi]$ by the same way, which is 
a  bounded Hamiltonian.
We claim that $F$ is as-admissible in our sense. 
Let us take any   non trivial asymptotic periodic solution $[x_i]$.
By the definition, this family satisfies:
$$|J \dot{x}_i + \nabla F(x_i)|C^0
= |J\dot{x}_i + 2f'(|x_i|^2)x_i|C^0 \to 0.$$
Let us put $G(x) = \frac{1}{2}|x|^2$. Then one has the  uniform estimate:
$$\frac{d}{dt}G(x_i(t)) = <\nabla G , \dot{x}_i> \sim 
<\nabla G , J\nabla F>(x_i(t)) =0.$$
Thus  there is a constant 
$0 < a \leq 2(\pi - \delta)$ with 
$\sup_t |2f'(|x_i(t)|^2) +a| \to 0$.

Let us put $a_i = -2f'(|x_i(0)|^2)$ with  $a_i \to a$.
For $y_i \equiv exp(-a_iJt)x_i(0)$, the family $[y_i]$ consists of
a periodic solution.
Then   $\sup_t |y_i(t) - x_i(t)| C^0  \to 0$ hold since 
  the Hamiltonian vector field
 is uniformly Lipschitz  (of completely bounded geometry).
 It also implies $\sup_t  |y_i(t) - x_i(t)| C^1 \to 0$, since
 they are (asymptotic) solutions.
So $\lim_i T(x_i) >1$ follows 
since $T(y_i) > 1 + \delta'$ hold for some  $\delta' >0$.
This completes the proof.
\vspace{3mm}

\begin{lem}
   $\text{As-cap}({\bf CP}^{\infty})  \geq  \text{As-cap}([D^{2i}]) $.
\end{lem}
$Proof:$
 There are standard symplectic embeddings
 $D^{2i}  \hookrightarrow {\bf CP}^i$ from the unit disks.
  Let $J_1$ be the induced almost complex structure, which 
  is different from the standard one $J_0$  on $D^{2i}$.
 Thus the  induced Riemannian metrics
 on  $D^{2i} $ are not standard.  
 Recall that as-admissibility of bounded Hamiltonians involves
 asymptotic behaviour of loops by the Riemannian norms.
 However this does not cause any problem for us, since  
 we use bounded Hamiltonians over $\cup_i D^{2i}$ which 
 vanish near boundary, and so 
 the  induced Riemannian  metrics are uniformly equivalent to the standard one
   on the support of $f$ in the following sense.
   Let us choose any pre-admissible Hamiltonian 
$f$ on  $  \cup_i D^{2i}$     with respect to $J_1$.
One may assume that for any small $\mu >0$ 
and any non trivial asymptotic periodic solution $[x_i]$,
there is a large $i_0$ so that
Supp $x_i $ are all contained in $\mu$ neighbourhoods of
  Supp $f_i$ with  $f_i = f| D^{2i}$
   for all $i \geq i_0$.
Then there is an equivalence of the induced Riemannian metrics,
$ C^{-1}  |\quad |_0 \leq |\quad |_1 \leq C|\quad |_0$,
where $C$ depends only  on $d( \partial \cup_i D^{2i}, \text{  Supp } f)$.
Thus if  $[x_i]$ is an asymptotic  periodic  solution for $J_0$, then it is
also the case  for $J_1$. The converse is the same.
 Thus  the inequality 
 $ \text{As-cap}({\bf CP}^{\infty})  \geq  \text{As-cap}([D^{2i}]) \geq \pi$ follows.
 
This completes the proof.

\vspace{3mm}

Let us see how the capacity invariant behave under small perturbation of the forms.
 Notice that the argument below works only for the products of disks.
 The situation is completely different for the torus case.
 \begin{sublem}
 Let us fix $i$ and let
 $\omega'$ be a  symplectic form  on $D(i)$ which is sufficiently 
 near the standard form $\omega$.
 Then the values of capacity $c(D(i), \omega') \sim c(D(i), \omega)$
 are sufficiently near from each other.
 \end{sublem}
 $Proof:$ This follows from the Darboux-chart  construction.
 Let us denote  $D_i(\delta) = D^{2N}(\delta) \times \dots \times D^{2N}(\delta)$, where
 $D^{2N}(\delta)$ is the $2N$ dimensional    $\delta$ disk.
It is enough to find  small $\epsilon >0$ 
 with $\epsilon \to 0$ as $\omega' \to \omega$, 
 and symplectic
 embeddings:
   $$I: (D_i(1- \epsilon) , \omega) \hookrightarrow (D(i), \omega'), \quad
   I': (D_i(1- \epsilon) , \omega') \hookrightarrow (D(i), \omega).$$
   The constructions of these are parallel, and we will check the first case only.
 
 Let us find a diffeomorphism $\varphi:  D\equiv D_i(1- \epsilon) 
  \hookrightarrow D(i)$ with $\varphi^*(\omega') = \omega$.
 Let us put $\omega_t = \omega + t(\omega' - \omega)$ for $t \in [0,1]$.
 By the assumption, 
 all $\omega_t$ give symplectic forms. 
 We will find a smooth family of diffeomorphisms $\varphi_t: D \hookrightarrow D(i)$ 
 given as the flow of vector field $X_t$, and 
 satisfying
 $\varphi_t^*(\omega_t) = \omega$. The last equality is
 equivalent to $d(i_{X_t} \omega_t) + \omega - \omega' =0$.
 
 By the Poincar\'e lemma, there is a smooth one form $\lambda$ on $D(i)$ with 
 $\omega - \omega' = d \lambda$ and $\lambda(0) =0$.
 The pointwise norm of $\lambda$ will be sufficiently small.
 Then we choose unique $X_t$ satisfying 
 $\omega_t(X_t , \quad ) = - \lambda$ and $X_t(0) =0$.
 Thus $|X_t|$ is also small pointwisely. In particular the flow of $X_t$ exists for all 
 $0 \leq t \leq 1$ with image in $D(i)$,
 if we start from any point $p \in D= D_i(1- \epsilon)$.
 Thus we have obtained the desired $\varphi_t$.
This completes the proof.
\vspace{3mm} \\
{\bf 4.D.2 Based admissibility:}
Let us equip the standard metric  on 
$ \tilde{D} \equiv \cup_i D(i)  \subset {\bf R}^{\infty}$, and let
   $f: \tilde{D}  \mapsto [0, \infty)$ be a  smooth and
bounded  function.
We say that 
a non trivial asymptotic periodic solution  $[x_i]$ is  {\em based}, if
$ \{ | x_i|C^0(S^1) \}$ is uniformly bounded, and 
 there is  some $i_0 \geq 1$
so that  lim inf$_i$ length $P_{i_0}(x_i) >0$, where:
$$P_{i_0}:  \tilde{D} = D^{2N} \times D^{2N} \times  \dots \mapsto  
D(i_0) = D^{2N} \times \dots D^{2N} \quad ( i_0 \text{times})$$
are the projections.
We will say that $f$ is  {\em based admissible},
if any based asymptotic  periodic solutions $[x_i]$
satisfy lim inf$_i T_i >1$.

Let us construct
based admissible functions by taking infinite products of some
$f: D^{2N} \mapsto [0,1]$.
Let us   choose a smooth function $g: [0,1] \to [0,1]$ such that  for some positive small $\delta, \epsilon >0$:
$$(1) \quad 0 \leq  -  g'  < 1+ \delta  , \quad 
 (2) \quad g |[0, \epsilon] \equiv 1, \quad
 g|[1  - \epsilon ,1 ] \equiv 0.$$
Then we  put $f: D ^{2N} \mapsto [0,1]$ by $f(p) = g(|p|^2)$.
$f$ is admissible by the proof of  lemma $4.2$.
For convenience, let us assign index on  each term: 
 $$\tilde{D} =
D^{2N} \times D^{2N} \times \dots
 = D_1^{2N} \times D_2^{2N} \times \dots \times D_l^{2N} \times \dots$$
We also denote $f_1,f_2, \dots$ where all $f_i$ are the same $f$
and we regard $f_i$ as  functions on $D_i^{2N}$. 
Let us put:
\begin{equation}
\begin{align*}
& F =    f_1 f_2 \dots f_l \dots: \tilde{D}
   \mapsto [0, 1] , \\
& F_i = f_1 \dots f_i: D(i)  = D_1^{2N} \times \dots D_i^{2N} \mapsto  [0,1]
\end{align*}
\end{equation}
where $f_if_j$ imply the pointwise multiplications.
Notice  $F(0) = 1$ at $0 \in \tilde{D}$.
We claim that $F$ is based admissible.
In fact, the Hamilton vector field is:
\begin{equation}
\begin{align*}
X_F & = f_2 f_3 \dots X_{f_1} +  f_1f_3f_4 \dots X_{f_2} + 
  f_1f_2f_4f_5 \dots X_{f_3} + \dots \\
& \equiv  G_1X_{f_1} + G_2 X_{f_2} + \dots
\end{align*}
\end{equation}
$G_i$ take the values in $[0, 1]$.

Let us take a based asymptotic   periodic solution $[l(i)]$ 
and denote its period $T=$ lim inf$_i T(l(i)) > 0$.
Let us fix  $i_0$,  and  denote by $l_i$ as 
 the projection of $l(i)$ on $D_{i_0}^{2N}$ component. 
Then both the convergences:
$$\lim_i | \dot{l_i} - G_{i_0}X_{f_{i_0}}|C^0 = 0 , \quad 
\lim_i | \dot{l(i)} - X_F|C^0 = 0  \quad (*)$$
hold.  One may assume  $\lim \inf_i$ length $l_i  >0$.

The following shows $F$ is based  admissible. 
\begin{sublem}
The periods $T(l_i) $  of $l_i$ are larger than $1$ for large $i$.
\end{sublem}
$Proof:$ 
We will use the properties that each $f_i$ is admissible on $D_i^{2N}$ and $G_i$ takes values less than $1$.
By choosing a subsequence, the family $\{ l_i\}_i$ converges to $l_{\infty}$ in 
$C^0(D_{i_0}^{2N})$. 
Let us put $g_i = G_{i_0}(l(i)(t)): [0,T(l(i))] \mapsto [0,1]$.
Then also a subsequence of $\{ g_i\}_i$ converges to $g: [0,T] \mapsto [0,1]$ 
in $C^0$ by $(*)$.
Thus $\{ l_i \}_i$ converges in fact  in $C^1(D_{i_0}^{2N})$ topology. 

Now the equality  $\dot{l}_{\infty}(t) = g(t)X_{f_{i_0}}$ holds.
Then there is  a small $\delta >0$ so that 
$ g(t) \geq \delta$ hold,
since $l_{\infty}$ is non trivial.
Let us regard  it as the periodic function
$g: {\mathbb R} \to (0, \infty)$.

Let us solve the ODE
 $\alpha: [0, T] \mapsto {\mathbb R}$ satisfying
$\dot{\alpha}(t)= g^{-1}(\alpha(t))$ with $\alpha(0)=0$.
Then put $l' : [0,T]   \mapsto D_{i_0}^{2N}$ by
$l'(t) = l_{\infty}(\alpha(t))$, which 
 satisfies the equation
 $\dot{l}'(t) = g^{-1}(\alpha(t))\dot{l}_{\infty}(\alpha(t))= X_{f_{i_0}}(l'(t))$.

The image of  $l'$ consists of the loop, and
  there is some $0< T'\leq  T$ with $l'(0)= l'(T')$,
since $\dot{\alpha}(t) \geq 1  $ for all $t$.
Since $f_{i_0}$ is admissible, it follows $ T \geq T' >1$.
This completes the proof.
\vspace{3mm} \\
{\bf 4.D.3  Infinite tori:}
Let $ T^{4 \infty} = [  (T^{4n}, \omega, J)]$ be the standard tori.
We verify the following:
 \begin{prop} $C([(T^{4n}, \omega, J)]) = \infty$ 
 \end{prop}
 $Proof:$ 
 Let us consider the standard torus embeddings:
$$T^4 \subset  T^8 \subset T^{12} \subset  \dots \subset 
T^{4n} \subset  \dots$$
Let us choose small real numbers $\alpha_1, \alpha_2 \in {\bf R}$.
Then one has a family of symplectic structures on $T^4$ by:
$$\omega(\alpha_1, \alpha_2) = dx_2 \wedge dx_1 + 
dx_4 \wedge dx_3 + \alpha_1 dx_3 \wedge dx_2 + \alpha_2 dx_1 \wedge dx_3.$$

\begin{sublem}[Z,  see HZ]
Suppose  $(\alpha_1, \alpha_2) \in {\bf R}^2$ is irrational.
  Then the corresponding 
symplectic manifold $(T^4, \omega(\alpha_1, \alpha_2))$
has infinite capacity.
\end{sublem}

{\em Proof of lemma:}
 Let us choose $\alpha_i$ as above,   put the corresponding
form by $\omega(\alpha_i) = \omega(\alpha_1, \alpha_2)$
and 
$\omega = \omega(0,0)$ be  the standard form on $T^4$.
Then we have a  symplectic sequence:
\begin{equation}
\begin{align*}
& T^{4\infty}( \alpha_1, \alpha_2)  =  (T^4, \omega(\alpha_1, \alpha_2)) \times 
  (T^4, \omega) \times  (T^4, \omega) \times \dots  \\
&  T^{4n}( \alpha_1, \alpha_2) = 
(T^4, \omega(\alpha_1, \alpha_2))    \times (T^4, \omega) \times
   \dots  \times (T^4, \omega)  \quad (n \text{ times}).
               \end{align*}
               \end{equation}
 Let $(T^4, \omega ,J)$ be the standard Kaehler structure.
Since 
$\omega(\alpha_1, \alpha_2)( V , J  V ) >0$ take positive values  for  $V \ne 0$
($\alpha_i$ are small),    one can construct another $J'$
which is compatible with $\omega(\alpha_1, \omega_2)$ and
is  near $J$.
Thus $(T^4, \omega(\alpha_i), J')$ gives another almost Kaehler data
which is sufficiently near the standard.

Let us choose  families of pairs $(\alpha_1^l, \alpha_2^l)$ 
and the compatible almost complex structures
 $J^l$ as above so that  $\alpha_i^l \to 0$ and $J^l \to J$
as $l \to \infty$. 
Then we have  a convergent family of product almost  Kaehler sequences 
$\{ [T^{4n},  \omega(\alpha_1^l, \omega_2^l ), J^l)] \}_l$
where  we choose  smooth embeddings 
$T^{4n} \hookrightarrow T^{4n}$
by the identity  for every $n$.
So  the family gives an element:
$$[[T^{4n},  \omega(\alpha_1^l, \omega_2^l ), J^l)] ]
 \in {\frak N}([(T^{4 n}, \omega, J)]).$$
Thus  it is enough to see 
$cap([T^{4n},  \omega(\alpha_1^l, \omega_2^l ), J^l)] ) = \infty$.
By sublemma $4.3$, there are  admissible functions 
$f: (T^4, \omega(\alpha_1^l, \omega_2^l)) \mapsto {\bf R}$
with   sufficiently large $m(f)$.
 $f$ induce the  bounded Hamiltonians on  $[T^{4n},  \omega(\alpha_1^l, \omega_2^l ), J^l)] $,
which are also  admissible over  
$[T^{4n}, \omega(\alpha_1^l, \omega_2^l) , J^l)]$. This implies that the capacity is infinite.
This completes the proof.

\section{Completion of spaces}
In this section we extend the class of the maps we treat to 
 $L^2$ loops over almost Kaehler sequences.
 \vspace{3mm} \\
{\bf 5.A Completion of spaces:}
Let $[(M_i, \omega_i, J_i)]$ be an almost Kaehler sequence, 
and take two $\epsilon$
complete almost Kaehler charts: 
$$\varphi(p_i): D(\epsilon)  \cong D_{p_i} \subset M\cup_{j \geq 0} M_j$$ 
for $i=1,2$ and 
$D(\epsilon) \subset {\mathbb R}^{\infty}$.
Let us  denote the completion of $D(\epsilon)$ by $\bar{D}(\epsilon) \subset H$ where 
$H$ is the  Hilbert space obtained by completion of ${\mathbb R}^{\infty}$ with the standard metric.

 Let us introduce an equivalent relation as follows.
 Two points $z_1, z_2 \in \bar{D}(\epsilon)$ are equivalent, if 
 there is a sequence:
 $$\{m_1,m_2, \dots \} \subset D_{p_1} \cap D_{p_2}\subset M$$
 so that the corresponding two sequences:
 $$w^i_k \equiv \varphi(p_i)^{-1}(m_k) \subset \bar{D}(\epsilon)$$
 converge to $z_1$ and $z_2$ in $\bar{D}(\epsilon)$ respectively.

 \begin{defn}
 The completion of $M=\cup_i M_i$ is defined by:
 $$M \ \subset  \ \bar{M} \ \equiv \ \coprod \ \bar{D}(\epsilon) / \sim .$$

The set of  points at infinity consists of:
$$\partial M \equiv \bar{M} \backslash M.$$
\end{defn}

\begin{lem}
Suppose $[(M_i, \omega_i, J_i)]$ is an almost  Kaehler sequence. Then 
$\bar{M}$ admits the structure of a Hilbert manifold.
\end{lem}
{\em Proof:}
This follows from lemma $1.8$.
This completes the proof.

\vspace{3mm}

Notice that by definition,  any bounded Hamiltonian
 $f: M \to {\bf R}$  extends to a smooth and bounded function $f: \bar{M} \to {\bf R}$.

 Later on we assume  that $\bar{M}$ admit  the Hilbert-manifold structure.
 Let $U_{\epsilon}(M_k)  \subset M$ be the $\epsilon$ neighbourhoods of $M_k$.
 We denote their completions by  $\bar{U}_{\epsilon}(M_k) \subset \bar{M}$.
 Then the holomorphic maps are extended as:
 $$\pi_k: \bar{U}_{\epsilon}(M_k) \to  M_k.$$
{\bf 5.B $L^2$ holomorphic curves:}
Let us introduce the Sobolev spaces of maps into the
completed spaces.
Notice that in order  to  define the full Sobolev spaces,
one has  to use cut off funcions over the local charts,
which play the role of the
locally finite partition of unities in the finite dimensional setting.
In our situation the Sacks-Uhlenbeck's estimates allow us
to obtain  uniform norms   which is necessary to perform 
 analysis of the moduli spaces we study.

Let us take  small $\delta >0$ and choose  $\frac{1}{2} \delta$ net
which is consisted by 
finite set of points $\{x_1, \dots,x_m\} \subset S^2$.
Let $ B_i(\delta) \subset S^2$ be $\delta$ balls with the center $x_i$,
so that $\{B_1(\frac{\delta}{2}), \dots, B_m(\frac{\delta}{2})\}$ 
consists of  an open covering of $S^2$.
Let us choose a smooth function
$\xi_i : B_i (\delta) \to [0,1] $ and have a partition of unity by:
$$\rho_i (x) = \frac{\xi_i(x) }{\Sigma_{i=1}^m \xi_i (x) }  \quad
(x \in S^2), \qquad
\xi_i (x) = 
\begin{cases}
0 & x \in   \partial B_i (\delta) \\
1 &  x \in  B_i(\frac{\delta}{2})
\end{cases}$$

Let 
$\varphi_i: \bar{D}(\epsilon) \cong \bar{D}_{p_i} \subset  \bar{M}$
be complete almost Kaehler charts for 
 finite sets $\{p_1, \dots, p_m\} \subset M$.
Let $u_0: S^2 \to \bar{M}$ be a map with:
$$u_0(B_i(\delta)) \subset \bar{D}_{p_i}$$
such that 
$ \rho_i   u_0  \in L^2_{l+1}(B_i(\delta) : \bar{D}_{\epsilon})$.
Then we obtain the  function spaces around $u_0$ defined by:
$$ L^2_{l+1}(S^2, \bar{M}; \{ p_i\}_i) =
\{ u : S^2 \to \bar{M}: \  u(B_i(\delta)) \subset \bar{D}_{p_i} \ 
  \rho_i   u  \in L^2_{l+1}(B_i(\delta) : \bar{D}_{\epsilon})\}$$
where we equip with the  norms:
$$||u||^2_{L^2_{l+1}} = \Sigma_{i=1}^l ||\rho_i (u)||^2_{L_{l+1}^2}.$$
These spaces admit  the Hilbert manifold structure on small  neighbourhoods of $u_0$.
\vspace{3mm} \\
{\em Remark 5.1:}
One will be able to  obtain the globally defined  function
spaces as Hilbert manifolds 
 by use of the  infinite dimensional version of the Levi-Civit\`a connection ([Kl]).
\vspace{3mm}

With the fixed data $\{x_i\}_i$,
we define the function spaces:
$$ L^2_{l+1}(S^2, \bar{M}) \equiv  \cup_{\{p_i\}_{i=1}^l \subset M}  L^2_{l+1}(S^2, \bar{M}; \{ p_i\}_i).$$

\vspace{3mm}

Let us introduce the functional spaces which are parallel to section $2$.
We define  the spaces of Sobolev maps:
\begin{equation}
\begin{align*}
  {\frak B}  & = \{  u   \in   L^2_{l+1}(S^2, \bar{M}) :  
 [u] = \alpha ,   \\
 & \int_{D(1)} u^*(\omega) = \frac{1}{2}<\omega, \alpha>, \quad
  u(*) = p_* \in M_0  , * \in \{ 0, \infty \}  \}.
  \end{align*}
  \end{equation}

 Let $E(J), \ F \mapsto S^2 \times \bar{M}$ be 
 vector bundles whose fibers are respectively:
\begin{equation} 
\begin{align*}
& E(J)(z,m) = \{ \phi: T_z S^2 \mapsto T_m\bar{M} : 
   \text{ anti complex linear map} \}, \\
& F(z,m) = \{ \phi: T_z S^2 \mapsto T_m\bar{M} : \text{ linear} \}.
\end{align*}
\end{equation}
Then we have   the Hilbert bundle and the  Cauchy-Riemann operators: 
 \begin{equation}
 \begin{align*}
&  {\frak E} = L^2_l({\frak B}^*(E(J))) = 
\cup_{u \in {\frak B}} \{u\} \times L^2_l(u^*(E(J))), \\
&  \bar{\partial}_J \in C^{\infty}({\frak E} \mapsto {\frak B}).
\end{align*}
\end{equation}

\begin{lem}
There is $m$  determined only by $[(M_i, \omega_i, J_i)]$, and
 for any  $l$  there is  $\epsilon =\epsilon(l)  >0$ so that 
 $\epsilon$ neighbourhoods of  any holomorphic curves
 $u: S^2 \to \bar{M}$ with $u(0)=p_0 $ and $u(\infty)=p_{\infty}$,
admit the Hilbert manifold structures in 
    $L^2_{l+1}(S^2, \bar{M}; \{p_i\}_{i=1}^m) $ 
 where $p_i=u(x_i)$.
\end{lem}
{\em Proof:} This follows from 
  lemma $2.2$ since 
the restrictions satisfy $u : B_i(\delta) \to D_{p_i}$ if $\delta >0$
 is sufficiently but uniformly small.
This completes the proof.
\vspace{3mm}

Now let us   define the moduli space of holomorphic curves by:
$$   \bar{\frak M}(\alpha, M,J)  = 
    \{ u \in C^{\infty}(S^2, \bar{M}) \cap {\frak B}    : \bar{\partial}_{J} (u) =0   \}.$$
    A priori inclusions hold:
    $$  {\frak M}(\alpha, M,J) \   \subset  \  \hat{\frak M}(\alpha, M,J)  \  \subset   \ \bar{\frak M}(\alpha, M,J) .$$
    
As before let us   say that $J$ is {\em regular}, if  the linealizations
 $D \bar{\partial}_{J}(u): T_u {\frak B} \mapsto {\frak E}_u$
 are onto  for all $u \in \bar{\frak M}$.

\vspace{3mm}

\begin{prop} Let 
$[(M_i, \omega_i, J_i)]$ be a minimal  and isotropic symmetric  Kaehler sequence.

If the 
moduli space of holomorphic curves is  non empty, regular and
 has $S^1$ freely $1$ dimension, then
 the equality holds:
$$  \bar{\frak M}(\alpha, M,J)  = {\frak M}(\alpha, M,J)  .$$
\end{prop}
{\em Proof:}
It follows from  lemma $2.2$ that for any $u \in  \bar{\frak M}(\alpha, M,J) $,
there is some $k$ so that the image of $u$ lies in $\epsilon$ neighbourhood of $M_k$.

Let $\pi_k: \bar{U}_{\epsilon}(M_k) \to M_k$ be the holomorphic projections.
Then $u_k \equiv \pi_k \circ u: S^2 \to M_k$ are also holomorphic curves.
By proposition $2.2$, the moduli spaces are strongly regular at  $u_k $.
Since  $ {\frak M}(\alpha, M,J)$ is compact by lemma $2.8$,
$u_k$ must be the unique holomorphic curve mod $S^1$ action,
in small neighbourhoods of $u_k$. 
So $u=u_k$ must hold.

This completes the proof.
\vspace{3mm} \\
{\bf 5.C Hamiltonian diffeomorphisms:}
Let $[(M_i, \omega_i, J_i)]$ be an almost Kaehler sequence, and take
a bounded Hamiltonian $f: M= \cup_i M_i \to {\bf R}$.

Let us choose  a uniformly bounded covering
$\{\varphi(p): D(\epsilon) \equiv \cup_i D^{2i}(\epsilon) \hookrightarrow \cup_i M_i\}_{p\in M}$.
By pulling back as $\varphi(p)^*(f):  D(\epsilon) \to {\bf R}$,
let us regard the restriction of the differential $df$ as the one form over $D(\epsilon)$.

The Hamiltonian vector field $X_f$ over $D(\epsilon)$ is defined as the unique
vector field which obey the equality:
$$df(Y) = \omega(X_f, Y)$$
for any  vector field $Y$ of completely bounded geometry over $D(\epsilon)$.
It follows from lemma  $1.8$ that 
$X_f$ is in fact determined globally over $M= \cup_i M_i$.
The following holds by applying the existence and uniqueness of 
ODE with  infinite dimensional targets:
\begin{lem} Let $[(M_i, \omega_i, J_i)]$ be an almost Kaehler sequence, and 
 $X_f$ be the globably defined vector field over $M$ as above.
 Then:
 
 (1)  There is the parametrized diffeomorphisms as its integral:
 $$D_t : \bar{M} \cong \bar{M}.$$
 
 (2) Let $U \subset \bar{M}$ be an open subset, and $f: D_t(U)  \to {\bf R}$ be pre admissible.
 Then the induced function $D_t^*(f) : U \to {\bf R}$ is also the case.
\end{lem}

\vspace{3mm}

We call $D=D_1$ as the {\em Hamiltonian diffeomorphisms}.

\begin{prop} Let $[(M_i, \omega_i, J_i)]$ be an almost Kaehler sequence.
Let us  take a bounded Hamiltonian $f: \cup_i M_i \to {\bf R}$
with the Hamiltonian diffeomorphisms
$D: \bar{M} \cong  \bar{M}$.

Let $U \subset  \bar{M}$ be an open subset. Then we have the equality:
$$\text{As-cap}(U) = \text{As-cap}(D(U)).$$
\end{prop}
{\em Proof:}
Let $f: D(U) \to [0, \infty)$ be an as-admissible function, and consider
the pull back $D^*(f): U \to [0, \infty)$.

Let $X_f$ and $X_{D^*(f)}$ be the Hamiltonian vector fields
over $D(U)$ and $U$ respectively. Then the push forward 
satisfies the equality:
$$D_*(X_{D^*(f)}) (m) = X_f(D(m))$$
for $m \in U$.

Let $x_i : [0, T_i] \to U$ be a non trivial  asymptotic periodic solution to $D^*(f)$.
It follows from the above equality that $D(x_i) : [0, T_i] \to D(U)$ must satisfy:
$$\sup_t |\dot{D(x_i)} - X_f(D(x_i))|(t) \to 0$$
as $i \to \infty$.

Notice that the images of $D(x_i)$ are certianly contained in $\bar{M}$,
but may not in $\cup_i M_i$.
Let $\pi_k: U_{\epsilon}(M_k) \to M_k$ be the family of holomorphic projections 
 in definition $1.3$.

For each $i$, there are some $k_i$ so that the followings hold,
since $\cup_iM_i$ is dense in $\bar{M}$:

(1) the  image of $D(x_i)$ are contained in $U_{\epsilon}(M_{k_i})$,

(2) $\sup_t |\dot{D(x_i)} - \dot{\pi_{k_i}(D(x_i))}| \to 0$ as $i \to \infty$.

\vspace{3mm}

One may assume $k_i <k_{i+1}$ for all $i$.
Let us put another family of non trivial asymptotic  loop:
$$y_l =
\pi_{k_i}(D(x_i)) $$
for  $ k_i   \leq  l \leq k_{i+1} -1$.
Then $y_l : [0, S_l] \to M_l$ is also an asymptotic periodic solution to $f$,
and so $\lim \inf T_i = \lim \inf_l S_l >1$ must hold by as-admissibility of $f$.
So $D^*(f)$ is also admissible.

In particular we have the inequality
$\text{As-cap}(U) \geq \text{As-cap}(D(U))$.
The same argument gives us the converse inequality.
So we must have the equality.
This completes the proof.

\large

\vspace{1cm}

Tsuyoshi Kato,
Department of Mathematics,
Faculty of Science

Kyoto University,
Kyoto 606-8502
Japan.
\vspace{5mm}

\enddocument